\documentclass[english,ruled, 10pt]{article}
\usepackage[T1]{fontenc}
\usepackage[utf8]{inputenc}
\usepackage{amsfonts, amsmath, amssymb, amsthm}
\usepackage{geometry}
\usepackage{algorithm, algpseudocode}
\usepackage{nicematrix}
\usepackage{bm}
\usepackage{xcolor}
\usepackage{optidef}
\usepackage{pstricks, pst-tree}
\usepackage{graphicx}
\usepackage[square,numbers]{natbib}
\usepackage{mathabx}
\usepackage{mathtools}
\usepackage{listings}
\usepackage{float}
\usepackage{caption}
\usepackage{pifont}
\usepackage{babel}
\usepackage{todo}
\usepackage{hyperref}
\hypersetup{
	colorlinks=true,
	linkcolor=red,
	citecolor = blue
}

\theoremstyle{plain}
\newtheorem*{assumption*}{\protect\assumptionname}

\theoremstyle{remark}
\newtheorem*{remark*}{\protect\remarkname}

\theoremstyle{plain}
\newtheorem{remark}{\protect\remarkname}

\theoremstyle{plain}
\newtheorem{thm}{\protect\theoremname}

\theoremstyle{definition}
\newtheorem{Defi}{\protect\definitionname}

\theoremstyle{plain}
\newtheorem{Assumption}{\protect\assumptionname}

\theoremstyle{plain}
\newtheorem{prop}{\protect\propositionname}

\theoremstyle{plain}
\newtheorem{lem}{\protect\lemmaname}

\theoremstyle{plain}
\newtheorem{coro}{\protect\corollaryname}

\newenvironment{Proof}{{\noindent\bf Proof.}}{\hfill$\blacksquare$\\[1mm]}

\providecommand{\assumptionname}{Assumption}
\providecommand{\corollaryname}{Corollary}
\providecommand{\definitionname}{Definition}
\providecommand{\lemmaname}{Lemma}
\providecommand{\propositionname}{Proposition}
\providecommand{\remarkname}{Remark}
\providecommand{\theoremname}{Theorem}

\numberwithin{equation}{section}

\newenvironment{keywords}{
  \par\noindent
  \textbf{Keywords:}
}{
  \par
}

\newif\ifsiam

\siamfalse 

\title{First-order Methods for Stochastic Variational Inequality Problems with Function Constraints
}
\author{
Digvijay Boob\thanks{dboob@smu.edu, Operations Research and Engineering Management, Southern Methodist University}
\hspace{6em}
Qi Deng \thanks{qideng@sufe.edu.cn, School of Information Management and Engineering, Shanghai University of Finance and Economics}
\hspace{6em}
Mohammad Khalafi\thanks{mohamadk@smu.edu, Operations Research and Engineering Management, Southern Methodist University}
}
\date{}

\begin{document}
\allowdisplaybreaks
\global\long\def\vertiii#1{\left\vert \kern-0.25ex  \left\vert \kern-0.25ex  \left\vert #1\right\vert \kern-0.25ex  \right\vert \kern-0.25ex  \right\vert }%
\global\long\def\matr#1{\bm{#1}}%
\global\long\def\til#1{\tilde{#1}}%
\global\long\def\wt#1{\widetilde{#1}}%
\global\long\def\wh#1{\widehat{#1}}%
\global\long\def\wb#1{\widebar{#1}}%
\global\long\def\mcal#1{\mathcal{#1}}%
\global\long\def\mbb#1{\mathbb{#1}}%
\global\long\def\mtt#1{\mathtt{#1}}%
\global\long\def\ttt#1{\texttt{#1}}%
\global\long\def\inner#1#2{\langle#1,#2\rangle}%
\global\long\def\binner#1#2{\big\langle#1,#2\big\rangle}%
\global\long\def\Binner#1#2{\Big\langle#1,#2\Big\rangle}%
\global\long\def\br#1{\left(#1\right)}%
\global\long\def\bignorm#1{\bigl\Vert#1\bigr\Vert}%
\global\long\def\Bignorm#1{\Bigl\Vert#1\Bigr\Vert}%
\global\long\def\setnorm#1{\Vert#1\Vert_{-}}%
\global\long\def\rmn#1#2{\mathbb{R}^{#1\times#2}}%
\global\long\def\deri#1#2{\frac{d#1}{d#2}}%
\global\long\def\pderi#1#2{\frac{\partial#1}{\partial#2}}%
\global\long\def\onebf{\mathbf{1}}%
\global\long\def\zero{\mathbf{0}}%

\global\long\def\norm#1{\lVert#1\rVert}%
\global\long\def\bnorm#1{\big\Vert#1\big\Vert}%
\global\long\def\Bnorm#1{\Big\Vert#1\Big\Vert}%

\global\long\def\brbra#1{\big(#1\big)}%
\global\long\def\Brbra#1{\Big(#1\Big)}%
\global\long\def\rbra#1{(#1)}%
\global\long\def\sbra#1{[#1]}%
\global\long\def\bsbra#1{\big[#1\big]}%
\global\long\def\Bsbra#1{\Big[#1\Big]}%
\global\long\def\cbra#1{\{#1\}}%
\global\long\def\bcbra#1{\big\{#1\big\}}%
\global\long\def\Bcbra#1{\Big\{#1\Big\}}%

\global\long\def\grad{\nabla}%
\global\long\def\Expe{\mathbb{E}}%
\global\long\def\rank{\text{rank}}%
\global\long\def\range{\text{range}}%
\global\long\def\diam{\text{diam}}%
\global\long\def\epi{\text{epi }}%
\global\long\def\inte{\operatornamewithlimits{int}}%
\global\long\def\cov{\text{Cov}}%
\global\long\def\argmin{\operatornamewithlimits{argmin}}%
\global\long\def\argmax{\operatornamewithlimits{argmax}}%
\global\long\def\tr{\operatornamewithlimits{tr}}%
\global\long\def\dis{\operatornamewithlimits{dist}}%
\global\long\def\sign{\operatornamewithlimits{sign}}%
\global\long\def\prob{\mathbb{P}}%
\global\long\def\st{\operatornamewithlimits{s.t.}}%
\global\long\def\dom{\text{dom}}%
\global\long\def\diag{\text{diag}}%
\global\long\def\and{\text{and}}%
\global\long\def\st{\text{s.t.}}%
\global\long\def\Var{\operatornamewithlimits{Var}}%
\global\long\def\raw{\rightarrow}%
\global\long\def\law{\leftarrow}%
\global\long\def\Raw{\Rightarrow}%
\global\long\def\Law{\Leftarrow}%
\global\long\def\vep{\varepsilon}%
\global\long\def\dom{\operatornamewithlimits{dom}}%

\global\long\def\Lbf{\mathbf{L}}%

\global\long\def\Ffrak{\mathfrak{F}}%
\global\long\def\Gfrak{\mathfrak{G}}%
\global\long\def\vfrak{\mathfrak{v}}%
\global\long\def\LGfrak{{\lambda\mathfrak{G}}}%
\global\long\def\gfrak{\mathfrak{g}}%
\global\long\def\sfrak{\mathfrak{s}}%
\global\long\def\xibar{\bar{\xi}}%
\global\long\def\Cbb{\mathbb{C}}%
\global\long\def\Ebb{\mathbb{E}}%
\global\long\def\Fbb{\mathbb{F}}%
\global\long\def\Nbb{\mathbb{N}}%
\global\long\def\Rbb{\mathbb{R}}%
\global\long\def\extR{\widebar{\mathbb{R}}}%
\global\long\def\Pbb{\mathbb{P}}%
\global\long\def\Acal{\mathcal{A}}%
\global\long\def\Bcal{\mathcal{B}}%
\global\long\def\Ccal{\mathcal{C}}%
\global\long\def\Dcal{\mathcal{D}}%
\global\long\def\Fcal{\mathcal{F}}%
\global\long\def\Gcal{\mathcal{G}}%
\global\long\def\Hcal{\mathcal{H}}%
\global\long\def\Ical{\mathcal{I}}%
\global\long\def\Kcal{\mathcal{K}}%
\global\long\def\Lcal{\mathcal{L}}%
\global\long\def\Mcal{\mathcal{M}}%
\global\long\def\Ncal{\mathcal{N}}%
\global\long\def\Ocal{\mathcal{O}}%
\global\long\def\Pcal{\mathcal{P}}%
\global\long\def\Ucal{\mathcal{U}}%
\global\long\def\Scal{\mathcal{S}}%
\global\long\def\Tcal{\mathcal{T}}%
\global\long\def\Xcal{\mathcal{X}}%
\global\long\def\Ycal{\mathcal{Y}}%
\global\long\def\Vcal{\mathcal{V}}%
\global\long\def\Ubf{\mathbf{U}}%
\global\long\def\Pbf{\mathbf{P}}%
\global\long\def\Ibf{\mathbf{I}}%
\global\long\def\Ebf{\mathbf{E}}%
\global\long\def\Abs{\boldsymbol{A}}%
\global\long\def\Qbs{\boldsymbol{Q}}%
\global\long\def\Lbs{\boldsymbol{L}}%
\global\long\def\Pbs{\boldsymbol{P}}%
\global\long\def\i{i}%
\global\long\def\Ibb{\mathbb{I}}

\DeclarePairedDelimiterX{\inprod}[2]{\langle}{\rangle}{#1, #2}
\DeclarePairedDelimiter\abs{\lvert}{\rvert}
\DeclarePairedDelimiter{\bracket}{ [ }{ ] }
\DeclarePairedDelimiter{\paran}{(}{)}
\DeclarePairedDelimiter{\braces}{\lbrace}{\rbrace}
\DeclarePairedDelimiterX{\gnorm}[3]{\lVert}{\rVert_{#2}^{#3}}{#1}
\DeclarePairedDelimiter{\floor}{\lfloor}{\rfloor}
\DeclarePairedDelimiter{\ceil}{\lceil}{\rceil}

\global\long\def\tsum{{\textstyle {\sum}}}

\newcommand{\opex}{\texttt{OpEx}}
\newcommand{\opconex}{\texttt{OpConEx}}
\newcommand{\sopconex}{\texttt{P-OpEx}}
\newcommand{\stochsopconex}{\texttt{S-P-OpEx}}
\newcommand{\fstochsopconex}{\texttt{FS-P-OpEx}}
\newcommand{\adopconex}{\texttt{AdOpEx}}
\newcommand{\gconex}{\texttt{GradConEx}}
\newcommand{\stopconex}{\texttt{StOpConEx}}
\newcommand{\fstopconex}{\texttt{F-StOpConEx}}

\maketitle

\begin{abstract}
The monotone Variational Inequality (VI) is a general model with important applications in various engineering and scientific domains. In numerous instances, the VI problems are accompanied by function constraints that can be data-driven, making the usual projection operator challenging to compute. 
This paper presents novel first-order methods for the function-constrained Variational Inequality (FCVI) problem in smooth or nonsmooth settings with possibly stochastic operators and constraints. Computing the projection operator is challenging for FCVI. We introduce the \adopconex{} method, which employs an operator extrapolation on the KKT operator of the FCVI in a smooth deterministic setting. Since this operator is not uniformly Lipschitz continuous in the Lagrange multipliers, we employ an adaptive two-timescale algorithm leading to bounded multipliers and achieving the optimal $O({1}/{T})$ convergence rate. For the nonsmooth and stochastic VIs, we introduce design changes to the \adopconex{} method and propose a novel \sopconex{} method that takes a partial extrapolation. It converges at the rate of $O({1}/{\sqrt{T}})$ when both the operator and constraints are stochastic or nonsmooth. This method has suboptimal dependence on the noise and Lipschitz constants of function constraints. We propose a constraint extrapolation approach leading to the \opconex{} method that improves this dependence by an order of magnitude.  
All our algorithms easily extend to saddle point problems with function constraints that couple the primal and dual variables while maintaining the same complexity results. To the best of our knowledge, all our complexity results are new in the literature.
 \vspace{-2mm}
\end{abstract}

\vspace{-1mm}
\section{Introduction}
In this paper, we consider the following monotone variational inequality problem
\begin{equation}\label{eq:FCVI}
	\text{Find } x^* \in \wt{X} : \quad \inprod{F(x^*)}{x -x^*} \ge 0, \quad\forall x \in \wt{X}, 
\end{equation}
where $\wt{X}\subset \Rbb^n $ is a closed convex set and $F$ is a {\em monotone} operator. We denote this problem as VI$(F, \wt{X})$. The variational inequality problem formed in~\eqref{eq:FCVI} has a variety of applications in finance, economics, transportation, communication, and many others. We refer to \cite{facchinei2003finite} for a review of the classic applications. There is also a growing interest in solving constrained min-max optimization, an important class of VI, that includes robust optimization~\citep{ben2002robust}, generative adversarial networks~\citep{goodfellow2020generative}, and reinforcement learning~\citep{bertsekas2011approximate}.  

The majority of works on VI problems assume that projection onto the set $\wt{X}$ is easy. However,  this paper considers a more general setting where the projection requirement onto set $\wt{X}$ is relaxed.
Specifically, the domain set $\wt{X}$ is formed by the intersection of two sets:
 \begin{equation}\label{eq:def_X_tilde}
 \wt{X}:= X \cap\{x: g_j(x) \le 0, j = 1, \dots, m\},    
 \end{equation}
 where $X \subset \Rbb^n$ is a ``simple'' convex compact set that allows for easy projection, and $g_j: X \rightarrow \Rbb$ is  a continuous convex function for all $j = 1, \dots, m$. Since projecting onto function the constraint set $\{x: g_j(x)\le 0,\ j=1, \dots, m\}$ is not possible, this setting leads to a \emph{function constrained VI} (FCVI) problem. We denote this problem as FCVI$(F, X, g_j(\cdot), j = 1, \dots, m\} = $ VI$(F, X \cap \{x:g_j(x) \le 0, j =1, \dots, m\})$. Throughout this paper, we assume that the operator $F:X\rightarrow \Rbb^n$ is a monotone operator, i.e., 
\[\inprod{F(x_1)-F(x_2)}{x_1 -x_2} \ge 0, \quad\forall x_1, x_2 \in X.\]  

We say that an FCVI problem is smooth if the operator $F(x)$ and gradients of constraint $g_j(x), j=1, \dots, m,$ are Lipschitz continuous. If any of these assumptions are not satisfied then we call it a nonsmooth FCVI problem. In this case, we assume that $F(x)$ is a bounded operator and $g_j(x), j =1, \dots, m$ are nonsmooth functions with bounded subgradients. We consider both smooth and nonsmooth FCVI problems in this paper.
Furthermore, we examine the stochastic variational inequalities (SVIs) problem, where  $F(x)$ is the expectation of a stochastic function:  $F(x) = \Ebb[\mathfrak{F}(x,\xi)]$. In this context, we assume the operator $F(x)$ is unavailable and can only be approximated through a stochastic estimator $\mathfrak{F}(x,\xi)$, with $\xi$ representing a random variable input. To allow the possibility of data-driven constraints in the FCVI problem, we consider a fully stochastic setting where both  $F(x)$ and  $g(x)$ are the expectation of stochastic functions $\Ffrak(x, \xi)$ and $\gfrak(x, \xi)$, respectively. Problems with general stochastic operators and function constraints have several important applications in science and engineering (See \cite{juditsky2011solving, Kotsalis2022, yang2024data}). In the next section, we review existing literature on VIs, SVIs, and some very recent developments in FCVIs.

\subsection{Literature review}
The classical study on VI problems focused on the structure of the problem, and quasi Newton-type and interior point methods that ensure a quadratic convergence in the local neighborhood of {\em regular} solutions \citep{pang1982iterative, josephy1979quasi, subramanian1985gauss, harker1990finite, ralph2000superlinear, qi2002smoothing}. The VI approach is also used in various engineering applications \cite{pang2010design,shanbhag2013stochastic}. See \cite{facchinei2003finite} for a review of these topics.
Owing to the growing interest in variational inequalities and the need for solving large-scale problems, a substantial body of literature has focused on developing first-order algorithms with convergence rate guarantees for solving VIs. We focus our literature review on this line of work. Throughout this discussion, we denote by $T$ the number of operator evaluations. For the deterministic and smooth VI problem, \cite{Nemirovski2004} proposed the Mirror-Prox method, which for the first time obtains a nontrivial convergence rate of  $\mathcal{O}(1/T)$. This rate significantly improves upon the $\mathcal{O}(1/\sqrt{T})$ rate of the standard projected gradient method. Since this method computes only two operator evaluations in all iterations, it leads to the optimal order of the oracle complexity due to lower bounds on bilinear saddle point problems \citep{ouyang2021lower}. Later, \cite{Nesterov2007} proposed a more simplified dual extrapolation method that obtains the same optimal convergence guarantee. In \cite{nesterov2011solving}, this rate is further improved to $\mathcal{O}(e^{-T})$ for strongly monotone VIs. 
This progress has spurred additional research efforts \citep{mokhtari2020convergence, GidelBVVL19} on the efficiency of classic methods such as the extragradient method~\citep{Korpelevich1976} and optimistic gradient~\citep{popov1980modification}. More simplified algorithms, like the projection-reflected gradient method~\citep{malitsky2015projected} and operator extrapolation~\citep{Kotsalis2022}, maintain a single sequence and only require one projection in each iteration. Some recent work focuses on auto-tuning of step-lengths \cite{malitsky2020golden, diakonikolas2020halpern}, and convergence in terms of strong VI solutions \citep{dang2015convergence, diakonikolas2020halpern}. However, despite the advances in first-order methods for VIs, the vast majority of the prior literature deals with VI$(F,\wt{X})$ problem where either $\wt{X} = \Rbb^n$ is an unconstrained domain or $\wt{X}$ is a ``simple'' set such that one can easily evaluate exact projection onto $X$.

The VI problem becomes significantly more challenging in the stochastic setting, where the exact evaluation of the operator $F(x)$ is not possible. In this case, we only have access to the noisy (unbiased) estimator $F(x, \xi)$ through a {\em stochastic oracle} where $\xi$ is the stochastic input. We focus the majority of our literature review on the line of work that establishes the sample complexity guarantee for stochastic VIs. Throughout the literature review of stochastic VI, we denote by $T$ the number of stochastic operator $F(x, \xi)$ evaluations. 
\cite{jiang2008stochastic} proposed solving SVIs using the well-known stochastic approximation (SA) method~\cite{robbins1951stochastic} and established asymptotic convergence to the global solution when the map is strongly monotone and smooth. Extensions of SA to merely monotone and non-Lipschitz SVI have been made in a series of follow-up works~\citep{koshal2012regularized, yousefian2017smoothing}.
Motivated by the success of mirror-prox (and other extragradient-type) methods in the deterministic setting,
\cite{juditsky2011solving} extended the mirror-prox method~\cite{Nemirovski2004}, which appears to be the first work that achieves a convergence rate of $\mathcal{O}(1/\sqrt{T})$ for SVIs. In subsequent work, \cite{yousefian2014optimal} developed the stochastic extragradient method for monotone VIs and obtained an improved $\mathcal{O}(1/T)$ rate under the additional weak sharpness assumption. Recently, \cite{boob2021optimal} developed optimal algorithms for the SVI problems under the constraint of differential privacy.
\cite{Kotsalis2022} developed the stochastic operator extrapolation method and demonstrated its ability to achieve optimal rates for a variety of SVIs, particularly for stochastic smooth and strongly monotone VIs.
\cite{huang2022new} proposed novel extra-point and momentum methods that achieve linear convergence to a small neighborhood of the optimal solution for stochastic smooth and strongly monotone VIs. They also adapted their techniques to develop zeroth-order methods for stochastic min-max saddle point problems. Despite significant progress, the existing literature for VI problem heavily relies on deterministic set constraint with access to the exact projection oracles. This assumption is restrictive if constraints have a functional form in \eqref{eq:def_X_tilde}. In particular, if $g$ is an expectation function, then there is no guarantee of even an approximate projection as we do not have access to $g$ itself.

There is one recent effort in the literature to develop first-order methods for FCVI$(F,X,g_j, j=1, \dots, m)$. 
\cite{Yang2022Solving} proposed a first-order interior point method called ACI for solving convex function-constrained monotone VIs. The key idea of ACI is to use the Alternating Direction Method of Multipliers to solve the penalty function subproblem within the interior point method. For a small enough barrier penalty, the authors demonstrated that the gap function over the last iterate converges at an $\mathcal{O}(1/\sqrt{T})$ rate, while over the ergodic mean it converges at an $\mathcal{O}(1/T)$ rate. Despite ACI's advantages over traditional interior point methods, it still requires a linear constraint projection, necessitating the iterative solution of linear systems of equations, thus making it difficult to generalize to large-scale problems. Moreover, ACI is a multi-loop algorithm that requires the exact solution of the penalty function subproblem, which is difficult to compute for general non-separable function constraints. Due to the assumptions and complicated multi-loop design, extending this method for the stochastic and fully-stochastic FCVI problem is impractical. 

From the point of view of saddle point problems with function constraints 
,  \cite{yang2024data} developed a new primal-dual algorithm for stochastic convex-concave saddle point problems with expectation function constraints and obtained the convergence rate of $\mathcal{O}(1/\sqrt{T})$. Despite the generality of this setting, 
their approach necessitates that the function constraints be separable for both the min and max variables. We will show that our framework of FCVI leads to a more general framework of saddle point problems and our proposed algorithms readily extend to this more general saddle point framework with either matching or better convergence rate guarantees.

\subsection{Contributions} 
In this paper, we propose several new first-order methods for various VIs problems with convex function constraints. Our contribution can be summarized as follows. 

First, we write the FCVI problem as a KKT system where $x \in X$ is the primal variable and $\lambda_j (\ge 0)$ denotes the dual (or KKT) multipliers associated with constraints $g_j(x)\le 0, j =1, \dots, m$. Denoting $\lambda = [\lambda_1,\dots, \lambda_m]^T$, we observe that the KKT system can be rewritten as VI$(\wt{F}(x, \lambda), X \times \{\lambda \ge \zero\})$ where $\wt{F}(x, \lambda)$ is jointly monotone in $(x, \lambda)$ over the feasible set is $X \times \{\lambda \ge \zero\}$. This motivates a recently proposed Operator Extrapolation type method \cite{Kotsalis2022} for solving FCVI$(F, X, g_j(\cdot), j = 1, \dots, m)$ applied on the joint operator $\wt{F}$ and feasible set $X \times \lambda$. However, this method is shown to converge only for Lipschitz continuous operators. This assumption is violated even for smooth FCVI problems since the dual set $\{\lambda \ge \zero\}$ is unbounded. To overcome this challenge, we propose an Adaptive Operator Extrapolation (\adopconex) method. This method has a two-timescale design for updates of $x$ and $\lambda$ that changes adaptively. It can ensure that $\lambda$ remains bounded for the smooth deterministic FCVI problems. After running for $T$ iterations, we show that \adopconex~method achieves the convergence rate of  $\mathcal{O}(1/T)$. Here, convergence for a candidate solution $\wh{x} \in X$ is measured by the well-known VI gap $\max_{x \in \wt{X}}\inprod{F(x)}{\wh{x} - x}$, and total infeasibility $\gnorm{\max\{g(\wh{x}),\zero\}}{}{}$. When we say the convergence rate is $\mathcal{O}(1/T)$, we imply that both the VI gap and total infeasibility are bounded above by $\mathcal{O}(1/T)$. Since \adopconex{} computes single new evaluation of $F$, $g$ and $\grad g$, it exhibits the optimal order of oracle complexity due to lower bounds of \cite{ouyang2021lower}.

Second, we observe that even if \adopconex~method achieves the optimal convergence rate in terms of $T$, it is only applicable for smooth deterministic FCVI problems. The same ideas cannot be readily applied for nonsmooth or stochastic FCVI problems as 
the dual update gets unstable and cannot be shown to remain bounded due to increased, possibly stochastic, errors during the update. To get around this issue, we propose a Partial Operator Extrapolation (\sopconex{})~method, which only extrapolates the primal (or the $F$-operator) part and foregoes any extrapolation for the dual update. We show that \sopconex{}~method converges for nonsmooth FCVI problem at the rate of $\mathcal{O}(\tfrac{1}{\sqrt{T}})$. We utilize a conjugate interpretation of $g$ in our analysis, though it is not required to know any properties of the conjugate of $g$ in the implementation of the \sopconex{} method. Subsequently, we also extend this method for stochastic and fully stochastic FCVI problems. In the stochastic FCVI, we only assume that the operator $F$ is given in the expectation and $g$ is a deterministic map. In the fully stochastic case, we assume both $F$ and $g$ are in the expectation form. In this setup, we cannot access $F$, $g$, or its subgradients. We only allow access to them through unbiased stochastic oracles of bounded variance. We also need to assume that $g(\cdot, \xi)$ is a convex function for all stochastic input $\xi$. For this setting, we achieve expected convergence at the rate of $\mathcal{O}(\tfrac{1}{\sqrt{T}})$ after getting $T
$ samples of $F$, $g$ and its subgradient through the stochastic oracle.

Third, we see that for the nonsmooth FCVI problem, the dependence of convergence rate of \sopconex{} method in terms of some Lipschitz constants and noise of $g$ is suboptimal. Essentially, the lack of extrapolation (or momentum) in the dual update of \sopconex{} leads to a large impact of the Lipschitz constant of $g$ and the noise of its stochastic oracle on the convergence rate. This results in the convergence rate of $\mathcal{O}(\tfrac{M_g + \sigma_g^2}{\sqrt{T}})$ where $M_g$ is Lipschitz constant of continuity of $g$ and $\sigma_g$ is the standard deviation of the stochastic oracle of $g$. On the other hand, the extrapolation term was removed to make the dual update stable, which ensures convergence. Recently, \cite{boob2023stochastic} showed a Constraint Extrapolation method that linearizes $g$ while extrapolating in the context of convex optimization problems. They show that it leads to a stable update for the dual. Motivated by this design, we propose an Operator Constraint Extrapolation (\opconex{}) method that 
employs linear approximations of $g$ to define the dual extrapolation (or acceleration/momentum) step. 
We show that \opconex{} obtains the optimal convergence rate of $\mathcal{O}(\tfrac{1}{\sqrt{T}})$ for nonsmooth FCVIs. Moreover, the dependence of convergence rate in terms of $M_g$ and $\sigma_g$ is of $\mathcal{O}(\tfrac{M_g}{T} + \tfrac{\sigma_g}{\sqrt{T}})$. We also show that the impact of other constants can be substantially reduced if a constant $B$ satisfying $B \ge \gnorm{\lambda^*}{}{}$ can be estimated where the $\lambda^*$ is the KKT multiplier (also called as the optimal dual). If this is not possible, the algorithm still converges and maintains $\mathcal{O}(\tfrac{M_g}{T} + \tfrac{\sigma_g}{\sqrt{T}})$ convergence though the impact of $\gnorm{\lambda^*}{}{}$ will be worse on the convergence of the feasibility criterion. We also show that 
Furthermore,\sopconex{} method $g(\cdot, \xi)$ is convex function for any stochastic input $\xi$. This assumption is not required for the convergence of \opconex{} method. It is important to note that \opconex{} method converges in expectation for stochastic FCVI problems 
while maintaining $\mathcal{O}(\tfrac{1}{T})$ convergence for smooth deterministic components and hence, obtains unified convergence guarantee that is best-known in the literature for the composite fully-stochastic FCVI problem under the mildest assumptions on the stochastic oracle. 

Fourth, as an important extension (and a result in its own right) of function-constrained VIs, we examine the convex-concave saddle point problem with a coupling constraint on both min and max players
\[ \min_{x \in X}\max_{y \in Y} f(x,y), \quad\text{s.t. }g(x, y)\le 0.\]
This seems to be the most general form of saddle point problem in the literature as it involves not only function constraints but also the coupling of choices for the min and max players.  
We demonstrate that all of our algorithms developed for FCVI problems apply to this coupled saddle point setting, which encompasses smooth, nonsmooth, or stochastic $f$ and smooth, nonsmooth, or stochastic $g(x)$. To implement these algorithms, we write $F(x,y):= [f_x'(x,y)^T, - f'_y(x,y)^T]^T$ and show that convergence in terms of the primal-dual gap is of the same order as that of VI. Here, $f'_x(x,y)$ and $ -f'_y(x,y)$ are subgradients of convex functions $f(\cdot,y)$ and $ -f(x, \cdot)$ at point $x$ and $y$, respectively. This seems to be the first convergence rate analysis for coupled convex-concave saddle point problems that nontrivially extends recent works, e.g., \cite{yang2024data}, which looked at stochastic function-constrained saddle point problems not involving coupled constraints. In fact, our convergence rate for \opconex{} method is tighter for this more general setting than those of \cite{yang2024data} as our algorithm effectively differentiates between the smooth and nonsmooth/stochastic components.


\begin{table}[t]
    \centering
    \caption{Convergence rate of the proposed methods for solving different FCVIs}\label{tab:conve_rates}
    \begin{tabular}{ccccc}
         Algorithm &Deterministic smooth & Nonsmooth &Stochastic &Fully stochastic \\\hline\hline
         \adopconex& $\mathcal{O}(\tfrac{1}{T})$ &-- &-- &--\\
         \sopconex &-- & $\mathcal{O}(\tfrac{M_g}{\sqrt{T}})$ & $\mathcal{O}(\tfrac{M_g + \sigma}{\sqrt{T}})$ & $\mathcal{O}(\tfrac{M_g+\sigma}{\sqrt{T}} + \tfrac{\sigma_g^2}{\sqrt{T}})$ \\
         \opconex{} &--&$\mathcal{O}(\tfrac{M_g}{T} + \tfrac{H_g}{\sqrt{T}})$ &$\mathcal{O}(\tfrac{M_g}{T} + \tfrac{H_g + \sigma}{\sqrt{T}})$ & $\mathcal{O}(\tfrac{M_g}{T} + \tfrac{H_g+\sigma+ \sigma_g}{\sqrt{T}})$\\\hline
    \end{tabular}
    \caption*{\footnotesize $M_g$: Lipschitz constant of $g$; $H_g$: Lipschitz constant of only the nonsmooth component of $g$, $\sigma$: Standard deviation of stochastic oracle for $F$, $\sigma_g$: standard deviation of stochastic oracles associated with $g$ and $\grad g$.}
\end{table}
\paragraph{Organization of the paper}
Section~\ref{sec:preliminaries} elucidates the notations and terminology pertinent to VIs. Section~\ref{sec:adaptive FCVI} introduces the \adopconex{}~method for solving the smooth deterministic FCVI problem~\eqref{eq:FCVI}. In Section~\ref{sec:Deterministic-simple-PD}, we propose the \sopconex{} method and show its convergence for nonsmooth 
and stochastic 
FCVIs. Section~\ref{sec:dstoch_opconex} proposes the \opconex{}~method for nonsmooth or stochastic FCVIs that significantly improves over the convergence rate guarantees of \sopconex{} method.
In Section~\ref{sec:saddle}, we exhibit our methods for convex-concave saddle point problems with coupling constraints and show the first complexity analysis for such problems. We wrap up with conclusions in Section \ref{sec:conclusion}.

\subsection{Notation and Preliminaries}\label{sec:preliminaries}
We use the following notation throughout the paper. Let $\gnorm{\cdot}{}{}$ represent the standard Euclidean norm. The nonnegative orthant of $\Rbb^k$ is denoted as $\Rbb^k_+$. For the given convex set $X$, we denote its normal cone at $x\in X$ as $N_X(x):=\{v:\inprod{v}{\wt{x}-x}\le 0, \ \forall\ \wt{x} \in X\}$, interior as $\text{int }{X}$ and relative interior as $\text{rint }X$. For the compact set $X$, we denote its diameter $D_X:=\max_{x_1, x_2 \in X}\gnorm{x_1-x_2}{}{}$. Let $[m]:=\{1, \dots, m\}$, $g(x) := [g_1(x), \dots, g_m(x)]^T$ and constraints in \eqref{eq:def_X_tilde} can be expressed as $g(x) \le \zero$. Here, bold $\zero$ denotes the vector with all elements being $0$. 
For any $x\in \Rbb$, we define $[x]_+$ as $\max\{x,0\}$. For a vector $x \in \Rbb^k$, we use $[x]_+$ to denote element-wise application of $[\cdot]_+$. Moreover, the ``+'' operation on sets denotes their Minkowski sum.

We denote the gradient of any convex function $h(x)$ by $\grad h(x)$. Moreover, we denote its subgradient by  $h'(x)$. We represent its subdifferential $\partial h$, which is defined as follows: at a point $x$ in the rint $X$, $\partial h(x)$ is comprised of all subgradients $h'(x)$ which are in the linear span of $X-X$. For a point $x \in X \setminus \text{rint }X$, $\partial h(x)$ consists of all vectors $v$, if any, such that there exists $x_i \in \text{rint }X$ and $v_i \in \partial h(x_i), i =1,2, \dots,$ with $x = \lim\limits_{i \to \infty}x_i$ and $v = \lim\limits_{i \to \infty}v_i$. Following this definition, it is well-known that if a convex function $h:X \to \Rbb$ is Lipschitz continuous with constant $M_h$, with respect to norm $\gnorm{\cdot}{}{}$, then the set $\partial h(x)$ is nonempty for any $x\in X$ and 
\[ h'(x) \in \partial h(x) \Rightarrow \abs{\inprod{h'(x)}{d}} \le M_h\gnorm{d}{}{}, \forall d \in \text{lin}(X-X), \]
which also implies 
\[\gnorm{h'(x)}{}{} = \sup_{\gnorm{z}{}{} \le 1} \inprod{ h'(x)}{z} \le M_h.\]
See \cite{ben2005non} for more details.
Using this definition of subgradients, we define $g'(x) := [g'_1(x), \dots, g'_m(x)] \in \Rbb^{n \times m}$ as the Jacobian of the constraint mapping $g(x)$. Similarly, when $g_j,  j  \in [m]$ are smooth functions, we define $\grad g(x) := [\grad g_1(x), \dots ,\grad g_m(x)] \in \Rbb^{n\times m}$ as the Jacobian of the differentiable constraint mapping. We also denote $\partial g(x) := \partial g_1(x) \times \dots \times \partial g_m(x)$.

We provide formal definitions for several concepts introduced earlier.  
We consider the general case where $F$ is possibly composed of Lipschitz continuous and bounded discontinuous monotone operators, satisfying
\begin{equation}\label{eq:F-Lipschitz-property}
	\gnorm{F(x_1) - F(x_2)}{}{} \le L\gnorm{x_1-x_2}{}{} + H, \quad \forall x_1, x_2\in X,
\end{equation}
where the $L$-term arises from the Lipschitz continuous component and the $H$-term from the discontinuous component. When $H = 0$, the above definition reduces to a Lipschitz continuous operator, also known as a ``smooth operator'' in some VI literature.
We assume that $g_j$ is comprised of smooth and nonsmooth components satisfying the following Lipschitz-type property:
\begin{equation}
    \label{eq:g_i-Lipschitz-property}
    g_j(x_1) -g_j(x_2) - \inprod{ g'_j(x_2)}{x_2-x_1} \le \tfrac{L_{g_j}}{2} \gnorm{x_1-x_2}{}{2} + H_{g_j}\gnorm{x_1-x_2}{}{}, \quad \forall x_1, x_2\in X, j \in [m],
\end{equation}
where $L_{g_j}$-term arises from the Lipschitz smooth component of $g_j$ and $H_{g_j}$-term arises from the nonsmooth component of $g_j$. 
Using \eqref{eq:g_i-Lipschitz-property}, Cauchy-Schwarz inequality and convexity of $g_j, j=1, \dots,m$, it is easy to see that
\begin{equation}
    \label{eq:g-Lipschitz-prop}
    \gnorm{g(x_1) - g(x_2) - \inprod{ g'(x_2)}{x_1-x_2}}{}{} 
    \le \tfrac{L_g}{2}\gnorm{x_1-x_2}{}{2} + H_g\gnorm{x_1-x_2}{}{},
\end{equation}
where 
$L_g := \paran{\tsum_{j=1}^m L_{g_j}^2}^{1/2}$ and $H_g := \paran{\tsum_{j=1}^m H_{g_j}^2}^{1/2}$. Moreover, we make the additional assumption that constraint functions are Lipschitz continuous. In particular, we have
\begin{equation}\label{eq:Lip_cont_g_j}
    g_j(x_1) - g_j(x_2) \le M_{g_j} \gnorm{x_1-x_2}{}{}, \quad \forall x_1, x_2 \in X, j \in [m].
\end{equation}
Note that the Lipschitz continuity assumption in \eqref{eq:Lip_cont_g_j} is common in the literature when $g_j$ is a nonsmooth function. If $g_j$ are Lipschitz smooth then their gradients are bounded due to compactness of $X$. Hence, \eqref{eq:Lip_cont_g_j} is not a strong assumption for the given setting. Also note that due to the definition of subgradient for convex function, we have $\gnorm{ g_j'(\cdot)}{}{} \le M_{g_j}$ which implies $\abs{  g_j'(x_2)^T(x_1-x_2) } \le \gnorm{ g_j'(x_2)}{}{} \gnorm{x_1-x_2}{}{} \le M_{g_j}\gnorm{x_1-x_2}{}{}$. Using this relation and \eqref{eq:Lip_cont_g_j}, we have for all $x_1,x_2 \in X$
\begin{equation}\label{eq:M_g_Lipschitz}
        \gnorm{g(x_1) - g(x_2)}{}{} \le M_g\gnorm{x_1 - x_2}{}{}; \qquad
        \gnorm{ g'(x_2)^T(x_1-x_2)}{}{} \le M_g\gnorm{x_1 - x_2}{}{},
\end{equation}
where $M_g := (\tsum_{i=1}^mM^2_{g_j})^{1/2}$.
Based on these definitions, we can classify FCVI problems into smooth and nonsmooth cases as follows:
\begin{enumerate}
    \item Smooth FCVI problem: We say that FCVI problem \eqref{eq:FCVI} is smooth if both operator $F$ and constraint functions $g_j, j =1, \dots, m$ are smooth, i.e., $H = H_{g_1} = \dots = H_{g_m} = 0$. This implies, for the smooth FCVI problem, we have $H = H_g  =0$.
    \item Nonsmooth FCVI problems: The FCVI problems which do not fall into the smooth case (and satisfies relations \eqref{eq:F-Lipschitz-property}, \eqref{eq:g_i-Lipschitz-property}) are called nonsmooth FCVI problems. Hence, at least one of $H > 0$ and $H_g > 0$ must hold for nonsmooth FCVI problems.
\end{enumerate} 
\paragraph{Convergence criteria for FCVIs}
Throughout this paper, we assume that the solution of \eqref{eq:FCVI}, denoted by $x^*$, exists. If the constraint $\{g(x) \le \zero\}$ were not present, then the resulting problem would be a standard VI problem with $\wt{X} = X$. 
Then, by definition of $x^*$, we have the following necessary and sufficient condition for the solution of VI
\[ F(x^*) + N_X(x^*) \ni \zero.\]
In the presence of function constraints $\{g(x) \le 0\}$, this condition can be extended to a KKT form under additional constraint qualifications such as Slater's condition:  There exists a vector $\lambda^* \in \Rbb^m_+$ such that
\begin{subequations}\label{eq:KKT_opt_FCVI}
	\begin{align}
		F(x^*) + \tsum_{j =1}^m\lambda^*_j \partial g_j(x^*) + N_X(x^*) &\ni \zero;\label{eq:opt_condn_FCVI}\\
		\lambda^*_jg_j(x^*) &= 0, \quad	\forall\ j =1, \dots,m,\label{eq:compl_slackness_condn_FCVI}
	\end{align}
\end{subequations}
where the sum in the first equation is the Minkowski sum of sets. See, e.g., \cite[Section 1.3.2 and 3.2]{facchinei2003finite} and survey article \cite{harker1990finite} that investigate various constraint qualifications that show existence of $\lambda^*$ in \eqref{eq:KKT_opt_FCVI}. We make the following assumption.
\begin{Assumption}\label{Assump:Assumption-1}
Throughout this paper, we assume that $(x^*, \lambda^*)$ satisfying \eqref{eq:KKT_opt_FCVI} exists.
\end{Assumption}
Consider the following Lagrangian function
\begin{equation}\label{eq:Lag_fun_FCVI}
    \Lcal(x, \lambda) := \inprod{F(x^*)}{x} + \inprod{\lambda} {g(x)}.  
\end{equation}
We have the following immediate consequences of KKT system \eqref{eq:KKT_opt_FCVI} for $\Lcal$, which motivates (to a large extent) some of the methods discussed in this paper. 
\begin{prop}\label{prop:FCVI_SPP_reformulation}
    Under Assumption \ref{Assump:Assumption-1}, we have that $(x^*, \lambda^*)$ is the solution of the the saddle-point problem 
    \begin{equation} \label{eq:FCVI_SPP_reformulation}
	  	\min_{x \in X} \max_{\lambda \ge \zero} \Lcal(x, \lambda).
  	\end{equation}%
\end{prop}
\begin{Proof}
    It is easy to see that $\Lcal(\cdot, \lambda)$ is a convex in $x$ for all $\lambda \in \Rbb^m_+$ and $\Lcal(x, \cdot)$ is concave in $\lambda$ for all $x\in X$. We denote the subgradient of $\Lcal$ with respect to the first argument $x$ as $ \Lcal'_x(x, \lambda)$. Similarly, we express the gradient with respect to the second argument $\lambda$  as $\grad_\lambda\Lcal(x, \lambda)$. We observe that \eqref{eq:KKT_opt_FCVI} is equivalent to
\begin{equation}\label{eq:sp_condition_Lag_func}
    \Lcal(x, \lambda^*) \ge \Lcal(x^*, \lambda^*) \ge \Lcal(x^*, \lambda), \quad \forall x \in X, \lambda \in \Rbb^m_+.
\end{equation}
Indeed, the first inequality is equivalent to \eqref{eq:opt_condn_FCVI} due to the convexity of  $\Lcal(\cdot, \lambda^*)$, and it is easy to see that the second inequality is equivalent to the complementary slackness condition in \eqref{eq:compl_slackness_condn_FCVI}. From the above relation, we see that the pair $(x^*, \lambda^*)$ is 
the saddle point solution of the following min-max problem
\end{Proof}
Based on this discussion, it is clear that the saddle point problem above is an equivalent reformulation of FCVI~\eqref{eq:FCVI} under Assumption \ref{Assump:Assumption-1}. 
Below, we formally define the $\epsilon$-solution of the FCVI problem. We also consider the stochastic case in the following definition.
\begin{Defi}\label{def:appx-soln}
	We say that $\wt{x}$ is an $\epsilon$-solution of \eqref{eq:FCVI} if
	\begin{equation}
		\max_{x\in \wt{X}}\,\inprod{F(x)}{\wt{x}-x} \le \epsilon , 
        \quad
		\gnorm{\bracket{g(\wt{x})}_+}{}{} \le \epsilon. 
        \label{eq:weak-VI}
	\end{equation}
    Similarly, we say that $\wt{x}$ is a stochastic $\epsilon$-solution of \eqref{eq:FCVI} if 
    \begin{equation}
		\Ebb\,[\max_{x\in \wt{X}}\,\inprod{F(x)}{\wt{x}-x}] \le \epsilon,
        \quad
		\Ebb\,[\gnorm{\bracket{g(\wt{x})}_+}{}{}] \le \epsilon, 
        \label{eq:stoch_weak-VI}
	\end{equation}
 where the expectation is taken over the randomness of the algorithm used to generate $\wt{x}$.
\end{Defi}
When $\epsilon = 0$, we have $g(\wt{x}) \le 0$ implying that $\wt{x}$ is feasible with respect to function constraints and $\inprod{F(x)}{x - \wt{x}} \ge 0$ for all $x \in \wt{X}$. This is known as the weak-gap criterion. For monotone operators, such $\wt{x}$ are identical with solution set of \eqref{eq:FCVI} \citep[Prop. 2.3.15]{facchinei2003finite}.
When the function constraint $g(x)\le \zero$ is not present (i.e., $g(x) \equiv 0$ for all $x$), the criteria in \eqref{eq:weak-VI} 
reduces to $\max_{x \in X} \inprod{F(x)}{\wt{x} -x} \le \epsilon$. This criterion is commonly referred to as the  ``approximate weak-gap'' in the VI literature for the set-constrained problems~\citep{Nemirovski2004}. Moreover, its stochastic counterpart in \eqref{eq:stoch_weak-VI} 
reduces to $\Ebb[\max_{x \in X} \inprod{F(x)}{\wt{x} -x}]$, which is a widely adopted stochastic weak-gap criterion in the SVI literature~\citep{Kotsalis2022, boob2021optimal}. When $g(x)$ is present, the approximate feasibility criterion $\gnorm{[g(\wt{x})]_+}{}{}$ is widely used in convergence rate analysis of function-constrained optimization \citep{aravkin2019level,boob2023stochastic}. 
\section{The Adaptive Operator Extrapolation method }\label{sec:adaptive FCVI}
In this section, we consider the smooth and deterministic FCVI$(F, X, g_j, j \in [m])$ problem. In view of the KKT system \eqref{eq:KKT_opt_FCVI}, we see that FCVI$(F, X, g_j, j \in [m])$ is equivalent to VI$(\wt{F}(x, \lambda), X \times \{\lambda \ge \zero\})$ where 
\[\wt{F}(x, \lambda) = \begin{bmatrix} F(x) + \grad g(x) \lambda \\ -g(x)\end{bmatrix}.\] 
Here, we note that $g$ is a smooth function; hence, $\grad g(x)$ is well-defined. It is easy to see that $\wt{F}(x, \lambda)$ is a monotone operator on the set $X\times \Rbb^m_+$. Indeed, $\inprod{\lambda}{g(x)}$ is convex-concave function on $X\times \Rbb^m_+$. Hence, its gradient operator $\begin{pmatrix}
    \grad g(x)\lambda\\-g(x)
\end{pmatrix}$ is a monotone operator on $X\times \Rbb^m_+$. Noting that $\begin{pmatrix}
    F(x)\\ \zero
\end{pmatrix}$ is also a monotone mapping on $X\times \Rbb^m_+$, we have that $\wt{F}(x, \lambda)$ is monotone. We can now apply the Operator Extrapolation (\opex{}) method \cite{Kotsalis2022} for VI$(\wt{F}, X\times \Rbb^m_+)$. However, there is a major challenge. The convergence of \opex{} method critically relies on the Lipschitz continuity of $\wt{F}$. Unfortunately, even for smooth FCVI, the operator $\wt{F}$ is not Lipschitz continuous. In fact, it is not even bounded. Indeed, even when $F$ is Lipschitz continuous and $g$ is smooth, i.e., $H = H_g = 0$, we have
\begin{align}
    &\gnorm{\wt{F}(x_1, \lambda_1) - \wt{F}(x_2, \lambda_2)}{}{}
    \le \gnorm{F(x_1) - F(x_2)}{}{} + \gnorm{\grad g(x_1)\lambda_1 - \grad g(x_2)\lambda_2}{}{} + \gnorm{g(x_1) - g(x_2)}{}{} \nonumber\\
    &\quad \le L\gnorm{x_1-x_2}{}{} + \gnorm{(\grad g(x_1) - \grad g(x_2))\lambda_1}{}{} + \gnorm{\grad g(x_2)(\lambda_1-\lambda_2)}{}{} + \gnorm{g(x_1) - g(x_2)}{}{} \nonumber\\
    &\quad \le (L+M_g+L_g\gnorm{\lambda_1}{}{})\gnorm{x_1-x_2}{}{} +  M_g\gnorm{\lambda_1-\lambda_2}{}{}, \label{eq:Lip_cont_F_tilde}
\end{align}
where the first inequality follows by \eqref{eq:F-Lipschitz-property} and the triangle inequality of $\gnorm{\cdot}{}{}$, and second inequality follows due to \eqref{eq:Lip_cont_g_j}. This implies for smooth FCVI, Lipschitz constant of $\wt{F}$, denoted by $\wt{L}$ is of $O(L+M_g + L_g\gnorm{\lambda}{}{})$. This is not a uniform bound that works for all $\lambda \ge \zero$. Hence, the operator $\wt{F}$ does not satisfy uniform Lipschitz continuity, which is a required assumption for the convergence of the operator extrapolation method \cite[Eq (1.2)]{Kotsalis2022}. 

The critical issue here is that the set of $\lambda$ is unbounded. In general, the $\lambda$ iterates of \opex{} algorithm applied for VI$(\wt{F}, X\times \Rbb^m_+)$ may also be unbounded if not handled carefully. To handle this challenging issue one can impose an artificial bound $\gnorm{\lambda}{}{} \le R$ where $R$ is some radius. As long as $R > \gnorm{\lambda^*}{}{}$, it is easy to argue that the saddle point of this modification of\eqref{eq:Lag_fun_FCVI} is $(x^*, \lambda^*)$ and we can still converge to the solution of FCVI \eqref{eq:FCVI}. However, finding such a bound $R$ can be challenging in the first place. To address the unboundedness of $\lambda$, we take a different approach. Note that we need to use the Lipschitz continuity of $\wt{F}$ in \eqref{eq:Lip_cont_F_tilde} not for all $\lambda \ge \zero$, but only for the $\lambda^t$ visited (or generated) by the algorithm. Hence, if we can ensure that the sequence of $\{\lambda^t\}_{t \ge 0}$ generated by the \opex{} method remains bounded, then we can ensure the optimal convergence rate of $O(\tfrac{1}{T})$, even if the $\wt{F}$ is not uniformly smooth. 

Keeping this issue of unbounded $\lambda$ in sight, we present a novel algorithm that adaptively searches for the optimal dual multiplier for smooth (deterministic) FCVIs. We refer to this algorithm as the Adaptive Operator Extrapolation (\adopconex) method (see Algorithm \ref{alg:alg2}). In particular, the Operator Extrapolation method applied to $\wt{F}$ has the following basic update scheme:
\[
(x^{t+1}, \lambda^{t+1}) = \argmin_{x\in X, \lambda \ge \zero} \inprod*{(1+\theta_t)\wt{F}(x^t, \lambda^t) - \theta_t\wt{F}(x^{t-1}, \lambda^{t-1})}{\begin{pmatrix}
    x\\ \lambda
\end{pmatrix}} + \frac{\eta_t}{2}\gnorm*{ \begin{pmatrix}
    x\\ \lambda
\end{pmatrix} - \begin{pmatrix}
    x^t\\ \lambda^t
\end{pmatrix} }{}{2}.\]
Due to the separability of the objective and constraint set $X\times \Rbb^m_+$ into two blocks of $x$ and $\lambda$, we can separate the above minimization problem into two separate updates for $x$ and $\lambda$:
\begin{align*}
    x^{t+1} &= \argmin_{x \in X} \inprod{(1+\theta_t)[F(x^t) + \grad g(x^t)\lambda^t] - \theta_t[F(x^{t-1}) + \grad g(x^{t-1})\lambda^{t-1}]}{x} + \tfrac{\eta_t}{2}\gnorm{x-x^t}{}{2}\\
    \lambda^{t+1} &= \argmin -\inprod{(1+\theta_t)g(x^t) - \theta_t g(x^{t-1})}{\lambda} + \tfrac{\eta_t}{2}\gnorm{\lambda - \lambda^t}{}{2}
\end{align*}
Algorithm \ref{alg:alg2} makes an important modification of the above updates where the stepsize parameter for the update of $\lambda^{t+1}$ is different from $\eta_t$. This leads to a two-timescale update scheme of the Operator Extrapolation method. Moreover, this method is called adaptive as we will set $\eta_t, \tau_t$ depending on the historical information $\{\lambda^i\}_{i = 0}^{t}$ obtained from the run of the algorithm. 
	\vspace{-1em}
	\begin{algorithm}[H]
	\begin{algorithmic}[1]
		\State {\bf Input: }$x^{0} \in X, \lambda^0 = \zero$. 
		\State Set $x^{-1} = x^{0}$.
		\For{$t=0,1,2,\dots,T-1$}
		\State $s^{t} \gets (1+\theta_t)g(x^t) - \theta_t g(x^{t-1})$
		\State $\lambda^{t+1} \gets \argmin_{\lambda\ge \zero} \inprod{-s^t}{\lambda} + \tfrac{\tau_t}{2}\gnorm{\lambda-\lambda^t}{}{2}$
		\State $u^t \gets (1+\theta_t)[F(x^t)+ \grad g(x^t)\lambda^t] - \theta_t[F(x^{t-1})+ \grad g(x^{t-1})\lambda^{t-1} ]$
		\State $x^{t+1} \gets \argmin_{x \in X} \inprod{u^t }{x} + \tfrac{\eta_t}{2}\gnorm{x-x^t}{}{2}$
		\EndFor
		\State {\bf Output: }$\wb{x}^T := \paran{\tsum_{t=0}^{T-1} \gamma_{t}x^{t+1}}\big/\paran{\tsum_{t=0}^{T-1}\gamma_t}$
	\end{algorithmic}
	\caption{Adaptive Operator Extrapolation (\adopconex) method} \label{alg:alg2}
\end{algorithm}
\vspace{-1.5em}
We mention a useful lemma that we use throughout this paper.
\begin{lem}
	\label{lem:3-point-mo}
	Assume that $h: S\to \Rbb$ 
	satisfies
	\begin{equation*} \label{eq:strong_conv_g}
		h(y) \ge h(x) + \inprod{\grad h(x)}{y-x} + \tfrac{\mu}{2}\gnorm{y-x}{}{2} , \quad \forall x, y \in S
	\end{equation*}
	for some $\mu \ge 0$, where $S$ is a closed convex set. If \[ \wb{x}_+ = \argmin_{x \in S} h(x) + \tfrac{1}{2}\gnorm{x-\wb{x}}{}{2},\]
	then 
	\[ h(\wb{x}_+) + \tfrac{1}{2}\gnorm{\wb{x}_+-\wb{x}}{}{2} + \tfrac{(\mu+1)}{2}\gnorm{x-\wb{x}_+}{}{2} \le h(x) + \tfrac{1}{2}\gnorm{x-\wb{x}}{}{2}, \quad \forall x \in S. \]
\end{lem}
Now, we present a crucial technical lemma to prove the convergence of the \adopconex{} method.
\begin{lem}\label{lem:adaptive_alg_basic_lemma}
    Let the parameters $\{\gamma_t, \theta_t, \eta_t, \tau_t\}_{t \ge 0}$ of Algorithm \ref{alg:alg2} satisfy
     \begin{equation}\label{eq:step_condn_1_adap}
			\gamma_t\eta_t \le \gamma_{t-1}\eta_{t-1},\qquad
			\gamma_{t}\tau_t \le \gamma_{t-1}\tau_{t-1},\qquad
			\gamma_t\theta_t =\gamma_{t-1}.
	\end{equation}
	Then, we have for all $x \in X,\ \lambda \ge \zero$,
	\begin{align}
		\allowdisplaybreaks
		&\Gamma_{t+1}[\inprod{F(x)}{\wb{x}^{t+1} -x} + \inprod{g(\wb{x}^{t+1})}{\lambda} - \inprod{g(x)}{\wb{\lambda}^{t+1}}] \nonumber\\
		&\le [\tfrac{\gamma_0\eta_0}{2} \gnorm{x-x^0}{}{2} + \tfrac{\gamma_0\tau_0}{2}\gnorm{\lambda- \lambda^0}{}{2} + \gamma_t\inprod{q_x^{t+1} + \Delta F_{t+1}}{x^{t+1}-x} - \gamma_t\inprod{q_\lambda^{t+1}}{\lambda^{t+1}-\lambda} ]\nonumber\\
		&- [ \tfrac{\gamma_t\eta_t}{2}\gnorm{x-x^{t+1}}{}{2} + \tfrac{\gamma_t\tau_t}{2} \gnorm{\lambda-\lambda^{t+1}}{}{2}] \nonumber\\
		&+ \tsum_{i=0}^t \bracket[\big]{\gamma_{i-1}\inprod{q_\lambda^i}{\lambda^{i+1}-\lambda^i} - \gamma_{i-1}  \inprod{q_x^i + \Delta F_i}{x^{i+1}-x^i} - \tfrac{\gamma_i\eta_i}{2} \gnorm{x^i-x^{i+1}}{}{2} -  \tfrac{\gamma_i\tau_i}{2}\gnorm{\lambda^{i+1}-\lambda^i}{}{2}} \label{eq:adaptive_conv_basic},
	\end{align}
	where $\wb{x}^{t+1} := \tfrac{1}{\Gamma_{t+1}}\tsum_{i=0}^t\gamma_ix^{i+1}$, $\Gamma_{t+1} := \tsum_{i=0}^t\gamma_i$, $q_\lambda^i := g(x^i) - g(x^{i-1})$, $q_x^i := \grad g(x^i) \lambda^i - \grad g(x^{i-1}) \lambda^{i-1}$ and $\Delta F_i := F(x^i)- F(x^{i-1})$.
\end{lem}
\begin{Proof}
		Note that $\lambda^{i+1} = \argmin\limits_{\lambda \ge \zero} \inprod{-s^i}{\lambda} + \tfrac{\tau_i}{2} \gnorm{\lambda-\lambda^i}{2}{2}$. Hence, using Lemma \ref{lem:3-point-mo} and definition of $q_\lambda^i$, we have for all $\lambda \ge \zero$,
	\begin{align}
		-\inprod{s^i}{\lambda^{i+1} -\lambda} &\le \tfrac{\tau_i}{2} \bracket*{\gnorm{\lambda-\lambda^i}{}{2} - \gnorm{\lambda^{i+1}-\lambda^i}{}{2} - \gnorm{\lambda-\lambda^{i+1}}{}{2}} \nonumber\\
		\Rightarrow
		\inprod{-g(x^{i+1})}{\lambda^{i+1}-\lambda} &\le \tfrac{\tau_i}{2} \bracket*{\gnorm{\lambda-\lambda^i}{}{2} - \gnorm{\lambda-\lambda^{i+1}}{}{2}} - \inprod{q_\lambda^{i+1}}{\lambda^{i+1}-\lambda} + \theta_i \inprod{q_\lambda^i}{\lambda^i-\lambda} \nonumber\\
		& \quad+ \theta_i \inprod{q_\lambda^i}{\lambda^{i+1}-\lambda^i}- \tfrac{\tau_i}{2}\gnorm{\lambda^{i+1}-\lambda^i}{}{2} \label{eq:int_rel21}.
	\end{align}
	Using optimality of $x^{i+1}$ along with Lemma \ref{lem:3-point-mo} and noting the definitions of $q_x^i, \Delta F_i$, we have for all $x \in X$
	\begin{align}
		\inprod{u^i}{x^{i+1}-x} &\le \tfrac{\eta_i}{2}\bracket[\big]{\gnorm{x-x^i}{}{2} - \gnorm{x-x^{i+1}}{}{2} - \gnorm{x^i-x^{i+1}}{}{2}} \nonumber\\
		\inprod{F(x^{i+1})}{x^{i+1}-x} &+ \inprod{\grad g(x^{i+1})\lambda^{i+1}}{x^{i+1}-x} \nonumber\\
		&\le \tfrac{\eta_i}{2}\bracket[\big]{\gnorm{x-x^i}{}{2} - \gnorm{x-x^{i+1}}{}{2}}
		+  \inprod{q_x^{i+1} + \Delta F_{i+1}}{x^{i+1}-x} - \theta_i \inprod{q_x^i + \Delta F_i}{x^i-x} \nonumber\\
		&\quad - \tfrac{\eta_i}{2} \gnorm{x^i-x^{i+1}}{}{2}-  \theta_i \inprod{q_x^i + \Delta F_i}{x^{i+1}-x^i} \label{eq:int_rel22}.
	\end{align}
	Moreover, using convexity of $g$, monotonicity of $F$ and the fact that $\lambda^{i+1} \ge \zero$, we have,
	\begin{subequations}\label{eq:lower_bds_potential_fun}
		\begin{align}
			\inprod{\grad g(x^{i+1})}{x^{i+1}-x} &\ge g(x^{i+1}) -g(x), &\forall \ x\in X, \nonumber\\
			\Rightarrow\inprod{\grad g(x^{i+1})\lambda^{i+1}}{x^{i+1}-x} &\ge \inprod{\lambda^{i+1}}{g(x^{i+1}) - g(x)},  &\forall \ x\in X, \label{eq:int_rel23}\\
			\inprod{F(x^{i+1})}{x^{i+1}-x} &\ge \inprod{F(x)}{x^{i+1}-x} , &\forall \ x\in X. \label{eq:int_rel24}
		\end{align}
	\end{subequations}
	Using relations in \eqref{eq:lower_bds_potential_fun} inside \eqref{eq:int_rel22}, we have for all $x \in X, \lambda \ge \zero$,
	\begin{align}
		\inprod{F(x)}{x^{i+1}-x}&+ \inprod{\lambda^{i+1}}{g(x^{i+1}) - g(x)}\nonumber\\
		&\le \tfrac{\eta_i}{2}\bracket[\big]{\gnorm{x-x^i}{}{2} - \gnorm{x-x^{i+1}}{}{2}}
		+  \inprod{q_x^{i+1} + \Delta F_{i+1}}{x^{i+1}-x} - \theta_i \inprod{q_x^i + \Delta F_i}{x^i-x} \nonumber\\
		&\quad - \tfrac{\eta_i}{2} \gnorm{x^i-x^{i+1}}{}{2}-  \theta_i \inprod{q_x^i + \Delta F_i}{x^{i+1}-x^i} \label{eq:int_rel25}.
	\end{align}
	Summing \eqref{eq:int_rel21} and \eqref{eq:int_rel25}, then multiplying the resulting relation by $\gamma_i$ and noting the third relation in \eqref{eq:step_condn_1_adap}, we have
	\begin{align}
		\gamma_{i}[&\inprod{F(x)}{x^{i+1}-x} + \inprod{g(x^{i+1})}{\lambda} - \inprod{g(x)}{\lambda^{i+1}}] \nonumber\\
		&\le [\tfrac{\gamma_i\eta_i}{2} \gnorm{x-x^i}{}{2} + \tfrac{\gamma_i\tau_i}{2}\gnorm{\lambda- \lambda^i}{}{2} + \gamma_i\inprod{q_x^{i+1} + \Delta F_{i+1}}{x^{i+1}-x} - \gamma_i\inprod{q_\lambda^{i+1}}{\lambda^{i+1}-\lambda} ]\nonumber\\
		&\quad- [ \tfrac{\gamma_i\eta_i}{2}\gnorm{x-x^{i+1}}{}{2} + \tfrac{\gamma_i\tau_i}{2} \gnorm{\lambda-\lambda^{i+1}}{}{2} + \gamma_{i-1} \inprod{q_x^i + \Delta F_i}{x^i-x}  - \gamma_{i-1}\inprod{q_\lambda^i}{\lambda^i-\lambda} ] \nonumber\\
		&\quad +  \gamma_{i-1}\inprod{q_\lambda^i}{\lambda^{i+1}-\lambda^i} - \gamma_{i-1}  \inprod{q_x^i + \Delta F_i}{x^{i+1}-x^i} - \tfrac{\gamma_i\eta_i}{2} \gnorm{x^i-x^{i+1}}{}{2} -  \tfrac{\gamma_i\tau_i}{2}\gnorm{\lambda^{i+1}-\lambda^i}{}{2}
	\end{align}
	Summing the above relation from $i = 0$ to $i = t$, using \eqref{eq:step_condn_1_adap} along with Jensen's inequality and noting that $q_x^0 = \Delta F_0 = \zero$ and $q_\lambda^0 = \zero$, we have \eqref{eq:adaptive_conv_basic}. Hence, we conclude the proof.
\end{Proof}
The above lemma provides a conceptual bound on the potential function used to analyze the VI problem. To convert this lemma into useful bounds, we need to analyze the potential function as well as additional inner product terms on the right-hand side of the above lemma. In the following lemma, we present some required conditions that lead to the desired bounds.
\begin{thm} \label{lem:adaptive_mu=0_final}
	Let $\gamma_t = \tfrac{\eta_0}{\eta_t}$ for $t \ge 0$, $\theta_t = \tfrac{\gamma_{t-1}}{\gamma_t}$, $\tau_t = \beta \eta_t$ and $\eta_t = c_1L +c_2L_g \max_{i \in [t]} \gnorm{\lambda^i}{}{}$ where $c_1, c_2, \beta$ satisfies the conditions 
	\begin{equation}\label{eq:elementary_ineq_constants}
			\tfrac{c_1}{3} \ge \tfrac{c_1}{c_2} + 1,\qquad
			\beta = \tfrac{12M_g^2}{c_1^2L^2}.
	\end{equation}
	Then, we have,
	\begin{align}
		\gnorm{\lambda^{t+1}}{}{} &\le B :=\sqrt{\tfrac{2}{\beta}}\gnorm{x^0-x^*}{}{}+ (\sqrt{2}+1) \gnorm{\lambda^*}{}{}, \qquad \forall t \ge 0,\label{eq:bound_lambda_adaptive}\\
		\Gamma_{T} \inprod{F(x)}{\wb{x}^T-x} &\le \tfrac{c_1L}{2} \gnorm{x-x^0}{}{2}, \qquad \forall \ x \in \wt{X},\label{eq:apx_opt_basic_adaptive}\\
		\Gamma_{T} \gnorm{[g(\wb{x}^T)]_+}{}{} &\le \tfrac{c_1L}{2} \gnorm{x^*-x^0}{}{2}  + \tfrac{\beta c_1L}{2}(\gnorm{\lambda^*}{}{}+1)^2 \label{eq:apx_feas_basic_adaptive}
	\end{align}
	while $\Gamma_{T} \ge \tfrac{c_1L}{c_1L + c_2L_gB} T$.
\end{thm}
\begin{Proof}
	Using the definitions of $\gamma_t$ and $\eta_t$, we have 
	\begin{align}\label{eq:basic_step_relations}
			\gamma_i\eta_i &= \eta_0 = c_1L, 
			\quad \gamma_i\tau_i = \beta \eta_0 = \beta c_1L, \quad \forall i \ge0,
            \nonumber\\
			\gamma_{i-1} &= \tfrac{\eta_0}{\eta_{i-1}} = \tfrac{c_1L}{c_1L + c_2L_g \max_{j \in [i-1]}\gnorm{\lambda^j}{}{}} \le \tfrac{c_1L}{c_2L_g \gnorm{\lambda^{i-1}}{}{}},\\
			\gamma_{i} &= \gamma_{i-1} \tfrac{\eta_{i-1}}{\eta_i}  \le \gamma_{i-1} \le \gamma_0 = 1.\nonumber
	\end{align}
	Using definitions of $q_\lambda^i, q_x^i$ and $\Delta F_i$ along with Cauchy-Schwarz inequality, we have
		\begin{align*}
			\gnorm{q_\lambda^i}{}{} &\le M_g\gnorm{x^i-x^{i-1}}{}{},\\
			\gnorm{q_x^i}{}{} &\le \gnorm{\grad g(x^i)[\lambda^i- \lambda^{i-1}] + [\grad g(x^i) - \grad g(x^{i-1})]\lambda^{i-1}}{}{} 
            \le M_g \gnorm{\lambda^i- \lambda^{i-1}}{}{} + L_g \gnorm{\lambda^{i-1}}{}{} \gnorm{x^i-x^{i-1}}{}{},\\
			\gnorm{\Delta F_i}{}{} &\le L\gnorm{x^i-x^{i-1}}{}{}.
		\end{align*}
	Using above relations along with \eqref{eq:basic_step_relations}, we have
	\begin{subequations}\label{eq:bounds_q_inner_prod}
		\begin{align}
			\gamma_{i-1} \inprod{q_\lambda^i}{\lambda - \lambda^i} &\le \gamma_{i-1} M_g\gnorm{x^i-x^{i-1}}{}{} \gnorm{\lambda - \lambda^i}{}{} \le M_g\gnorm{x^i-x^{i-1}}{}{} \gnorm{\lambda - \lambda^i}{}{}, \label{eq:bounds_q_inner_prod-1}\\
			\gamma_{i-1}  \inprod{q_x^i}{x^i-x} &\le \gamma_{i-1}M_g \gnorm{\lambda^i- \lambda^{i-1}}{}{} \gnorm{x-x^i}{}{} + \gamma_{i-1}L_g \gnorm{\lambda^{i-1}}{}{} \gnorm{x^i-x^{i-1}}{}{}\gnorm{x-x^i}{}{} \nonumber\\
			&\le M_g \gnorm{\lambda^i- \lambda^{i-1}}{}{}\gnorm{x-x^i}{}{} + \tfrac{c_1L}{c_2}\gnorm{x^i-x^{i-1}}{}{}\gnorm{x-x^i}{}{}, \label{eq:bounds_q_inner_prod-2}\\
			\gamma_{i-1}  \inprod{\Delta F_i}{x^i - x} &\le  \gamma_{i-1}L\gnorm{x^i-x^{i-1}}{}{}\gnorm{x-x^i}{}{} \le L\gnorm{x^i-x^{i-1}}{}{}\gnorm{x-x^i}{}{}. \label{eq:bounds_q_inner_prod-3}
		\end{align}
	\end{subequations}
	Summing \eqref{eq:bounds_q_inner_prod-1} with $\lambda = \lambda^{i+1}$, \eqref{eq:bounds_q_inner_prod-2}-\eqref{eq:bounds_q_inner_prod-3} with $x = x^{i+1}$ and then summing the resulting relation from $i = 0$ to $t$, we have,
	\begin{align}
		&\tsum_{i=0}^t\big[\gamma_{i-1}[\inprod{q_\lambda^i}{\lambda^{i+1} - \lambda^i} - \inprod{q_x^i + \Delta F_i}{x^{i+1}-x^i}] - \tfrac{\eta_0}{2} \gnorm{x^{i+1} -x^i}{}{2} - \tfrac{\beta\eta_0}{2} \gnorm{\lambda^{i+1} - \lambda^i}{}{2} \big]\nonumber\\ 
		&\le \tsum_{i=0}^t \Big[\paran[\big]{ \tfrac{c_1}{c_2} +1}L\gnorm{x^i-x^{i-1}}{}{}\gnorm{x^{i+1}-x^i}{}{} +M_g\gnorm{x^i-x^{i-1}}{}{} \gnorm{\lambda^{i+1} - \lambda^i}{}{}  \nonumber\\
		&\hspace{4em}+ M_g \gnorm{\lambda^i- \lambda^{i-1}}{}{}\gnorm{x^{i+1}-x^i}{}{} - \tfrac{c_1L}{2} \gnorm{x^{i+1} -x^i}{}{2} - \tfrac{\beta c_1L}{2} \gnorm{\lambda^{i+1} - \lambda^i}{}{2} \Big]\nonumber\\
		&\le \tsum_{i=0}^t \Big[ \paran[\big]{ \tfrac{c_1}{c_2} +1}L\gnorm{x^i-x^{i-1}}{}{}\gnorm{x^{i+1}-x^i}{}{} - \tfrac{c_1L}{6}\gnorm{x^i-x^{i-1}}{}{2} - \tfrac{c_1L}{6} \gnorm{x^{i+1}-x^i}{}{2} \nonumber\\
		&\quad +M_g\gnorm{x^i-x^{i-1}}{}{} \gnorm{\lambda^{i+1} - \lambda^i}{}{}  - \tfrac{c_1L}{12}\gnorm{x^i-x^{i-1}}{}{2} - \tfrac{\beta c_1L}{4} \gnorm{\lambda^{i+1} - \lambda^i}{}{2} \nonumber\\
		&\quad +M_g\gnorm{\lambda^i- \lambda^{i-1}}{}{}\gnorm{x^{i+1}-x^i}{}{}  - \tfrac{c_1L}{12}\gnorm{x^{i+1}-x^i}{}{2} - \tfrac{\beta c_1L}{4} \gnorm{\lambda^i- \lambda^{i-1}}{}{2} \Big] \nonumber\\
        &\quad- \tfrac{c_1L}{4} \gnorm{x^{t+1}-x^t}{}{2} -\tfrac{\beta c_1L}{4} \gnorm{\lambda^{t+1} - \lambda^t}{}{2} \nonumber\\
		&\le - \tfrac{c_1L}{4} \gnorm{x^{t+1}-x^t}{}{2} -\tfrac{\beta c_1L}{4} \gnorm{\lambda^{t+1} - \lambda^t}{}{2} \label{eq:int_rel26},
	\end{align}
	where last inequality follows due to \eqref{eq:elementary_ineq_constants} and Young's inequality. 
	
	Note that the step-size policy $\{\gamma_{t}, \eta_t, \tau_t\}$ satisfies the relation \eqref{eq:step_condn_1_adap} and all requirements of Lemma \ref{lem:adaptive_alg_basic_lemma}. Hence, using the final result of Lemma \ref{lem:adaptive_alg_basic_lemma} along with \eqref{eq:bounds_q_inner_prod} for $i = t+1$, we have
	\begin{align}
		\Gamma_t&[\inprod{F(x)}{\wb{x}^{t+1} -x} + \inprod{g(\wb{x}^{t+1})}{\lambda} - \inprod{g(x)}{\wb{\lambda}^{t+1}}] \nonumber\\
		&\le \tfrac{c_1L}{2} \gnorm{x-x^0}{}{2} + \tfrac{\beta c_1L}{2} \gnorm{\lambda - \lambda^0}{}{2} +  \gamma_t\inprod{q_x^{t+1} + \Delta F_{t+1}}{x^{t+1}-x} - \gamma_t\inprod{q_\lambda^{t+1}}{\lambda^{t+1}-\lambda}\nonumber\\
		&\quad - \tfrac{c_1L}{2}\gnorm{x-x^{t+1}}{}{2} - \tfrac{\beta c_1L}{2} \gnorm{\lambda-\lambda^{t+1}}{}{2} - \tfrac{c_1L}{4} \gnorm{x^{t+1}-x^t}{}{2} -\tfrac{\beta c_1L}{4} \gnorm{\lambda^{t+1} - \lambda^t}{}{2} \nonumber\\
		& \le \tfrac{c_1L}{2} \gnorm{x-x^0}{}{2} + \tfrac{\beta c_1L}{2} \gnorm{\lambda - \lambda^0}{}{2} \nonumber\\
		&\quad + \big[ M_g \gnorm{\lambda^{t+1}- \lambda^t}{}{}\gnorm{x-x^{t+1}}{}{} -\tfrac{\beta c_1L}{4} \gnorm{\lambda^{t+1} - \lambda^t}{}{2} - \tfrac{c_1L}{4}\gnorm{x-x^{t+1}}{}{2} \big]\nonumber\\
		&\quad + \big[\big( \tfrac{c_1}{c_2}+1 \big)L\gnorm{x^{t+1}-x^t}{}{}\gnorm{x-x^{t+1}}{}{}  - \tfrac{c_1L}{8} \gnorm{x^{t+1}-x^t}{}{2} - \tfrac{c_1L}{4}\gnorm{x-x^{t+1}}{}{2} \big] \nonumber\\
		&\quad + \big[M_g\gnorm{x^{t+1}-x^t}{}{} \gnorm{\lambda - \lambda^{t+1}}{}{} - \tfrac{c_1L}{8} \gnorm{x^{t+1}-x^t}{}{2}  - \tfrac{\beta c_1L}{4} \gnorm{\lambda-\lambda^{t+1}}{}{2} \big] 
        - \tfrac{\beta c_1L}{4} \gnorm{\lambda-\lambda^{t+1}}{}{2} 
        \nonumber\\
		&
        \le  \tfrac{c_1L}{2} \gnorm{x-x^0}{}{2} + \tfrac{\beta c_1L}{2} \gnorm{\lambda - \lambda^0}{}{2} -  \tfrac{\beta c_1L}{4} \gnorm{\lambda-\lambda^{t+1}}{}{2} , \label{eq:main_conve_adaptive_technical}
	\end{align}
	where the last inequality follows due to \eqref{eq:elementary_ineq_constants} and Young's inequality. Using $x =x*, \ \lambda = \lambda^*$ in the above relation and noting in view of $(x^*,\ \lambda^*)$ being a saddle point of \eqref{eq:FCVI_SPP_reformulation}, we have $\inprod{F(x^*)}{\wb{x}^{t+1} -x^*} + \inprod{g(\wb{x}^{t+1})}{\lambda^*} - \inprod{g(x^*)}{\wb{\lambda}^{t+1}} \ge 0$, and we obtain, 
	\begin{equation*}
		\tfrac{\beta c_1L}{4} \gnorm{\lambda^*-\lambda^{t+1}}{}{2} \le  \tfrac{c_1L}{2} \gnorm{x^*-x^0}{}{2} + \tfrac{\beta c_1L}{2} \gnorm{\lambda^*}{}{2},
	\end{equation*}
	which implies \eqref{eq:bound_lambda_adaptive}. 
	
	Now, since $\{\lambda^{t+1}\}$ remains bounded, we have $\gamma_{t} \ge \tfrac{c_1L}{c_1L + c_2L_gB}$ for all $t \ge 0$, implying that $\Gamma_{t+1} \ge \tfrac{c_1L}{c_1L + c_2L_gB} (t+1)$. Hence, after running the algorithm from $t = 0$ to $T-1$ and noting \eqref{eq:main_conve_adaptive_technical}, we have
	\begin{align}\label{eq:key_relation_adopex}
		\Gamma_{T}[\inprod{F(x)}{\wb{x}^T -x} + \inprod{g(\wb{x}^T)}{\lambda} - \inprod{g(x)}{\wb{\lambda}^T}] \le \tfrac{c_1L}{2} \gnorm{x-x^0}{}{2} + \tfrac{\beta c_1L}{2} \gnorm{\lambda - \lambda^0}{}{2}. 
	\end{align}
 	Using $x \in \wt{X}$ and $\lambda = \lambda^0 =\zero$, in \eqref{eq:key_relation_adopex} and noting that $\inprod{\wb{\lambda}^T}{-g(x)} \ge 0$ due to $\wb{\lambda}^T \ge \zero$, $-g(x) \ge \zero$ for all $x \in \wt{X}$, we have \eqref{eq:apx_opt_basic_adaptive}.
    Note that 
	\begin{align}
		&\inprod{F(x^*)}{\wb{x}^T-x^*} + (\gnorm{\lambda^*}{}{}+1)\gnorm{[g(\wb{x}^T)]_+}{}{} - \inprod{\wb{\lambda}^T}{g(x^*)} \nonumber\\
        &\ge \inprod{F(x^*)}{\wb{x}^T-x^*} + \inprod{\lambda^*}{g(\wb{x}^T)} - \inprod{\wb{\lambda}^T}{g(x^*)}
        + \gnorm{[g(\wb{x}^T)]_+}{}{} \nonumber\\
		&\ge \gnorm{[g(\wb{x}^T)]_+}{}{}, \label{eq:int_rel8_new}
	\end{align}
	where the first inequality follows due to Cauchy-Schwarz inequality and  $[g(\wb{x}^T)]_+ \ge g(\wb{x}^T)$, and the second inequality follows due to Proposition \ref{prop:FCVI_SPP_reformulation} (in particular, \eqref{eq:sp_condition_Lag_func}) yielding $\Lcal(\wb{x}^T, \lambda^*) - \Lcal(x^*, \wb{\lambda}^T) \ge 0$. 
    Using $x = x^*$ and $\lambda = (\gnorm{\lambda^*}{}{}+1) \tfrac{[g(\wb{x}^T)]_+}{\gnorm{[g(\wb{x}^T)]_+}{}{}}$ in \eqref{eq:key_relation_adopex}, adding the resulting relation with \eqref{eq:int_rel8_new} and noting $\lambda^0 = \zero$, we obtain \eqref{eq:apx_feas_basic_adaptive}. Hence, we conclude the proof.
\end{Proof}
The following corollary is in order. \vspace{-0.8em}
\begin{coro}\label{cor:conv_adaptive_OE}
	Using $c_1 = c_2 = 6$, we obtain an $O(1/T)$ convergence rate for Algorithm \ref{alg:alg2}.
\end{coro}
\begin{Proof}
	We only need to verify $c_1, c_2$ satisfies the first requirement of \eqref{eq:elementary_ineq_constants}. Then, by using Lemma~\ref{lem:adaptive_mu=0_final}, we can obtain the result. Hence, we conclude the proof.
\end{Proof}
\vspace{-1.8em}

\section{Partial Operator Extrapolation method}\label{sec:Deterministic-simple-PD}
Observe from Theorem \ref{lem:adaptive_mu=0_final} that to get a bound on the dual iterates in $\adopconex{}$, we need to set the stepsize policy in each iteration adaptively. In particular, the stepsizes $\eta_t$ and $\tau_t$ depends on the history of $\{\lambda_i\}_{i =1}^t$. The convergence of such a method critically relies on a descent type condition in \eqref{eq:main_conve_adaptive_technical}. However, it may not be possible to apply this idea directly to nonsmooth or stochastic cases. This section presents a different algorithm that uses a fixed stepsize policy and can still obtain a bound on the dual iterates. We refer to this algorithm as the Partial Operator Extrapolation (\sopconex) method since it only partially extrapolates operator $\wt{F}$. 
We present a unified analysis of \sopconex{} for the nonsmooth and fully stochastic FCVIs.

\subsection{The \sopconex{} method for the fully-stochastic FCVI problem}\label{sec:full-stoch} 
In this section, we assume the operator $F$, and function constraints $g$ 
are given in expectation form. Recall that we refer to this 
problem as the fully stochastic FCVIs. We use a stochastic oracle (SO) to generate random vectors estimating $F$, $g$, and $ g'(x)$. 
In particular, for a given $x$, the SO can generate a random vector  $\Ffrak(x,\xi)$ such that 
\begin{equation}\label{eq:SO_F_operator-mo}
	\forall x \in X, \quad F(x) = \Ebb[\Ffrak(x,\xi)], \quad \Ebb[\gnorm{\Ffrak(x, \xi) - F(x)}{}{2}] \le \sigma^2.
\end{equation}
  Essentially, $\Ffrak(x, \xi)$ is an unbiased estimator of $F(x)$ with a bounded variance. Moreover, SO can generate $\Gfrak(x,\xi)=[\Gfrak_1(x,\xi),\dots,\Gfrak_m(x,\xi)],\ \gfrak(x, \xi) = [\gfrak_1(x, \xi), \dots, \gfrak_m(x, \xi)]^\top$ which for all $x \in X$ satisfies
 \begin{equation} \label{eq:SO_g_operator-mo}
 	\begin{aligned}
 		\Ebb[\gfrak_j(x, \xi)] &= g_j(x), \quad \forall j \in [m],\\ 
 		\Ebb[\Gfrak_j(x, \xi)] &= \ g'_j(x), \quad \forall j \in [m],\\
 		\Ebb[\gnorm{\Gfrak_j(x, \xi) -  g'_j(x)}{}{2}] &\le \sigma_j^2, \quad \forall j \in [m],\\
 		\Ebb[\gnorm{\gfrak(x, \xi) - g(x)}{}{2}] &\le \sigma_\gfrak^2, 
 	\end{aligned}
 \end{equation}
where $\xi$ is (conditionally) independent of search point $x$ and the expectation is taken with respect to $\xi$. Note that the last relation in \eqref{eq:SO_g_operator-mo} is satisfied if the individual SO $\gfrak_j(x, \xi), \ j \in [m]$, satisfies $\Ebb[(\gfrak_j(x, \xi) - g_j(x))^2] \le \sigma_{\gfrak_j}^2$. In particular, we can set $\sigma_\gfrak^2 = \tsum_{j=1}^m \sigma_{\gfrak_j}^2$. 
We also assume that SO in this section satisfies
\begin{equation}\label{eq:SO_g_convexity}
	\gfrak_j(\cdot, \xi) \text{ is a convex lower semi-continuous function } \forall \xi.
\end{equation}
We will remove this assumption in Section \ref{sec:dstoch_opconex}. In the $t$-th iteration of \sopconex{} method with the search point $x^t$, it generates random vectors $\Ffrak^t:= \Ffrak(x^t, \xi^t), \gfrak(x^t, \xibar^t)$ and $\Gfrak(x^t, \xi^t)$ using SO in \eqref{eq:SO_F_operator-mo} -\eqref{eq:SO_g_convexity}.
We are now ready to describe the \sopconex{} method explicitly.
\vspace{-2mm}
\begin{algorithm}[H]
	\begin{algorithmic}[1]
		\State {\bf Input: }$x^{0} \in X, \lambda^0 = \zero$. 
		\State Set $x^{-1} = x^{0}$.
		\For{$t=0,1,2,\dots,T-1$}
		\State $\lambda^{t+1} \gets \argmin_{\lambda\ge 0} -\inprod{ \gfrak(x^t,\bar{\xi}^t)}{\lambda} + \tfrac{\tau_t}{2}\gnorm{\lambda-\lambda^t}{}{2}$
		\State $x^{t+1} \gets \argmin_{x \in X} \inprod{(1+\theta_t)\Ffrak^{t} - \theta_t\Ffrak^{t-1} +{\Gfrak(x^t,\xi^t)}\lambda^{t+1} }{x} + \tfrac{\eta_t}{2}\gnorm{x-x^t}{}{2}$
		\EndFor
		\State {\bf Output: }$\wb{x}_T := \paran{\tsum_{t=0}^{T-1} x^{t+1}}\big/T$
	\end{algorithmic}
	\caption{Partial Operator Extrapolation Method for fully-stochastic FCVI }\label{alg:algpdfullstoch}
\end{algorithm}
\vspace{-4mm}
We highlight some important differences between $\adopconex{}$ and $\sopconex{}$. First, unlike $\adopconex{}$, we do not extrapolate in $g$ (Line 4 of Algorithm \ref{alg:algpdfullstoch} is a special case of Line 5 in Algorithm \ref{alg:alg2} when we set $\theta_t = 0$). Second, we only apply extrapolation in the operator $F(\cdot)$, not in $ g'(\cdot)\lambda$ (Cf. Line 5 of Algorithm \ref{alg:algpdfullstoch} with Line 6 of Algorithm \ref{alg:alg2}). Effectively, we forego the extrapolation of terms of $\wt{F}$ that are generated due $\inprod{\lambda}{g(x)}$ which deals with the unbounded $\Rbb^m_+$ set. Third, unlike $x^{t+1}$ update in $\adopconex{}$ which can be done independently from $\lambda^{t+1}$, $\sopconex{}$ uses $\lambda^{t+1}$ in update of $x^{t+1}$ which is the most up-to-date information about $\lambda$.
In the next section, we present a stochastic version of the conjugate representation of $g$ which we will use to find a bound on the dual iterates $\{\lambda^t\}$ of \sopconex{} method. We now discuss the convergence guarantee of Algorithm~\ref{alg:algpdfullstoch}, and later, provide the convergence analysis.

\subsection{Convergence guarantees of \sopconex{} method}
We have the following convergence guarantee of Algorithm \ref{alg:algpdfullstoch}.
\begin{thm}\label{thm:feasib-fully-stoch}
	Suppose Algorithm~\ref{alg:algpdfullstoch} generates $\{x^{t+1},\lambda^{t+1}\}$ by setting  $\theta_t=1, \eta_t = \eta $ and $\tau_t = \tau$
	where
	\begin{equation}\label{eq:conditions}
	    \eta =\tfrac{\sqrt{2T}[B(3M_g+4\|\sigma_\Gfrak\|(1+\sigma_\gfrak))+H+\sigma]}{D_X} + \tfrac{25L}{3}, \quad \tau =\tfrac{2\sqrt{2T}(2M_g+5\|\sigma_\Gfrak\|+\sigma_\gfrak)D_X}{B}, 
	\end{equation}
    where $B\geq 1$ is a constant. Then, we have 
	\begin{subequations}\label{eq:feasib-opt-full}
		\begin{equation}\label{eq:opt-full}
			\begin{aligned}
				\Ebb[\sup_{x \in \wt{X}} \langle F(x),\bar{x}_T-x \rangle]
				& \leq D_X \Big[\tfrac{3\sqrt{2T}[B(3M_g+4\|\sigma_\Gfrak\|(1+\sigma_\gfrak))+H+\sigma]}{2\sqrt{T}} +\tfrac{\sigma_\gfrak^2B}{2\sqrt{2T}(2M_g+5\|\sigma_\Gfrak\|+\sigma_\gfrak)D_X^2}\\
				&\quad   + \tfrac{25(H^2+2\sigma^2)}{6\sqrt{2T}[B(3M_g+4\|\sigma_\Gfrak\|(1+\sigma_\gfrak))+H+\sigma]}
				 + \tfrac{H^2+2\sigma^2}{T\sqrt{2T}[B(3M_g+4\|\sigma_\Gfrak\|(1+\sigma_\gfrak))+H+\sigma]}\\
				 &\quad + \tfrac{\sigma^2}{2\sqrt{2T}[B(3M_g+4\|\sigma_\Gfrak\|(1+\sigma_\gfrak))+H+\sigma]}+ \tfrac{(15M_g^2+ 28\|\sigma_\Gfrak\|^2) (\|\lambda^*\|^2 + 2R_fe)}{3\sqrt{2T}[B(3M_g+4\|\sigma_\Gfrak\|(1+\sigma_\gfrak))+H+\sigma]} \Big]
				+ \tfrac{25LD_X^2}{2T} ,
			\end{aligned}
		\end{equation}
		\begin{equation}\label{eq:feasib-full}
			\begin{aligned}
				\Ebb\|g(\bar{x}_T)_+\|  &\leq \tfrac{B(3M_g+4\|\sigma_\Gfrak\|(1+\sigma_\gfrak))+H+\sigma}{\sqrt{2T}D_X} \|x^*-x^0\|^2 + \tfrac{25L}{6T} \|x^*-x^0\|^2+ \tfrac{\sqrt{2T}(2M_g+5\|\sigma_\Gfrak\|+\sigma_\gfrak)(\|\lambda^*\|+1)^2D_X}{\sqrt{T}B}\\
				&\quad  + \tfrac{\sigma_\gfrak^2B}{2\sqrt{2T}(2M_g+5\|\sigma_\Gfrak\|+\sigma_\gfrak)D_X} +  \tfrac{(H^2+2\sigma^2)D_X}{T\sqrt{2T}[B(3M_g+4\|\sigma_\Gfrak\|(1+\sigma_\gfrak))+H+\sigma]}+ \tfrac{15M_g^2D_X(\|\lambda^*\|+1)^2}{\sqrt{2T}[B(3M_g+4\|\sigma_\Gfrak\|(1+\sigma_\gfrak))+H+\sigma]}\\
				&\quad  +
				\tfrac{10M_g^2D_X(\|\lambda^*\|^2 + 2R_fe)}{\sqrt{2T}[B(3M_g+4\|\sigma_\Gfrak\|(1+\sigma_\gfrak))+H+\sigma]} +  \tfrac{25D_X[2(H^2+2\sigma^2) + \|\sigma_\Gfrak\|^2(\|\lambda^*\|^2 + 2R_fe)]}{3\sqrt{2T}[B(3M_g+4\|\sigma_\Gfrak\|(1+\sigma_\gfrak))+H+\sigma]} + \tfrac{(\|\lambda^*\|+1)\sigma_\gfrak}{\sqrt{T}},
			\end{aligned}
		\end{equation}
	\end{subequations}
 where 
 \[
   R_f :=\tfrac{\eta}{\tau}\|x^\ast-x^0\|^2+3\|\lambda^\ast\|^2 +2\sigma_\gfrak^2 + 9(H^2+2\sigma^2),
 \]
 and  $\sigma_\Gfrak := [\sigma_1, \dots, \sigma_m]^\top$.
\end{thm}

Theorem \ref{thm:feasib-fully-stoch} states the convergence result of \sopconex{} method for the fully-stochastic FCVI$(F, X, g_j, j \in [m])$. One can derive the convergence results for nonsmooth deterministic and stochastic (where we only use a stochastic approximation for the operator $F$) FCVI from the fully stochastic case. First, let us derive the convergence results for the nonsmooth deterministic case (i.e., $\sigma=\sigma_\gfrak = 0, \sigma_\Gfrak = \zero$). Thus the error convergence rate of $\sopconex{}$ after $T$ iterations is of 
\[
\mathcal{O}\big( \tfrac{LD_X^2}{T} + \tfrac{M_gD_X}{\sqrt{T}}(B + \tfrac{(\|\lambda^*\|+1)^2}{B}) + \tfrac{HD_X}{\sqrt{T}}\big) 
\]
Assuming $B\geq \|\lambda^*\|+1$ we have the following convergence rate for $\sopconex{}$
\begin{equation}\label{eq:deter}
    \mathcal{O}(\tfrac{B(M_g+H)D_X}{\sqrt{T}} + \tfrac{LD_X^2}{T} ).
\end{equation}
Note that the convergence in $M_g$ is of $\mathcal{O}(\tfrac{M_g}{\sqrt{T}})$ whether or not $B\geq \|\lambda^*\|+1$. We will show that this convergence rate dependence on $M_g$ can be improved significantly.

For the stochastic case where $\sigma >  0$ and $\sigma_\gfrak = 0, \sigma_\Gfrak = \zero$, the expected error convergence after $T$ iterations for $\sopconex{}$ is of 
\[
\mathcal{O}\big(\tfrac{LD_X^2}{T}  + \tfrac{M_gD_X}{\sqrt{T}}(B + \tfrac{(\|\lambda^*\|+1)^2}{B}) + \tfrac{(H+\sigma)D_X}{\sqrt{T}} \big)
\]
One can observe that the effect of the operator noise $(\delta^t)$ in each iteration $t\leq T$ of $\sopconex{}$ is the same as the nonsmoothness of the operator. Moreover, assuming $B\geq \|\lambda^*\|+1$, we balance the effect of $\gnorm{\lambda^*}{}{}$ and get similar results as we had in \eqref{eq:deter}. Regarding the fully stochastic case, we have the convergence rate of
\[
\mathcal{O}\big(\tfrac{LD_X^2}{T}  + \tfrac{M_gD_X}{\sqrt{T}}(B + \tfrac{(\|\sigma_\Gfrak\|+\sigma_\gfrak)(\|\lambda^*\|+1)^2}{B}) + \tfrac{(H+\sigma)D_X}{\sqrt{T}} + \tfrac{BD_X\|\sigma_\Gfrak\|\sigma_\gfrak}{\sqrt{T}} \big).
\]
Moreover, assuming that $B\geq \|\lambda^*\|+1$, 
we have the following overall expected error convergence rate after $T$ iterations of $\sopconex{}$
\[
\mathcal{O}\big(\tfrac{LD_X^2}{T}  + \tfrac{BM_gD_X}{\sqrt{T}}(\gnorm{\sigma_\Gfrak}{}{} + \sigma_\gfrak) + \tfrac{(H+\sigma)D_X}{\sqrt{T}} + \tfrac{BD_X\sigma_\gfrak\gnorm{\sigma_\Gfrak}{}{}}{\sqrt{T}} \big).
\]
The overall convergence rate is of $\mathcal{O}(\tfrac{1}{\sqrt{T}})$. We note that the dependence on $\sigma_\gfrak, \sigma_\Gfrak$ is of $\mathcal{O}(\tfrac{B\sigma_\gfrak\gnorm{\sigma_\Gfrak}{}{}}{\sqrt{T}})$ which can lead to slower convergence rate. We will significantly improve this dependence in Section \ref{sec:dstoch_opconex}.

\subsection{Convergence analysis of \sopconex{} method} 
To perform the convergence analysis of Algorithm \ref{alg:algpdfullstoch}, we need the framework of the conjugate duality of convex functions for the stochastic setup as presented below.
\subsubsection{Stochastic conjugate duality framework}
We present a stochastic conjugate duality (SCD) framework that utilizes the Fenchel conjugate representation of constraint functions $g_j, j \in [m]$. Denoting $g_j^*(v_j) : V_j \to \Rbb$ as the Fenchel conjugate of $g_j$, we define 
\begin{equation}\label{eq:cojugate}
	\Lcal_g(x,v) = v^\top x - g^*(v),
\end{equation} 
where $v = [ v_{1},\cdots, v_{m}] \in \mathbb{R}^{n\times m}$ and $g^*(v) = [g_1^*(v_1),  \dots, g_m^*(v_m)]^\top$. We also denote $\mathcal{V} = [V_{1},\cdots, V_{m}] \subset \mathbb{R}^{n\times m}$ as 
the domain of $g^*$. Then, it is clear that $v \in \Vcal$. From the definition of $g^*, \mathcal{V}$ and $\Lcal_g(x,v)$ and using the fact that $g$ is a convex continuous mapping, we have $g(x) = \sup_{v \in \mathcal{V}} \Lcal_g(x,v)$. Then, for a $\lambda \in \mathbb{R}^m_+$ , we define the following Lagrangian function 
\begin{equation}\label{eq:overall-SPP}
	\Lcal(x,\lambda,v) =  \langle F(x^\ast),x\rangle\ + \langle  \lambda,\Lcal_g(x,v)\rangle\ .
\end{equation} 
Proposition \ref{prop:FCVI_SPP_reformulation}, nonnegativity of the set $\Rbb^m_+$ and convexity of $g$ ensures the existence of $v^* \in \partial g(x^*)$ such that $(x^*,\lambda^*,v^*)$ is a solution of the problem $\min_{x\in X}\max_{\lambda \in \Rbb^m_+}\max_{v \in \mathcal{V}} \Lcal(x, \lambda,v)$, i.e., for any $ (x,\lambda,v) \in X \times \mathbb{R}^m_+\times \mathcal{V}$, we have 
\begin{equation*}
	\Lcal(x,\lambda^*,v^*)\geq  \Lcal(x^*,\lambda^*,v^*) \geq  \Lcal(x^*,\lambda,v).
\end{equation*}
Meanwhile, for $g$, we consider its Fenchel conjugate and let
	\begin{equation}\label{eq:Fenchel-Cons}
		v^{t+1} = \argmax_{v\in \mathcal{V}} \{\Lcal_g(x^{t},v)\}, \quad v^{t+1} \in \partial g(x^t),\quad g(x^t)  = \Lcal_g(x^{t},v^{t+1}).
	\end{equation} 
Using \eqref{eq:SO_g_convexity}, we can extend the conjugate duality framework for the stochastic setting. 
Let us define $ \vfrak_1^{t+1}$ and $\vfrak_2^{t+1}$ corresponding to $v^{t+1}$ as below
	\begin{equation*}\label{eq:subgrad_vt+1-1-stoch}
		\vfrak_1^{t+1} \ \in \partial \gfrak(x^t,\bar{\xi}^t),\quad \vfrak_2^{t+1} \in \partial\gfrak(x^t,\xi^t).
	\end{equation*}
Note that from \eqref{eq:SO_g_operator-mo}, we have $\Ebb[\vfrak_1^{t+1}]  = \Ebb[\vfrak_2^{t+1}] = v^{t+1}$. 
Given \eqref{eq:Fenchel-Cons} and the above definitions of $\vfrak^{t+1}_i,  i = 1,2$, we can restate updates in lines 4 and 5 of Algorithm \ref{alg:algpdfullstoch} as follows 
\begin{subequations}\label{eq:conj-update}
	\begin{equation}\label{eq:conj-update-lambda}
		\lambda^{t+1} = \argmin_{\lambda\ge 0} -\inprod{ \Lcal_g(x^t,\vfrak_1^{t+1})}{\lambda} + \tfrac{\tau_t}{2}\gnorm{\lambda-\lambda^t}{}{2},
	\end{equation}
	\begin{equation}\label{eq:conj-update-x}
		x^{t+1} = \argmin_{x \in X} \inprod{(1+\theta_t)\Ffrak^{t} - \theta_t\Ffrak^{t-1} +\vfrak_2^{t+1}\lambda^{t+1} }{x} + \tfrac{\eta_t}{2}\gnorm{x-x^t}{}{2}.
	\end{equation}
\end{subequations}
Let us denote $z^{t+1} := (x^{t+1},\lambda^{t+1},v^{t+1})$ and $ z  := (x,\lambda,v)  \in \wt{X} \times \mathbb{R}^m_+\times \mathcal{V}  $.
Then, we define $Q(z^{t+1}, z)$ as,
\begin{equation*}\label{eq:Gap-PD}
	Q(z^{t+1}, z)  :=  \langle F(x),x^{t+1}-x\rangle + \langle  \lambda,\Lcal_g(x^{t+1},v)\rangle - \langle  \lambda^{t+1},\Lcal_g(x,v^{t+1})\rangle. 
\end{equation*}
We are ready to present the convergence analysis for Algorithm \ref{alg:algpdfullstoch}. We first mention two useful propositions. 
\vspace{-2em}
\begin{prop}\label{prop:lemma2lan}
	(Lemma 2 of \cite{zhang2020optimal}). Let $\mathcal{F}_t$ be a filtration and $\{\delta_t\}_{t=1}^N$ be a martingale noise sequence such that $\delta_j \in \mathcal{F}_{t+1}$ for $j = 1, 2, \ldots, t$, $\mathbb{E}[\delta_t | \mathcal{F}_t] = 0$, and $\mathbb{E}[\|\delta_t\|^2 | \mathcal{F}_t] \leq \sigma^2$. For any random variable $\pi \in \Pi$ correlated with $\{\delta_t\}_{t=1}^N$, suppose it is bounded such that $\|\pi\| \leq M_\pi$ uniformly, then
	\[
	\mathbb{E}\left[\tsum_{t=1}^{N} \hat{\pi}^\top \delta_t\right] \leq \sqrt{N}M_\pi \sigma.
	\]	
\end{prop}
\vspace{-1.5em}
\begin{prop}\label{prop:tech_res1-mo}
	Let $\rho_0, \dots, \rho_j$ be a sequence of elements in $\Rbb^n$ and let $S$ be a convex set in $\Rbb^n$. Define the sequence $v_t, t = 0,1,\dots$, as follows: $v_0 \in S$ and 
	$
	v_{t+1} = \argmin_{x \in S} \inprod{\rho_t}{x} + \tfrac{1}{2}\gnorm{x-v_t}{}{2}.
	$ 
	Then, for any $x \in S$ and $t \ge 0$, the following inequalities hold
	\begin{equation}\label{eq:int_rel114_1-mo}
		\inprod{\rho_t}{v_t-x} \le \tfrac{1}{2}\gnorm{x-v_t}{}{2}-\tfrac{1}{2} \gnorm{x-v_{t+1}}{}{2}+ \tfrac{1}{2}\gnorm{\rho_t}{}{2},
	\end{equation}
\end{prop} 
\begin{Proof}
	Using Lemma \ref{lem:3-point-mo} with $g(x) =  \inprod{\rho_t}{x}$ 
	and $\mu = 0$, we have, due to the optimality of $v_{t+1}$,
	\begin{equation*}
		\inprod{\rho_t}{v_{t+1}-x} + \tfrac{1}{2}\gnorm{v_{t+1}-v_t}{2}{2} + \tfrac{1}{2} \gnorm{x-v_{t+1}}{}{2}\le  \tfrac{1}{2}\gnorm{x-v_t}{}{2},
	\end{equation*}
	is satisfied for all $x \in S$. The above relation and the fact
	$
	\inprod{\rho_t}{v_t-v_{t+1}} -\tfrac{1}{2}\gnorm{v_{t+1}-v_t}{}{2} \le \tfrac{1}{2}\gnorm{\rho_t}{}{2}$,
	imply \eqref{eq:int_rel114_1-mo}. 
	Hence, we conclude the proof.
\end{Proof}
Stochasticity of $\Ffrak, \gfrak$ and $\Gfrak$ gives rise to three types of errors below:
\begin{align*}
	\delta^t &:= \Ffrak^t - F(x^t)\\
	\delta_\gfrak^t &:= g(x^t, \xibar^t)  - g(x^t)\\
	\delta_\Gfrak^t &:= \tsum_{j \in [m]} \lambda^{t+1}_j[\Gfrak_j(x^t,\xi^t) - g'_j(x^t)]
\end{align*}
The following lemma bounds these errors.
%
\vspace{-1em}
\begin{lem}\label{lem:dstoch_useful_bounds-mo}
	For stochastic oracles satisfying \eqref{eq:SO_F_operator-mo} and \eqref{eq:SO_g_operator-mo}, we have 
	\begin{subequations}\label{eq:bd_var_dstoch-mo}
		\begin{align}
			\Ebb[\gnorm{\delta^{i} - \delta^{i-1}}{}{2}] 
			&\le 2\sigma^2, \label{eq:bd_var_F-mo}\\
			\Ebb[\gnorm{\delta_\Gfrak^i}{}{2}] 
			&\le \gnorm{\sigma_\Gfrak}{}{2}\Ebb[\gnorm{\lambda^{i+1}}{}{2}],
			\label{eq:bd_var_Gfrak-mo}\\
            \Ebb[\gnorm{\delta_\gfrak^i}{}{2}]&\leq \sigma_\gfrak^2
		\end{align}
	\end{subequations}
	where recall that $\sigma_\Gfrak = [\sigma_1, \dots, \sigma_m]^\top$.
	When $x,\lambda$ is fixed (i.e., non-random), we have
	\begin{equation}\label{eq:useful-expectation2-mo}
        \Ebb[\inprod{\delta^{i+1}}{x^{i+1}-x}]=0,\quad
		\Ebb[\inprod{\delta_\Gfrak^i}{x^{i}-x}] = \zero,\quad
        \Ebb[\inprod{\delta_\gfrak^i}{\lambda^i- \lambda}] = 0.
	\end{equation}
	Moreover, consider the sequences $\{x^i_a\}$, $\{\wt{x}^i_a\}$ as follows
	\begin{subequations}
		\begin{align}
			x^1_a = x^0, 
			\quad &	x^{i+1}_a := \argmin_{x \in X} -\inprod{\tfrac{\delta^i}{\eta_{i-1}}}{x} + \tfrac{1}{2}\gnorm{x-x^i_a}{}{2},
			\quad\forall i\geq 1,
			\label{eq:x-aug-seq_def-mo}\\
			\wt{x}^0_a = x^0, 
			\quad &\wt{x}^i_a := \argmin_{x\in X} -\inprod{\tfrac{\delta_\Gfrak^{i-1}}{\eta_{i-1}}}{x} 
			+ \tfrac{1}{2}\gnorm{x- \wt{x}^{i-1}_a}{}{2},
			\quad \forall  i\ge1.
			\label{eq:x-aug-seq_def2-mo}
		\end{align}
	\end{subequations}\vspace{-1em}
	Then, we have 
	\begin{subequations}
		\begin{align}
			\forall x\in X,
			\quad \tsum_{i=0}^t \gamma_{i}\inprod{\delta^{i+1}}{x-{x}^{i+1}_a} 
			&\le \tfrac{\gamma_0\eta_0}{2}\gnorm{x- x^0}{}{2} 
			+ \tsum_{i=0}^t \tfrac{\gamma_i}{2\eta_i}\gnorm{\delta^{i+1}}{}{2},
			\label{eq:stoch_innerprod_delta-mo}\\
   			\forall x\in X, \quad \hspace{1.4em}
			\tsum_{i=0}^t \gamma_i\inprod{\delta_\Gfrak^i}{x-\wt{x}^i_a} 
			&\le \tfrac{\gamma_0\eta_0}{2}\gnorm{x- x^0}{}{2} 
			+ \tsum_{i=0}^t \tfrac{\gamma_{i}}{2\eta_i}\gnorm{\delta_\Gfrak^i}{}{2}. \label{eq:stoch_innerprod_delta_G-mo}
		\end{align}
	\end{subequations}
\end{lem}
\begin{Proof}
	Denoting $\xi^{[i]} := (\xi^0, \dots, \xi^i)$ and borrowing a similar notation for $\xibar^{[i]}$, it is easy to see that $x^{i+1}$ depends on $(\xi^{[i]}, \xibar^{[i]})$, and as well as $\lambda^{i+1}$ depends on $(\xi^{[i-1]}, \xibar^{[i]})$. Then, we have
	\begin{align*}
             \Ebb[\delta_\gfrak^i] &= \Ebb[ [\Ebb[\gfrak(x^i, \xibar^t)|\xi^{[i-1]}, \xibar^{[i]}] -g(x^i)]] = \zero,\\
		\Ebb[\delta_\Gfrak^i] &= \Ebb\big[ \tsum_{j=1}^m\lambda^{i+1}_j [\Ebb[\Gfrak_j(x^i, \xi^i)|\xi^{[i-1]}, \xibar^{[i]}] - g'_j(x^i)]\big] = \zero,\\
		\Ebb[\delta^{i+1}] &= \Ebb[\Ffrak^{i+1} - F(x^{i+1} )] = \Ebb[ \Ebb[\Ffrak(x^{i+1}, \xi^{i+1})|\xi^{[i]},\xibar^{[i]}] -F(x^{i+1})] = \zero,
	\end{align*}
	the final equality in each of the above relations follows from the properties of stochastic oracles in \eqref{eq:SO_g_operator-mo}, and the rest of the equalities follow the tower law of the expectation. Using the above relations, we have
	\begin{align*}
		\Ebb[\inprod{\delta_\gfrak^i}{\lambda^{i}-\lambda}]&=\Ebb[\inprod{\Ebb[\delta_\gfrak^i|\xi^{[i-1]}, \xibar^{[i]}]}{\lambda^i-\lambda}] = \zero,\\
		\Ebb[\inprod{\delta_\Gfrak^i}{x^{i}-x}]&=\Ebb[\inprod{\Ebb[\delta_\Gfrak^i|\xi^{[i-1]}, \xibar^{[i]}]}{x^i-x}] = \zero,\\
		\Ebb[\inprod{\delta^{i+1}}{x^{i+1}-x}]&=\Ebb[\inprod{\Ebb[\delta^{i+1}|\xi^{[i]}, \xibar^{[i]}]}{x^{i+1}-x}] = \zero,
	\end{align*}
	for any non-random $x$. Hence, we prove \eqref{eq:useful-expectation2-mo}.
	
	Now, noting the definition of the sequence $\{\wt{x}^i_a\}_{i\ge 0}$ in \eqref{eq:x-aug-seq_def2-mo} 
	and applying Proposition \ref{prop:tech_res1-mo}, we have
	\[
	-\inprod{\tfrac{\delta_\Gfrak^i}{\eta_i}}{\wt{x}^i_a-x} 
	\le \tfrac{1}{2}\gnorm{x- \wt{x}^{i}_a}{}{2} 
	- \tfrac{1}{2}\gnorm{x- \wt{x}^{i+1}_a}{}{2} 
	+ \tfrac{1}{2\eta^2_i}\gnorm{\delta_\Gfrak^i}{}{2}.
	\]
	Multiplying the above relation by $\gamma_{i}\eta_i$, summing it from $i=0$ to $t$, given that $\gamma_{i}\eta_i \leq \gamma_{i-1}\eta_{i-1}$  we get \eqref{eq:stoch_innerprod_delta_G-mo}. 
	Also, considering the augmented sequence $\{x^i_a\}$ in \eqref{eq:x-aug-seq_def-mo} and using the similar development, we have \eqref{eq:stoch_innerprod_delta-mo}.
	
	We now prove an upper bound on $\Ebb[\gnorm{\delta_\Gfrak^i}{}{2}]$ as follows.
	\begin{align*}
		\Ebb[\gnorm{\delta_\Gfrak^i}{}{2}] 
		&= \Ebb[\gnorm{\tsum_{j=1}^m \lambda^{i+1}_j\delta^i_{\Gfrak_j}}{}{2}] 
		\le \Ebb\big[ \big(\tsum_{j=1}^m\gnorm{\lambda^{i+1}_j\delta^i_{\Gfrak_j}}{}{}\big)^2
		\big]
		\le \Ebb\bracket[\big]{\gnorm{\lambda^{i+1}}{}{2} \paran[\big]{\tsum_{j=1}^m \gnorm{\delta^i_{\Gfrak_j}}{}{2}} }
		\nonumber\\
		&= \Ebb\bracket[\big]{\gnorm{\lambda^{i+1}}{}{2} \paran[\big]{\tsum_{j=1}^m \Ebb[\gnorm{\delta^i_{\Gfrak_j}}{}{2}|\xi^{[i-1]},\xibar^{[i]}] } } 
		\le \gnorm{\sigma_\Gfrak}{}{2}\Ebb[\gnorm{\lambda^{i+1}}{}{2}],
	\end{align*}
	where the first inequality follows from the relation $\gnorm{\tsum_i a_i}{}{2} \le (\tsum_{i}\gnorm{a_i}{}{})^2$, the second inequality is a consequence of Cauchy-Schwarz inequality, and the last inequality follows from the definition of $\sigma_\Gfrak$.
	Finally, the bounds on $\Ebb[\gnorm{\delta_\gfrak^{i}}{}{2}]$ and $\Ebb[\gnorm{\delta^{i}- \delta^{i-1}}{}{2}]$ follows by the following relations 
	\begin{equation*}\label{eq:useful_expectations-mo}
		\begin{aligned}
            &\Ebb[\gnorm{\delta_\gfrak^i}{}{2}] = \Ebb[\Ebb[\gnorm{\delta_\gfrak^i}{}{2}|\xi^{[i-1]},\xibar^{[i]}]] \le \sigma_\gfrak^2,\\
			&\Ebb[\gnorm{\delta^i}{}{2}] = \Ebb[\Ebb[\gnorm{\delta^i}{}{2}|\xi^{[i-1]},\xibar^{[i]}]] \le \sigma^2,\\
			&\Ebb[\gnorm{\delta^i - \delta^{i-1}}{}{2}] 
			= \Ebb[\gnorm{\delta^i}{}{2}] + \Ebb[\gnorm{\delta^{i-1}}{}{2}] 
			+ \Ebb\inprod{\delta^i}{\Ebb[\delta^{i+1}|\xi_{[i]}]}
			= \Ebb[\gnorm{\delta^i}{}{2}] + \Ebb[\gnorm{\delta^{i-1}}{}{2}] 
			\le 2\sigma^2,
		\end{aligned}
	\end{equation*}
 where the first and second relation follows from \eqref{eq:SO_g_operator-mo} and \eqref{eq:SO_F_operator-mo}, respectively. We conclude the proof.
\end{Proof}
Now, we are ready to present an important convergence lemma. 
\begin{lem}\label{lemma:fully-Stoch-Partial}
	Suppose Algorithm \ref{alg:algpdfullstoch} generates $\{x^{i+1},\lambda^{i+1}, v^{i+1}\}$ by setting $\theta_i = 1$, $\eta_i=\eta\geq \tfrac{25L}{3}$, and $\tau_t = \tau$. Let $\lambda \in \mathbb{R}^m_+$ be a nonnegative bounded variable such that $\|\lambda\|\leq M_\lambda$. Then, for any integer $ t \leq T$, we have
	\begin{align}
		\tsum_{i=0}^{t-1}Q(z^{i+1}, z) + \tfrac{\tau}{2}\|\lambda-\lambda^{t}\|^2\leq & 	\tfrac{\eta}{2}\|x-x^0\|^2  +\tfrac{\tau}{2} \|\lambda-\lambda^0\|^2 +\tfrac{H^2+\| \delta^{t}-\delta^{t-1}\|^2}{\eta}+\tsum_{i=0}^{t-1} \tfrac{25(H^2 +\| \delta^{i}-\delta^{i-1}\|^2)}{6\eta} \nonumber\\ 
		&\quad+\tsum_{i=0}^{t-1}(\tfrac{10M_\lambda^2M_g^2}{\eta}+\tfrac{5M_g^2\|\lambda-\lambda^{i+1}\|^2}{2\eta})+\tsum_{i=0}^{t-1}\tfrac{\|\delta_\gfrak^i\|^2}{\tau} + \tfrac{25\|\delta_\Gfrak^i\|^2}{6\eta} \nonumber\\
		&\quad + \tsum_{i=0}^{t-1}\inprod{\delta^{i+1}}{x-x^{i+1}} + \inprod{\delta_\gfrak^i}{\lambda-\lambda^i} - \inner{\delta_\Gfrak^i}{x^i-x},\label{eq:fully-stoch-partial}
	\end{align}
	where recall that 
	$\delta_\Gfrak^i = \tsum_{j=1}^m \lambda^{i+1}_j [\Gfrak_j(x^i, \xi^i) -  g'_j(x^i)]$ and 
	$\delta_\gfrak^i = \gfrak(x^i,\bar{\xi}^i) -g(x^i)$.
\end{lem}
\begin{Proof}
Note that $Q$ can be written as 
	\begin{equation}\label{eq:gap-conj}
		\begin{aligned}
			Q(z^{i+1}, z)& = \langle F(x),x^{i+1}-x\rangle + \langle  \lambda,\Lcal_g(x^{i+1},v)\rangle - \langle  \lambda^{i+1},\Lcal_g(x,v^{i+1})\rangle\\
			& = \underbrace{ \langle F(x),x^{i+1}-x\rangle\ +\langle  \lambda^{i+1},\Lcal_g(x^{i+1},v^{i+1})-\Lcal_g(x,v^{i+1})\rangle}_{T_1}\\
			&\quad + \underbrace{\langle  \lambda,\Lcal_g(x^{i+1},v)-\Lcal_g(x^{i+1},v^{i+1})\rangle}_{T_2}+ \underbrace{\inner{\lambda-\lambda^{i+1}}{\Lcal_g(x^{i+1},v^{i+1})}}_{T_3}.
		\end{aligned}
	\end{equation}
For $T_2$ in \eqref{eq:gap-conj} and in view of \eqref{eq:cojugate}, we have 
	\begin{align}
		\langle  \lambda,\Lcal_g(x^{i+1},v)-\Lcal_g(x^{i+1},v^{i+1})\rangle  &= \langle  \lambda,\inner{x^{i+1}}{v-v^{i+1}}-(g^*(v)-g^*(v^{i+1}))\rangle  \nonumber\\
		& = \langle  \lambda,(v-v^{i+1})^\top(x^{i+1}-x^i)-(g^*(v)-g^*(v^{i+1})- (v-v^{i+1})^\top x^i)\rangle \nonumber\\
		& \leq  \langle  \lambda,(v-v^{i+1})^\top(x^{i+1}-x^i)\rangle \nonumber\\
		& \leq \tfrac{5\|v-v^{i+1}\|^2\|\lambda\|^2}{2\eta_i} + \tfrac{\eta_i\|x^{i+1}-x^i\|^2}{10} 
		\leq \tfrac{10M_\lambda^2M_g^2}{\eta_i} + \tfrac{\eta_i\|x^{i+1}-x^i\|^2}{10},\label{eq:firstterm}
	\end{align}
	where the last relation follows from $ \|v-v^{i+1}\|^2 \leq 2\|v\|^2 + 2\|v^{i+1}\|^2\leq 4M_g^2$ and $\|\lambda\|\leq M_\lambda$. 
\\
Next, $T_3$ in \eqref{eq:gap-conj} will be decomposed further as below
	\begin{align}
			&(\lambda-\lambda^{i+1})^\top\Lcal_g(x^{i+1},v^{i+1}) \nonumber\\
			& = (\lambda-\lambda^{i+1})^\top(\Lcal_g(x^{i+1},v^{i+1}) - \Lcal_g(x^i,v^{i+1})) 
			+ (\lambda-\lambda^{i+1})^\top (g(x^i) - \gfrak(x^i,\bar{\xi}^i)) + (\lambda-\lambda^{i+1})^\top \gfrak(x^i,\bar{\xi}^i) \nonumber \\
			& = \underbrace{\langle v^{i+1}(\lambda - \lambda^{i+1}), x^{i+1}-x^i\rangle}_{ T_{3.1}} +\underbrace{(\lambda-\lambda^{i+1})^\top  (g(x^i) - \gfrak(x^i,\bar{\xi}^i))}_{ T_{3.2}}
			+ \underbrace{ (\lambda-\lambda^{i+1})^\top \gfrak(x^i,\bar{\xi}^i)}_{ T_{3.3}}. \label{eq:second-term-full-stoch}
	\end{align}
	From Young's inequality, we can write $T_{3.1}$ in \eqref{eq:second-term-full-stoch} as follows 
	\begin{equation*}\label{eq:secondterm1-full-stoch}
		\langle  v^{i+1}(\lambda - \lambda^{i+1}), x^{i+1}-x^i\rangle \leq \tfrac{5M_g^2\|\lambda-\lambda^{i+1}\|^2}{2\eta_i} + \tfrac{\eta_i\|x^{i+1}-x^i\|^2}{10}.
	\end{equation*}
	 Considering $T_{3.2}$ in \eqref{eq:second-term-full-stoch}, and from the definition of $ \delta_\gfrak^i $, we have
	\begin{equation*}\label{eq:3.2stoch1}
		\begin{aligned}
			(\lambda-\lambda^{i+1})^\top(g(x^i)-\gfrak(x^i,\bar{\xi}^i))& = (\lambda-\lambda^{i})^\top\delta_\gfrak^i + (\lambda^i-\lambda^{i+1})^\top\delta_\gfrak^i.
		\end{aligned}
	\end{equation*}
	Moreover, 
	\begin{equation*}\label{eq:3.2stoch2}
		\begin{aligned}
			(\lambda^i-\lambda^{i+1})^\top\delta_\gfrak^i &\leq \tfrac{\tau\|\lambda^i-\lambda^{i+1}\|^2}{4} + \tfrac{\|\delta_\gfrak^i\|^2}{\tau}.
		\end{aligned}
	\end{equation*}
	Now, consider $T_{3.3}$ in \eqref{eq:second-term-full-stoch}. In view of \eqref{eq:conj-update-lambda}, from the optimality of $\lambda^{i+1}$ and $\tau_i = \tau$, we have 
	\begin{equation*}\label{eq:lambdaoptfullstoch}
		\langle \lambda-\lambda^{i+1},\gfrak(x^i,\bar{\xi}^i)\rangle \leq \tfrac{\tau}{2} [\|\lambda-\lambda^i\|^2  -\|\lambda-\lambda^{i+1}\|^2 - \|\lambda^i-\lambda^{i+1}\|^2] . 
	\end{equation*}
	Next, let us handle $T_1$ in \eqref{eq:gap-conj}, we have 
	\begin{equation*}
		\begin{aligned}
			\langle F(x),x^{i+1}-x\rangle\ +\langle  \lambda^{i+1},\Lcal_g(x^{i+1},v^{i+1})-\Lcal_g(x,v^{i+1})\rangle &= \langle F(x),x^{i+1}-x\rangle\ +\langle  v^{i+1}\lambda^{i+1} , x^{i+1}- x \rangle\\
			&  = \langle F(x)+ v^{i+1}\lambda^{i+1},x^{i+1}-x\rangle.
		\end{aligned}
	\end{equation*}
	In view of \eqref{eq:conj-update-x} and from optimality of $x^{i+1}$ we have 
	\begin{equation}\label{eq:xoptfullstoch}
		\begin{aligned}
			\langle \Ffrak^{i+1} + v^{i+1}\lambda^{i+1} , x^{i+1}-x \rangle &\leq  \tfrac{\eta_i}{2}[\|x-x^i\|^2 - \|x-x^{i+1}\|^2-\|x^i-x^{i+1}\|^2] - \langle \delta_\Gfrak^{i},x^{i+1}-x \rangle\\
			&\quad  + \langle \Delta \Ffrak_{i+1},x^{i+1}-x \rangle  - \theta_i \langle \Delta \Ffrak_{i},x^{i}-x \rangle - \theta_i \langle \Delta \Ffrak_{i},x^{i+1}-x^i \rangle.
		\end{aligned}
	\end{equation}
	Note that 
\begin{equation}\label{eq:monotone-full-stoch}
	\begin{aligned}
		\inprod{\Delta \Ffrak_t}{x^{i+1}-x^i} & =  \inprod{\Delta F_t}{x^{i+1}-x^i} + \inprod{\delta^{i}- \delta^{i-1}}{x^{i+1}-x^i},\\
		\inprod{\Ffrak^{i+1}}{x^{i+1}-x} & = \inprod{F(x^{i+1})+\delta^{i+1}}{x^{i+1}-x} + 
		\ge \inprod{F(x)}{x^{i+1}-x} + \inprod{\delta^{i+1}}{x^{i+1}-x}.
	\end{aligned}
\end{equation}
	Hence, from \eqref{eq:xoptfullstoch}-\eqref{eq:monotone-full-stoch} one can conclude 
	\begin{equation*}\label{eq:xopt4fullstoch}
		\begin{aligned}
			&\langle F(x) + \lambda^{i+1} v^{i+1}, x^{i+1}-x \rangle + \langle \delta^{i+1},x^{i+1}-x \rangle\\
            & \leq 
			\tfrac{\eta_i}{2}\big[\|x-x^i\|^2 - \|x-x^{i+1}\|^2-\|x^i-x^{i+1}\|^2\big] + \langle \Delta \Ffrak_{i+1},x^{i+1}-x \rangle  - \theta_t \langle \Delta \Ffrak_{i},x^{i}-x \rangle\\
			&\quad  - \theta_i [\langle \Delta F_{i},x^{i+1}-x^i \rangle + \langle \delta^i-\delta^{i-1},x^{i+1}-x^i \rangle] - \langle \delta^{i+1}+\delta_\Gfrak^{i},x^{i+1}-x \rangle.
		\end{aligned}
	\end{equation*}
	From the bounds we had for $T_1$, $T_2$, and $T_3$ in \eqref{eq:gap-conj}, letting $\theta_i =1,\eta_i=\eta, \tau_i=\tau$, we have 
	\begin{equation}\label{eq:fullstoch1gap}
		\begin{aligned}
			Q(z^{i+1},z)&\leq 	\tfrac{\eta}{2}[\|x-x^i\|^2 - \|x-x^{i+1}\|^2] + \tfrac{\tau}{2}[ \|\lambda-\lambda^i\|^2  -\|\lambda-\lambda^{i+1}\|^2]\\
			&\quad +  \langle \Delta \Ffrak_{i+1},x^{i+1}-x \rangle  -  \langle \Delta \Ffrak_{i},x^{i}-x \rangle -[ \langle \Delta F_{i},x^{i+1}-x^i \rangle + \langle \delta^i-\delta^{i-1},x^{i+1}-x^i \rangle]\\
			&\quad  - \langle \delta^{i+1}+\delta_\Gfrak^{i},x^{i+1}-x \rangle  - \tfrac{3\eta}{10}\|x^{i+1}-x^i\|^2 - \tfrac{\tau}{4}\|\lambda^i-\lambda^{i+1}\|^2\\
			&\quad +\tfrac{10M_g^2M_\lambda^2}{\eta}+ \tfrac{5M_g^2\|\lambda-\lambda^{i+1}\|^2}{2\eta} + \inner{\lambda-\lambda^i}{\delta_\gfrak^i} + \tfrac{\|\delta_\gfrak^i\|^2}{\tau} .
		\end{aligned}
	\end{equation} 
	Note that 
		\begin{equation}\label{eq:cauchy-full}
		\langle \Delta \Ffrak_{t},x^{t}-x \rangle \leq L \|x^t-x^{t-1}\|\|x^{t}-x\| + H\|x^{t}-x\| + \langle \delta^t-\delta^{t-1},x^{t}-x \rangle
	\end{equation}
	Thus in view of \eqref{eq:cauchy-full}, rearranging $- \tsum_{i=0}^{t-1}\tfrac{3\eta}{10}\|x^{i+1}-x^i\|^2$, and summing \eqref{eq:fullstoch1gap} from $i=0$ to $i=t-1$ we have 
	\begin{align}
			\tsum_{i=0}^{t-1}Q(z^{i+1},z)&\leq 	\tfrac{\eta}{2}[\|x-x^0\|^2 - \|x-x^{t}\|^2] + \tfrac{\tau}{2}[ \|\lambda-\lambda^0\|^2  -\|\lambda-\lambda^{t}\|^2]+  L \|x^{t}-x^{t-1}\|\|x^{t}-x\| + H \|x^{t}-x\|  \nonumber\\
			& \quad +\tsum_{i=0}^{t-1}[ \tfrac{10M_g^2M_\lambda^2}{\eta}+ \tfrac{5M_g^2\|\lambda-\lambda^{i+1}\|^2}{2\eta}]+\tsum_{i=0}^{t-1}\langle \delta_\Gfrak^i,x^{i}-x^{i+1} \rangle - 	\tfrac{3\eta}{50}\|x^{i+1}-x^i\|^2\nonumber\\
			&\quad+ \tsum_{i=0}^{t-1} L \|x^i-x^{i-1}\|\|x^{i+1}-x^i\| - \tfrac{3\eta}{50}\|x^{i+1}-x^i\|^2 -\tfrac{3\eta}{50}\|x^{i}-x^{i-1}\|^2 \nonumber\\
			& \quad +\tsum_{i=0}^{t-1} [H \|x^{i+1}-x^i\|  - \tfrac{3\eta}{50}\|x^{i+1}-x^i\|^2] - \tfrac{3\eta}{50}\|x^{t}-x^{t-1}\|^2 + \langle 	\delta^{t}-\delta^{t-1},x^{t}-x \rangle \nonumber\\
			&\quad +\tsum_{i=0}^{t-1}\langle \delta^i-\delta^{i-1},x^{i+1}-x^i \rangle - 	\tfrac{3\eta}{50}\|x^{i+1}-x^i\|^2 \nonumber\\
            &\quad  + \tsum_{i=0}^{t-1} \langle \delta_\gfrak^i,\lambda-\lambda^i \rangle-\langle 	\delta_\Gfrak^i,x^{i}-x \rangle - \langle \delta^{i+1},x^{i+1}-x \rangle + \tsum_{i=0}^{t-1}\tfrac{\|\delta_\gfrak^i\|^2}{\tau}.\label{eq:fullstoch3gap}
	\end{align} 
	Now, using $ \tfrac{25L}{3}\leq \eta$ and from Young's inequality, 
    we have
	\begin{equation}\label{eq:relations}
		\begin{aligned}
            L \|x^i-x^{i-1}\|\|x^{i+1}-x^i\| - \tfrac{3\eta}{50}\|x^{i+1}-x^i\|^2 -\tfrac{3\eta}{50}\|x^{i}-x^{i-1}\|^2 &\le 0,\\
			H \|x^{i+1}-x^i\|  - \tfrac{3\eta}{50}\|x^{i+1}-x^i\|^2 &\leq  \tfrac{25H^2}{6\eta} ,\\
			\langle \delta^i-\delta^{i-1},x^{i+1}-x^i \rangle - \tfrac{3\eta}{50}\|x^{i+1}-x^i\|^2 &\leq\tfrac{25}{6\eta}\|  \delta^i-\delta^{i-1}\|^2,\\
			\langle \delta_\Gfrak^i,x^{i}-x^{i+1} \rangle - \tfrac{3\eta}{50}\|x^{i+1}-x^i\|^2 &\leq \tfrac{25\|\delta_\Gfrak^i\|^2}{6\eta},\\
			L \|x^{t}-x^{t-1}\|\|x^{t}-x\|  -\tfrac{3\eta}{50}\|x^{t}-x^{t-1}\|^2 - \tfrac{\eta}{8}\|x-x^{t}\|^2 &\leq 0, \\
			\langle \delta^{t}-\delta^{t-1},x^{t}-x \rangle  - \tfrac{\eta}{4}\|x-x^{t}\|^2 &\leq \tfrac{\| \delta^{t}-\delta^{t-1}\|^2}{\eta},\\
			H\|x^{t}-x\|  - \tfrac{\eta}{8}\|x-x^{t}\|^2 &\leq \tfrac{H^2}{\eta}.
        \end{aligned}
	\end{equation}
	Using the relations in \eqref{eq:relations} inside \eqref{eq:fullstoch3gap}, we get \eqref{eq:fully-stoch-partial}. Hence, we conclude the proof.
\end{Proof}
\vspace{-2em}
\subsubsection{Boundedness of the Lagrange multiplier
}
Our next goal is to show that $\lambda^t$ is bounded in expectation. 
Before formally stating this lemma, we mention the following technical proposition. 
\begin{prop}[Lemma 5 of \cite{boob2023stochastic}]\label{prop:recursive_bound-mo}
	Let $\{a_t\}_{t\ge 0}$ be a nonnegative sequence, $m_1,m_2 \ge 0$ be constants such that $a_0 \le m_1$ and the following relation holds for all $t \ge 1$:
	\[a_t \le m_1  + m_2\tsum_{ k =0}^{t-1}a_k.\]
	Then, we have $a_t \le m_1 (1+m_2)^t$.
\end{prop}
\vspace{-1em}
\begin{lem}\label{lemma:boundedness3}
	Let $\{x^t,\lambda^t,v^t\}$ be the sequence generated by Algorithm \ref{alg:algpdfullstoch} with parameters $ \theta_t=1, \eta_t = \eta$ and $\tau_i =\tau$ such that $\eta\tau \geq 2(5M_g^2+17\gnorm{\sigma_\Gfrak}{}{2})T$ and $\tau = \Theta(\sqrt{T})$ for all $t \leq T$. Then, for any $t \le T$, we have 
	\begin{equation*}
		\mathbb{E}[\|\lambda^*-\lambda^{t}\|^2]\leq 2 R_fe,
	\end{equation*}
	where 
	\[
	R_f =\tfrac{\eta}{\tau}\|x^\ast-x^0\|^2+3\|\lambda^\ast\|^2 +2\sigma_\gfrak^2 + 9(H^2+2\sigma^2).
	\]
\end{lem}
\begin{Proof}
According to \eqref{eq:useful-expectation2-mo} for fixed (non-random) $x = x^*$ and $\lambda = \lambda^*$, we have 
	\begin{equation}\label{eq:expectation-full}
		\begin{aligned}
			\mathbb{E}[\langle \delta_\gfrak^i,\lambda^{i}-\lambda^* \rangle] &= 0, \\
			\mathbb{E}[\langle \delta_\Gfrak^i,x^{i}-x^*\rangle] &= \zero, \\
			\mathbb{E}[\langle \delta^{i+1},x^{i+1}-x^* \rangle] &= 0 .
		\end{aligned}
	\end{equation}
	We now apply Lemma \ref{lemma:fully-Stoch-Partial} for  $ z^\ast = (x^\ast,\lambda^\ast,v^\ast)$, set $M_{\lambda^\ast} = \|\lambda^\ast\|$, observe that $Q(z^{i+1},z^*)\geq 0 $ can be dropped from the LHS, then multiply both sides of the relation by $\tfrac{2}{\tau}$, take expectation and apply bounds in \eqref{eq:bd_var_dstoch-mo} and \eqref{eq:expectation-full} along with the relation
	\[
	\mathbb{E}[\|\delta_\Gfrak^i\|^2] \leq \gnorm{\sigma_\Gfrak}{}{2}\mathbb{E}[\|\lambda^{i+1}\|^2] \leq 2 \gnorm{\sigma_\Gfrak}{}{2} [\|\lambda^\ast\|^2 + \mathbb{E}[\|\lambda^{i+1}- \lambda^\ast\|^2] ],
	\] 
	we obtain, 
	\begin{align*}\label{eq:boundfullstoch1}
			\mathbb{E}[\|\lambda^\ast-\lambda^{t}\|^2]&\leq \tfrac{\eta}{\tau}\|x^\ast-x^0\|^2+ \|\lambda^\ast-\lambda^{0}\|^2 + \tfrac{2t\sigma_\gfrak^2}{\tau^2}+  \tfrac{25t(H^2+2\sigma^2)}{3\eta\tau} + \tfrac{2(H^2+2\sigma^2)}{\eta\tau}\\
			&\quad + \tsum_{i=0}^{t-1}[ \tfrac{2(10M_g^2 + \tfrac{25\gnorm{\sigma_\Gfrak}{}{2}}{3})\|\lambda^\ast\|^2}{\eta\tau}+ \tfrac{(5M_g^2+\tfrac{50\gnorm{\sigma_\Gfrak}{}{2}}{3})\mathbb{E}[\|\lambda^\ast-\lambda^{i+1}\|^2]}{\eta\tau}].
	\end{align*}
\vspace{-0.1em}
	Hence, we have
	\begin{equation*}\label{eq:boundfullstoch3}
		\begin{aligned}
			(1-\tfrac{5M_g^2+\tfrac{50\gnorm{\sigma_\Gfrak}{}{2}}{3}}{\eta\tau})&\mathbb{E}[\|\lambda^\ast-\lambda^{t}\|^2]\leq \tfrac{\eta}{\tau}\|x^\ast-x^0\|^2 + \|\lambda^\ast-\lambda^{0}\|^2+ \tfrac{2t\sigma_\gfrak^2}{\tau^2}+  \tfrac{25t(H^2+2\sigma^2)}{3\eta\tau} \\
			&\quad + \tsum_{i=0}^{t-1} \tfrac{2(10M_g^2 + \tfrac{25\gnorm{\sigma_\Gfrak}{}{2}}{3})\|\lambda^\ast\|^2}{\eta\tau}+ \tsum_{i=0}^{t-2}\tfrac{(5M_g^2+\tfrac{50\gnorm{\sigma_\Gfrak}{}{2}}{3})\mathbb{E}[\|\lambda^\ast-\lambda^{i+1}\|^2]}{\eta\tau}+ \tfrac{2(H^2+2\sigma^2)}{\eta\tau}.
		\end{aligned}
	\end{equation*}
	Note that $\lambda^0 = \zero$, and define $R^t_f$ as    
	\[
	R^t_f  := \tfrac{\eta}{\tau}\|x^\ast-x^0\|^2 + \|\lambda^\ast\|^2+ \tfrac{2t\sigma_\gfrak^2}{\tau^2}+  \tfrac{(25t+6)(H^2+2\sigma^2)}{3\eta\tau} +  \tfrac{2t(10M_g^2 + \tfrac{25\gnorm{\sigma_\Gfrak}{}{2}}{3})\|\lambda^\ast\|^2}{\eta\tau}.
	\]
	Note that since $\eta\tau \geq 2(5M_g^2+17\gnorm{\sigma_\Gfrak}{}{2})T$ and $\tau = \Theta(\sqrt{T})$, we have $ R^t_f \leq R_f$.
	Thus 
	\begin{equation*}\label{eq:bounded4fullstoch}
		(1-\tfrac{15M_g^2+50\gnorm{\sigma_\Gfrak}{}{2}}{3\eta\tau})\mathbb{E}[\|\lambda^\ast-\lambda^{t}\|^2] \leq R_f + \tfrac{15M_g^2+50\gnorm{\sigma_\Gfrak}{}{2}}{3\eta\tau}\tsum_{i=0}^{i=t-2}\mathbb{E}[\|\lambda^\ast-\lambda^{i+1}\|^2].
	\end{equation*}
	Note that $	(1-\tfrac{15M_g^2+50\gnorm{\sigma_\Gfrak}{}{2}}{3\eta\tau}) \geq 1/2$, thus
	\begin{equation*}\label{eq:bounded6}
		\Ebb[\|\lambda^{t} - \lambda^\ast\|^2] \leq 2[R_f  + \tfrac{15M_g^2+50\gnorm{\sigma_\Gfrak}{}{2}}{3\eta\tau}\tsum_{i=0}^{i=t-2}\mathbb{E}[\|\lambda^\ast-\lambda^{i+1}\|^2]\leq 2R_f + \tfrac{1}{T}\tsum_{i=0}^{t-2} \|\lambda^*-\lambda^{i+1}\|^2. 
	\end{equation*}
	From Proposition \ref{prop:recursive_bound-mo}, we have 
	\begin{equation*}\label{eq:bounded8}
		\Ebb[\|\lambda^{t} - \lambda^\ast\|^2] \leq 2R_f (1+\tfrac{1}{T})^t \leq 2R_fe, \ \forall t \le T.
	\end{equation*} 
    Hence,  we conclude the proof.
\end{Proof}
\noindent Notice that this bound holds deterministically for nonsmooth deterministic FCVI. We are ready to present the proof of Theorem \ref{thm:feasib-fully-stoch}.\\
{\bf Proof of Theorem \ref{thm:feasib-fully-stoch}.}\\
First, note that $\eta\tau \geq 2(5M_g^2+17\gnorm{\sigma_\Gfrak}{}{2})T$ which is a requirement for bounding $\{\lambda^t\}$ in Lemma \ref{lemma:boundedness3}. Second, by setting $ \hat{z} = (x,\zero,\hat{v})$  where $x \in \wt{X}$ and $\bar{z}_T = (\bar{x}_T,\bar{\lambda}_T,\bar{v}_T)$ where $\hat{v}\in \partial g(\bar{x}_T)$ and using concavity of $Q$ in $v$ and $\lambda$, we have 
	\begin{equation}\label{eq:opt1}
			\tfrac{1}{T}\tsum_{t=0}^{T-1}Q(z^{t+1}, \hat{z})\geq Q(\bar{z}_T,\hat{z})
			 = \langle F(x),\bar{x}_T-x \rangle - \inner{\bar{\lambda}_T}{\Lcal_g(x,\bar{v}_T)}\geq \langle F(x),\bar{x}_T-x \rangle,
	\end{equation}
	where the last inequality in \eqref{eq:opt1} follows from the fact that $\bar{\lambda}_T\geq 0$, $ \Lcal_g(x,\bar{v}_T) \leq g(x)$ and $g(x)\leq 0 $ for all $x\in \tilde{X}$.  
	Therefore, using Lemma \ref{lemma:fully-Stoch-Partial} and taking expectation on both sides, we have 
 \begin{equation*}\label{eq:opt-full-1}
		\begin{aligned}
			\Ebb[\sup_{x \in \wt{X}} \langle F(x),\bar{x}_T-x \rangle] \leq & \tfrac{1}{T}\big[\tfrac{\eta}{2}\Ebb[\|x-x^0\|^2] + \tsum_{t =0}^{T-1}\tfrac{5M_g^2\Ebb[\|\lambda^{t+1}\|^2]}{2\eta}+ \tfrac{25T(H^2 + 2\sigma^2)}{6\eta} + \tsum_{t =0}^{T-1} \tfrac{25\|\sigma_\Gfrak\|^2\Ebb[\|\lambda^{t+1}\|^2]}{6\eta}\\
			&+\tfrac{H^2+2\sigma^2}{\eta} +  \tfrac{T\sigma_\gfrak^2}{\tau} +\tsum_{t=0}^{T-1}\Ebb[\inprod{\delta^{t+1}}{x-x^{t+1}} + \inprod{\delta_\gfrak^t}{\lambda-\lambda^t} - \inner{\delta_\Gfrak^t}{x^t-x}] \big].
		\end{aligned}
	\end{equation*}
	In view of \eqref{eq:bd_var_Gfrak-mo}, 
    \eqref{eq:stoch_innerprod_delta-mo}, \eqref{eq:stoch_innerprod_delta_G-mo} and the last relation in \eqref{eq:useful-expectation2-mo} with $\lambda = \zero$
    , we have 
	\begin{equation}\label{eq:opt-full-2}
		\begin{aligned}
			\Ebb&[\sup_{x \in \wt{X}} \langle F(x),\bar{x}_T-x \rangle]\leq  \tfrac{1}{T}\big[\tfrac{3\eta}{2}\Ebb[\|x-x^0\|^2] + \tsum_{t =0}^{T-1}\tfrac{5M_g^2\Ebb[\|\lambda^{t+1}\|^2]}{2\eta}+ \tfrac{25T(H^2 + 2\sigma^2)}{6\eta} + \tsum_{t =0}^{T-1} \tfrac{25\|\sigma_\Gfrak\|^2\Ebb[\|\lambda^{t+1}\|^2]}{6\eta}\\
			& + \tfrac{T\sigma_\gfrak^2}{\tau}+\tfrac{H^2+2\sigma^2}{\eta} - \tsum_{t=0}^{T-1}\Ebb[\inprod{\delta^{t+1}}{x^{t+1}-x_a^{t+1}} + \inner{\delta_\Gfrak^t}{x^t-\wt{x}^t_a}] + \tsum_{t=0}^{T-1}\Ebb[ \tfrac{\|\delta_\Gfrak^t\|^2 + \|\delta^{t+1}\|^2}{2\eta}]  \big].
		\end{aligned}
	\end{equation}
 Moreover, note that from Lemma \ref{lemma:boundedness3} we have
	\begin{equation}\label{eq:bound}
		\Ebb[\|\lambda^{t+1}\|^2] \leq 2\|\lambda^\ast\|^2 + 2\Ebb[\|\lambda^{t+1}-\lambda^\ast\|^2]\leq 2\|\lambda^\ast\|^2 + 4R_fe.
	\end{equation}
 Thus, given  \eqref{eq:bound}, \eqref{eq:bd_var_dstoch-mo} and \eqref{eq:conditions} along with $\Ebb[\inprod{\delta^{t+1}}{x^{t+1}-x_a^{t+1}} + \inner{\delta_\Gfrak^t}{x^t-\wt{x}^t_a}]=0$, we obtain \eqref{eq:opt-full}. 
 
 Now, we prove the feasibility gap. Let $\hat{v}$ be the dual variable associated with $\bar{x}_T$ such that 
	\begin{equation}\label{eq:feasib1}
			\hat{v} \in \arg \max_{v\in V} \Lcal_g(\bar{x}_T,v) \quad \text{and}\quad  g(\bar{x}_T) =  \Lcal_g(\bar{x}_T,\hat{v})\geq  \Lcal_g(\bar{x}_T,v^\ast).
	\end{equation}
	Since $\lambda^\ast\geq 0$ and by \eqref{eq:feasib1} and the overall Lagrangian Function in \eqref{eq:overall-SPP}, one can see 
	\begin{equation*}
		\begin{aligned}
			\Lcal(\bar{x}_T,\lambda^\ast,\hat{v}) = & \langle F(x^\ast),\bar{x}_T\rangle + \langle \lambda^\ast,\Lcal_g(\bar{x}_T,\hat{v}) \rangle\\
			\geq&  \langle F(x^\ast),\bar{x}_T\rangle + \langle \lambda^\ast,\Lcal_g(\bar{x}_T,v^\ast) \rangle\\
			= &\Lcal(\bar{x}_T,\lambda^\ast,v^\ast)\geq  \Lcal(x^\ast,\lambda^\ast,v^\ast) = \inprod{F(x^\ast)}{x^\ast}.
		\end{aligned}
	\end{equation*}
	Hence, we get
	\[
	\langle F(x^\ast),\bar{x}_T\rangle  + \langle \lambda^\ast,g(\bar{x}_T) \rangle - \langle F(x^\ast),x^\ast\rangle\geq 0.
	\]
	Meanwhile, due to the fact that $\lambda^\ast\geq 0$ and $g(\bar{x}_T)\leq g(\bar{x}_T)_+$, we have 
	$
	\langle \lambda^\ast,g(\bar{x}_T) \rangle \leq  \langle \lambda^\ast,g(\bar{x}_T)_+ \rangle 
	$
	which leads to 
	\begin{equation}\label{eq:feasib1a}
		\langle F(x^\ast),\bar{x}_T\rangle + \|\lambda^\ast\|\|g(\bar{x}_T)_+\|-\langle F(x^\ast),x^\ast\rangle\geq 0
	\end{equation}
	Now, consider $\tilde{\lambda} = \tfrac{(\|\lambda^*\|+1)g(\bar{x}_T)_+}{\|g(\bar{x}_T)_+\|}$. Therefore,
	\[
	\langle \tilde{\lambda},g(\bar{x}_T) \rangle =(\|\lambda^\ast\|+1) \|g(\bar{x}_T)_+\| .
	\] 
	From Jensen's inequity  we have 
	\begin{equation}\label{eq:feasb2}
		\begin{aligned}
			\tfrac{1}{T}\tsum_{t=0}^{T-1}\big(\Lcal( x^{t+1},\tilde{\lambda},\hat{v})-\Lcal(x^\ast, \lambda^{t+1}, v^{t+1})\big)&\geq \Lcal( \bar{x}_T,\tilde{\lambda},\hat{v})-\Lcal(x^\ast,\lambda^\ast,v^\ast)\\
			& = \langle F(x^\ast),\bar{x}_T-x^\ast\rangle + (\|\lambda^\ast\|+1) \|g(\bar{x}_T)_+\| .
		\end{aligned}
	\end{equation}
	Also note that for $\tilde{z} = (x^*,\til{\lambda},\hat{v})$ 
	\begin{equation}\label{eq:Q-and-L}
		Q(z^{t+1},\til{z}) = \Lcal( x^{t+1},\tilde{\lambda},\hat{v})-\Lcal(x^\ast, \lambda^{t+1}, v^{t+1}).
	\end{equation}
	Therefore, from \eqref{eq:feasb2}-\eqref{eq:Q-and-L}, 
    we have 
	\begin{equation}\label{eq:feasib3}
		\|g(\bar{x}_T)_+\|  \leq 	\tfrac{1}{T}\tsum_{t=0}^{T-1}\big(\Lcal( x^{t+1},\tilde{\lambda},\hat{v})-\Lcal(x^\ast, \lambda^{t+1}, v^{t+1})\big) - (\langle F(x^\ast),\bar{x}_T-x^\ast\rangle + \|\lambda^\ast\| \|g(\bar{x}_T)_+\|)
	\end{equation}
	Note that from \eqref{eq:feasib1a}, \eqref{eq:Q-and-L} and \eqref{eq:feasib3}, we have 
	\begin{equation}\label{eq:feasib4}
		\|g(\bar{x}_T)_+\|  \leq 	\tfrac{1}{T}\tsum_{t=0}^{T-1}\Lcal( x^{t+1},\tilde{\lambda},\hat{v})-\Lcal(x^\ast, \lambda^{t+1}, v^{t+1}) = \tfrac{1}{T} \tsum_{t=0}^{T-1}Q(z^{t+1}, \wt{z}).
	\end{equation}
	Note that  
\begin{equation}\label{eq:tilde1}
	\|\tilde{\lambda} -\lambda^{t+1}\|^2 \leq 2\|\tilde{\lambda}\|^2 + 2\|\lambda^{t+1}\|^2
\end{equation}
	Moreover, from \eqref{eq:feasib4}, using Lemma \ref{lemma:fully-Stoch-Partial} with $z = \wt{z}$, taking expectation on both sides of the resulting relation, noting that  $\lambda^0 =\zero$, using Lemma \ref{lem:dstoch_useful_bounds-mo} and stepsize policy \eqref{eq:conditions}, 
    and the following relation 
 \begin{equation*}
     \Ebb[\tsum_{i=0}^{T-1} \inprod{\delta_\gfrak^i}{\tilde{\lambda}-\lambda^i}]  = \Ebb[\tsum_{i=0}^{T-1} \inprod{\delta_\gfrak^i}{\tilde{\lambda}}] \leq \sqrt{T}M_{\tilde{\lambda}} \sigma_\gfrak = \sqrt{T}(\|\lambda^*\|+1) \sigma_\gfrak,
 \end{equation*}
 we get \eqref{eq:feasib-full}. Finally, the equality in the above relation follows 
 from the fact that $\Ebb[\inprod{\delta_\gfrak^i}{\lambda^i}] = 0$ and the inequality holds due to Proposition \ref{prop:lemma2lan}. Hence, we conclude the proof. \hfill $\blacksquare$

\section{Operator Constraint Extrapolation method for FCVI problem} \label{sec:dstoch_opconex}
In the previous section, we presented the \sopconex{} method for convergence for nonsmooth, fully stochastic FCVI problems. In this section, we significantly improve the convergence rate guarantee of this method by making a critical change to the constraint extrapolation term. Hence, we call this method the Operator Constraint Extrapolation method. 

We denote the linear approximation of $g(\cdot)$ at point $x^t$ as 
\[ \ell_g(x^t) = g(x^{t-1}) +  g'(x^{t-1})^T(x^t-x^{t-1}),\]
where $ g'(x^{t-1}) = [ g'_1(x^{t-1}), \dots, g'_m(x^{t-1})]$ as defined earlier. Note that for any point $x^t$, we always linearize at the previous point $x^{t-1}$. We use a stochastic version of $\ell_g$ in our algorithm. In particular, using
the SO in \eqref{eq:SO_g_operator-mo}, we have
\begin{equation}
    \ell_\gfrak^t(x^t):= \gfrak(x^{t-1},\xibar^t) + \Gfrak(x^{t-1}, \xibar^t)^T(x^t-x^{t-1}), 
\end{equation}
where $\Gfrak(x^{t-1}, \xibar^t) := [\Gfrak_1(x^{t-1}, \xibar^t), \dots, \Gfrak_m(x^{t-1}, \xibar^t)] \in \Rbb^{n\times m}$. Here, we used $\xibar^t$ as a copy of the random variable $\xi$ independent of $\xi^t$. 
Hence, when $x^{t-1}$ is fixed, $\gfrak_j(x^{t-1}, \xibar^t),\Gfrak_j(x^{t-1}, \xi^t)$ and $\gfrak_j(x^{t-1}, \xi^t), \Gfrak_j(x^{t-1}, \xibar^t)$ are conditionally independent and unbiased estimators of $g_j(x^{t-1}), g'_j(x^{t-1})$, respectively. As we show later, the independent sample $\xibar^t$ ensures that $\ell_\gfrak^t(x^t)$ is an unbiased estimator of $\ell_g(x^t)$. To signify the use of random sample $\xibar^t$, we add a superscript $t$ in the notation of $\ell_\gfrak^t$.

We are ready to formally describe the Operator Constraint Extrapolation (\opconex)~method. This algorithm, while borrowing some elements of \sopconex{}, introduces a new constraint linearization term. The inspiration for this design comes from \cite{boob2023stochastic}, where they show that the linearization of constraints before taking extrapolation leads to a stable and accelerated update of $\lambda$. Accordingly, this single-loop method 
extrapolates constraints after their linearization (line 4 of Algorithm \ref{alg:alg4}). This extrapolation term is missing in \sopconex~method (compare with line 4 of Algorithm \ref{alg:algpdfullstoch}). Moreover, it can handle the stochastic operator as well as constraints using SO in \eqref{eq:SO_F_operator-mo}, \eqref{eq:SO_g_operator-mo} and does not need \eqref{eq:SO_g_convexity}. Finally, note that \opconex{} method does not require convexity and continuity assumption on function $\gfrak(\cdot, \xi)$. Hence, our convergence guarantees in this section hold under milder assumptions than those in Section \ref{sec:full-stoch}.
The final algorithm is presented below.
\begin{algorithm}[H]
	\begin{algorithmic}[1]
		\State {\bf Input: }$x^{0} \in X$, $\lambda^0 = \zero$, $\gfrak(x^0, \wb{\xi}^0)$ and $\Ffrak^0$. 
		\State Set $x^{-1} = x^{0}$, $\Ffrak^{-1} = \Ffrak^0$ and $\ell_\gfrak^0(x^0) = \ell_\gfrak^0(x^{-1}) = \gfrak(x^0, \wb{\xi}^0)$.
		\For{$t=0,1,2,\dots,T-1$}
		\State $\sfrak^{t} \gets (1+\theta_t)\ell_\gfrak^t(x^t) - \theta_t \ell_\gfrak^t(x^{t-1})$
		\State $\lambda^{t+1} \gets \argmin_{\lambda\ge \zero} \inprod{-\sfrak^t}{\lambda} + \tfrac{\tau_t}{2}\gnorm{\lambda-\lambda^t}{}{2}$
		\State $x^{t+1} \gets \argmin_{x \in X} \inprod{(1+\theta_t)\Ffrak^t - \theta_t\Ffrak^{t-1} +\tsum_{i =1}^m \lambda^{t+1}_i\Gfrak_i(x^t, \xi^t)  }{x} + \tfrac{\eta_t}{2}\gnorm{x-x^t}{}{2}$
		\EndFor
		\State {\bf Output: }$\wb{x}^T := \paran{\tsum_{t=0}^{T-1} \gamma_{t}x^{t+1}}\big/\paran{\tsum_{t=0}^{T-1}\gamma_t}$
	\end{algorithmic}
	\caption{Operator Constraint Extrapolation (\opconex) method}\label{alg:alg4}
\end{algorithm}
In the following subsection, we discuss the convergence guarantees of Algorithm \ref{alg:alg4} and compare it with the guarantees of \sopconex{} method. Then, we will present the convergence analysis.

\subsection{Convergence guarantees of \opconex{} method}
We have the following convergence rate guarantee for Algorithm \ref{alg:alg4}.\vspace{-0.8em}
\begin{thm}\label{thm:full_stoch_opconex_conv}
    Let $B \ge 1$ and  $\sigma_{X, g}:= \sqrt{\sigma_\gfrak^2 + D_X^2\gnorm{\sigma_\Gfrak}{}{2}}$. Suppose we set
    \begin{equation}\label{eq:step_policy_dstoch}
			\gamma_t = \theta_t  = 1,
            \qquad
			\eta_t =  L_gB + \eta,
            \qquad
			\tau_t = \tau,
	\end{equation}
	where $\eta = 8L+\tfrac{8M_gB}{D_X} + \tfrac{2(H + H_gB + \sqrt{2}\sigma + 4B\gnorm{\sigma_\Gfrak}{}{})}{D_X}\sqrt{T}$ and $\tau= \tfrac{9D_X}{B}\max\{M_g, \gnorm{\sigma_\Gfrak}{}{}\} + \tfrac{8\sigma_{X, g}}{B}\sqrt{T}$. Then, we have
    \begin{align}
        \Ebb[\sup_{x\in \wt{X}}\inprod{F(x)}{\wb{x}^T -x}] 
        &\le   \tfrac{3[(8L+L_gB)D_X + 8M_gB]D_X}{2T} 
        + \tfrac{3(H+\sqrt{2}\sigma + H_gB+ 4\gnorm{\sigma_\Gfrak}{}{}B)D_X}{\sqrt{T}} 
        \nonumber\\
        &\qquad+ \tfrac{(9\omega^2 + 8H^2 + 9\sigma^2)D_X}{4(H+ H_gB+ \sqrt{2}\sigma + B\gnorm{\sigma_\Gfrak}{}{} + 4M_gB)\sqrt{T}} 
        + \tfrac{9\sigma_{X,g}B}{8\sqrt{T}},
        \label{eq:conv_dstoch_opt}\\
        \Ebb[\gnorm{[g(\wb{x}^T)]_+}{}{}]
        &\le \tfrac{[(8L+L_gB)D_X + 8M_gB]D_X}{2T} + \tfrac{9\max\{M_g, \gnorm{\sigma_\Gfrak}{}{}\}D_X(\gnorm{\lambda^*}{}{}+1)^2}{BT}
        + \tfrac{(H+\sqrt{2}\sigma + H_gB+ 4\gnorm{\sigma_\Gfrak}{}{}B)D_X}{\sqrt{T}} 
        \nonumber\\
        &\qquad\tfrac{2(H^2 + \mathcal{H}_*^2 + 2\sigma^2+ \omega^2)D_X}{(H+H_gB+ \sqrt{2}\sigma + B\gnorm{\sigma_\Gfrak}{}{}) \sqrt{T}} + \tfrac{2\sigma_{X,g}}{\sqrt{T}}\big(B+ \tfrac{4(\gnorm{\lambda^*}{}{}+1)^2}{B}\big),
        \label{eq:conv_dstoch_feas}
    \end{align}
     where $\omega = 2.2\bracket[\big]{ 
    (6\gnorm{\lambda^*}{}{2} + 18B^2 + \tfrac{(8L+L_gB)BD_X}{M_g})\gnorm{\sigma_\Gfrak}{}{2} + 4(H + \sqrt{2}\sigma + \tfrac{\mathcal{H}_*^2}{2(H+ H_gB + \sqrt{2}\sigma)}  + \tfrac{H_gB}{2})B\gnorm{\sigma_\Gfrak}{}{}}^{1/2}$ and $\mathcal{H}_* :=  H_g(\gnorm{\lambda^*}{}{}+1) + \tfrac{L_gD_X}{2}[\gnorm{\lambda^*}{}{} + 1 -B]_+$.
\end{thm}
Theorem \ref{thm:full_stoch_opconex_conv} provides a unified convergence result for function-constrained variational inequalities. We divide the convergence rate analysis into two cases where the problem is either nonsmooth deterministic or general stochastic.  
First, when $\sigma = \sigma_\gfrak = \sigma_\Gfrak = 0$ (deterministic case), the error convergence rate of this method after $T$ iterations is of 
\begin{equation}\label{eq:opconex_main_conv}
    O\big(\tfrac{(L+L_gB)D_X^2}{T} + \tfrac{M_gD_X}{T}(B + \tfrac{(\gnorm{\lambda^*}{}{}+1)^2}{B}) + \tfrac{1}{\sqrt{T}}(H+H_gB + \tfrac{\mathcal{H}_*^2}{H+H_gB})\big).
\end{equation}  
When $B \ge \gnorm{\lambda^*}{}{}+1$, we have $\mathcal{H}_* \le H_gB$ and the convergence error in \eqref{eq:opconex_main_conv} is of
\[O\big(\tfrac{(L+L_gB)D_X^2}{T} + \tfrac{M_gD_XB}{T} + \tfrac{H+H_gB}{\sqrt{T}}\big),\]
and we can see that $B$ serves as the de facto radius of the dual set $\{\lambda \ge \zero\}$ though it is not explicitly stated. If $B < \gnorm{\lambda^*}{}{}+1$, we have $\mathcal{H}_* = H_g(\gnorm{\lambda^*}{}{}+1) + \tfrac{L_gD_X}{2} [\gnorm{\lambda^*}{}{}+1 - B]_+$. This may lead to a nontrivial impact of the Lipschitz constant $L_g$ on the convergence rate. Even then, the rate of convergence is of the same order and the algorithm does not need explicit knowledge of $\gnorm{\lambda^*}{}{}$ to ensure convergence. 

Finally, note that convergence rate in terms of $M_g$ is of $O(\tfrac{M_g}{T})$ irrespective of whether $B \ge \gnorm{\lambda^*}{}{}+1$. This is in sharp contrast with the \sopconex~method whose convergence is of $O(\tfrac{M_g}{\sqrt{T}})$ as observed in \eqref{eq:deter} (a consequence of the main convergence Theorem \ref{thm:feasib-fully-stoch}). This order of magnitude difference can be quite costly when $M_g$ is large. Consider the situation where $g(x) = x^TAx + b\gnorm{x}{}{} - c$ is a sum of a smooth quadratic and a nonsmooth norm function. Then, we can set $M_g = \gnorm{A}{}{}\max_{x \in X} \gnorm{x}{}{} + b$ whereas $H_g = b$. Now, if $X$ has a large radius and $b$ is a reasonable number then clearly $O(\tfrac{M_g}{T} + \tfrac{H_g}{\sqrt{T}}) \ll O(\tfrac{M_g}{\sqrt{T}})$ and hence, the \opconex~method would converge much faster than \sopconex~method. The large radius of $X$ is quite common. E.g., if $X$ is a $\ell_\infty$-ball of radius 1, then its $\max_{x\in X}\gnorm{x}{}{} = \sqrt{n}$. The strength of \opconex~method vis-à-vis \sopconex~method comes from the fact that its primal update parameter $\eta$ depends on $M_g$ as $O(\tfrac{M_g}{D_X})$. On the contrary, $\eta$ in \sopconex~method is of $\Theta(\tfrac{M_g\sqrt{T}}{D_X})$. Given that the number of iterations $T$ is large for nonsmooth problems, the primal step-length $\tfrac{1}{\eta}$ is an order of magnitude smaller for \sopconex~compared to \opconex~method. This is, however, necessary for \sopconex~method to ensure that $\lambda$ remains bounded. That is not the case for \opconex~method since it has the modified constraint extrapolation update in the dual (line~4 of Algorithm~\ref{alg:alg4}). 

Now, we discuss the convergence of \opconex~method for the stochastic case. When the problem is semi-stochastic, i.e., $\sigma > 0$ and $\sigma_\gfrak = \sigma_\Gfrak = 0$, the expected error convergence rate after $T$ iterations of \opconex~method is of
\begin{equation*}
    O\big(\tfrac{(L+L_gB)D_X^2}{T} + \tfrac{M_gD_X}{T}(B + \tfrac{(\gnorm{\lambda^*}{}{}+1)^2}{B}) + \tfrac{1}{\sqrt{T}}(H+H_gB + \sigma + \tfrac{\mathcal{H}_*^2}{H+H_gB})\big).
\end{equation*}
Compared with \eqref{eq:opconex_main_conv}, this term has an additional $O(\tfrac{\sigma}{\sqrt{T}})$ in the convergence while retaining the benefits of \opconex~mentioned earlier.  For the fully-stochastic case, we have convergence rate of
\begin{equation}\label{eq:fstopconex-conv}
    O\big(\tfrac{(L+L_gB)D_X^2}{T} + \tfrac{M_gD_X}{T}(B + \tfrac{(\gnorm{\lambda^*}{}{}+1)^2}{B}) + \tfrac{1}{\sqrt{T}}(H+H_gB + \sigma + \tfrac{\mathcal{H}_*^2}{H+H_gB}) + \tfrac{1}{\sqrt{T}}(\gnorm{\sigma_\Gfrak}{}{}B + \tfrac{\omega^2}{\gnorm{\sigma_\Gfrak}{}{}B}) +  \tfrac{\sigma_{X,g}}{\sqrt{T}}(B + \tfrac{(\gnorm{\lambda^*}{}{}+1)^2}{B})\big).
\end{equation}
Note that if $B \ge \gnorm{\lambda^*}{}{}+1$ then, we have $\mathcal{H}_* = O(H_gB)$. Furthermore, under mild assumptions of $L_gD_X = O(M_g)$ and $H+H_gB+\sigma = O(B\gnorm{\sigma_\Gfrak}{}{})$, we can see that $\omega = O(B\gnorm{\sigma_\Gfrak}{}{})$. This implies the convergence rate in \eqref{eq:fstopconex-conv} for \opconex~method is of 
\[ O\paran[\Big]{\tfrac{((L+L_gB)D_X^2 + M_gD_XB}{T} + \tfrac{(H+H_gB + \sigma)D_X + (\sigma_\gfrak + \gnorm{\sigma_\Gfrak}{}{}D_X)B}{\sqrt{T}}}.\]
The overall convergence rate above is $O(\tfrac{1}{\sqrt{T}})$. However, the dependence of convergence rate on $M_g$ is still of $O(\tfrac{1}{T})$. Hence, \opconex~method can effectively reduce the impact of Lipschitz constant $M_g$. 
If $B < \gnorm{\lambda^*}{}{}+1$, then the multipliers in the convergence rate bound can be worse. Even then, the convergence is of the same order in terms of $T$, i.e., $O(\tfrac{1}{\sqrt{T}})$.

\subsection{Convergence analysis of \opconex{} method}
We now present the convergence analysis of Algorithm \ref{alg:alg4} that leads up to Theorem \ref{thm:full_stoch_opconex_conv}. The following lemma states an important intermediate result 
for Algorithm \ref{alg:alg4}. \vspace{-0.5em}
\begin{lem}\label{prop:dstoch_1st_result}
    Suppose that the step-size parameter $\{\gamma_t, \theta_t, \eta_t, \tau_t\}$ satisfy
    \begin{equation}\label{eq:step_condn_1}
			\gamma_t\eta_t \le \gamma_{t-1}\eta_{t-1},\qquad
			\gamma_{t}\tau_t \le \gamma_{t-1}\tau_{t-1},\qquad
			\gamma_t\theta_t =\gamma_{t-1},
	\end{equation}as well as
    \begin{equation}\label{eq:step_condn_dstoch-1}
		\begin{aligned}
			64(\gamma_t\theta_t  L)^2 &\le \gamma_t\gamma_{t-1}(\eta_t - L_gB)(\eta_{t-1} - L_gB),\\
			72(\gamma_t\theta_t M_g)^2 &\le \gamma_t\gamma_{t-1} \tau_t(\eta_{t-1} - L_gB),\\
		\end{aligned}		
	\end{equation}
	for all $t \ge 1$, and satisfy
	\begin{equation}\label{eq:step_condn_dstoch-2}
		72(\gamma_t\theta_t M_g)^2\le \gamma_t\gamma_{t-2} \tau_t(\eta_{t-2} - L_gB),\hspace{7.5em}
	\end{equation}
	for all $t \ge 2$. Moreover, the parameters also satisfy
    \begin{equation}\label{eq:step_condion_dstoc-3}
        \begin{aligned}
            \tfrac{6M_g^2\gamma_{t}}{\tau_{t}} - \tfrac{\gamma_{t-1}(\eta_{t-1} - L_gB)}{12} &\le 0, \\
            \tfrac{2\gamma_{t}L^2}{\eta_{t}} + \tfrac{3\gamma_tM_g^2}{\tau_{t}} - \tfrac{3\gamma_{t}(\eta_{t} - L_gB)}{16} &\le 0.
        \end{aligned}
    \end{equation}
    Then, we have for all $x \in X$,
    \begin{align}
        \Gamma_{t+1}[&\inprod{F(x)}{\wb{x}^{t+1}-x} + \inprod{\lambda}{g(\wb{x}^{t+1})} -\inprod{\wb{\lambda}^{t+1}}{g(x)}]
        \nonumber\\
        &\le\tfrac{\gamma_0\eta_0}{2} \gnorm{x-x^0}{}{2} 
        + \big[\tfrac{\gamma_0\tau_0}{2}\gnorm{\lambda-\lambda^{0}}{}{2}
        - \tfrac{\gamma_{t}\tau_{t}}{4}\gnorm{\lambda-\lambda^{t+1}}{}{2}\big]
        \nonumber\\
        &\quad 
        + \tsum_{i=0}^{t}\Big[
        \gamma_i\inprod{\delta_\sfrak^i}{\lambda^i- \lambda} 
        -\gamma_i\inprod{\delta_\Gfrak^i}{x^{i}-x}
        -\gamma_i\inprod{\delta^{i+1}}{x^{i+1}-x}
        \Big] + \tfrac{2\gamma_tH^2}{\eta_t} + \tfrac{\gamma_t}{\eta_t}\gnorm{\delta^{t+1} - \delta^{t}}{}{2}
        \nonumber\\
         &\quad + \tsum_{i= 0}^t \bracket[\big]{
        \tfrac{\gamma_i}{\tau_i} \gnorm{\delta_\sfrak^i}{}{2} 
        + \tfrac{4\gamma_i}{\eta_i-L_gB}\gnorm{\delta_\Gfrak^i}{}{2}
        + \tfrac{4\gamma_i \mathcal{H}(\lambda)^2}{\eta_i - L_gB}} + \tsum_{i=1}^t\bracket[\big]{\tfrac{4\gamma_i\theta_i^2 H^2}{\eta_i - L_gB}
        + \tfrac{4\gamma_i\theta_i^2}{\eta_i - L_gB}\gnorm{\delta^i-\delta^{i-1}}{}{2} }.\label{eq:main_rel_dstoch}
    \end{align}
    where
    $\delta_\Gfrak^t := \tsum_{j=1}^m \lambda^{t+1}_j [\Gfrak_j(x^t, \xi^t) -  g'_j(x^t)]$, 
    $\delta_\sfrak^t:= \sfrak^t -[(1+\theta_t)\ell_g(x^t) -\theta_t\ell_g(x^{t-1})]$,
    $\mathcal{H}(\lambda) := \tfrac{L_gD_X}{2}[\gnorm{\lambda}{}{}-B]_+ + H_g\gnorm{\lambda}{}{}$, $\Gamma_{t+1} := \tsum_{i=0}^t\gamma_i$, $\wb{x}^{t+1} := \tfrac{1}{\Gamma_{t+1}} \tsum_{i=0}^t\gamma_ix^{i+1}$
    and recall that $\delta^{t+1} = \Ffrak^{t+1} -F(x^{t+1})$.
\end{lem}
\begin{Proof}
    Using the optimality of $x^{t+1}$ and Lemma \ref{lem:3-point-mo}, we have for all $x \in X$,
    	\begin{align*}
		&\inprod{(1+\theta_t)\Ffrak^t -\theta_t\Ffrak^{t-1}}{x^{t+1}-x} 
        + \inprod{\lambda^{t+1}}{ \Gfrak(x^t,\xi^t)(x^{t+1}-x)} 
        \le \tfrac{\eta_t}{2}\bracket[\big]{\gnorm{x-x^t}{}{2} - \gnorm{x-x^{t+1}}{}{2} - \gnorm{x^t-x^{t+1}}{}{2}} 
        \end{align*}
        which, after adding and subtracting $\inprod{\Ffrak^{t+1}}{x^{t+1}-x}$ in the first inner product, denoting $\Delta \Ffrak_t := \Ffrak^t  - \Ffrak^{t-1}$, rearranging terms and noting definitions of $\delta_\Gfrak^t, \delta^{t+1}$, can be written as 
        \begin{align}
        \inprod{F(x^{t+1})}{x^{t+1}-x} &+ \inprod{\lambda^{t+1}}{ g'(x^t)^\top(x^{t+1}-x)} 
        \nonumber\\
		& \le \tfrac{\eta_t}{2}\bracket[\big]{\gnorm{x-x^t}{}{2} - \gnorm{x-x^{t+1}}{}{2} - \gnorm{x^t-x^{t+1}}{}{2}} 
        + \inprod{\Delta \Ffrak_{t+1}}{x^{t+1}-x} 
        - \theta_t\inprod{ \Delta \Ffrak_t}{x^t-x} 
        \nonumber\\
        &\quad - \theta_t\inprod{\Delta \Ffrak_t}{x^{t+1}-x^t} 
        - \inprod{\delta^{t+1} + \delta_\Gfrak^t}{x^{t+1}-x}. \label{eq:int_rel42}
	\end{align}
    Moreover, from convexity of $g$, monotonicity of $F$ and definition of $\delta^t$, we have,
	\begin{align}
		\inprod{g'(x^t)}{x^{t+1}-x} &= \inprod{ g'(x^t)}{x^{t+1}-x^t} + \inprod{ g'(x^t)}{x^t-x} \nonumber\\
		&\ge\inprod{ g'(x^t)}{x^{t+1}-x^t} + g(x^t) -g(x) = \ell_g(x^{t+1}) -g(x) \label{eq:int_rel3}\\
		\inprod{F(x^{t+1})}{x^{t+1}-x} &\ge \inprod{F(x)}{x^{t+1}-x}, \qquad \forall \ x\in X\label{eq:int_rel4}\\
        \inprod{\Delta \Ffrak_t}{x^{t+1}-x^t} & =  \inprod{\Delta F_t}{x^{t+1}-x^t} + \inprod{\delta^{t}- \delta^{t-1}}{x^{t+1}-x^t} \label{eq:int_rel5-1}
	\end{align}
    where in the last relation recall that $\Delta F_t = F(x^t) - F(x^{t-1})$.
    Using \eqref{eq:int_rel3} - \eqref{eq:int_rel5-1} in \eqref{eq:int_rel42} and the fact that $\lambda^{t+1} \ge \zero$, 
    we have for all $x \in X$,
    \begin{align}
        \inprod{F(x)}{x^{t+1}-x} &+ \inprod{\lambda^{t+1}}{\ell_g(x^{t+1})-g(x)} 
        \nonumber\\
		& \le \tfrac{\eta_t}{2}\bracket[\big]{\gnorm{x-x^t}{}{2} - \gnorm{x-x^{t+1}}{}{2} - \gnorm{x^t-x^{t+1}}{}{2}} 
        + \inprod{\Delta \Ffrak_{t+1}}{x^{t+1}-x} 
        - \theta_t\inprod{ \Delta \Ffrak_t}{x^t-x} 
        \nonumber\\
        &\quad - \theta_t[\inprod{\Delta F_t}{x^{t+1}-x^t} 
        + \inprod{\delta^{t+1}- \delta^t}{x^{t+1}-x^t}]
        - \inprod{\delta^{t+1} + \delta_\Gfrak^t}{x^{t+1}-x}.
        \label{eq:int_rel43}
    \end{align}
    Note that $\lambda^{t+1} = \argmin\limits_{\lambda \ge \zero} \inprod{-\sfrak^t}{\lambda} + \tfrac{\tau_t}{2} \gnorm{\lambda-\lambda^t}{2}{2}$. Hence, using Lemma \ref{lem:3-point-mo}, we have for all $\lambda \ge \zero$,
	\begin{equation}\label{eq:int_rel1-1} 
		-\inprod{\sfrak^t}{\lambda^{t+1} -\lambda} \le \tfrac{\tau_t}{2} \bracket*{\gnorm{\lambda-\lambda^t}{}{2} - \gnorm{\lambda^{t+1}-\lambda^t}{}{2} - \gnorm{\lambda-\lambda^{t+1}}{}{2}}.
	\end{equation}
    Adding relations \eqref{eq:int_rel43} and \eqref{eq:int_rel1-1}, 
    defining $s^t := (1+\theta_t)\ell_g(x^t) - \theta_t\ell_g(x^{t-1})$ and noting that $\delta_\sfrak^t = \sfrak^t - s^t$ and, we have for all $x\in X$ and $\lambda \ge \zero$
    \begin{align}
         \inprod{F(x)}{x^{t+1}-x} &+ \inprod{\lambda^{t+1}}{\ell_g(x^{t+1})-g(x)}  - \inprod{s^t}{\lambda^{t+1}-\lambda}
        \nonumber\\
        & \le \tfrac{\eta_t}{2}\bracket[\big]{\gnorm{x-x^t}{}{2} - \gnorm{x-x^{t+1}}{}{2} - \gnorm{x^t-x^{t+1}}{}{2}} 
        + \inprod{\Delta \Ffrak_{t+1}}{x^{t+1}-x} 
        - \theta_t\inprod{ \Delta \Ffrak_t}{x^t-x} 
        \nonumber\\
        &\quad + \tfrac{\tau_t}{2}\big[\gnorm{\lambda-\lambda^t}{}{2} 
        - \gnorm{\lambda-\lambda^{t+1}}{}{2} 
        - \gnorm{\lambda^t-\lambda^{t+1}}{}{2} \big] 
        + \inprod{\delta_\sfrak^t}{\lambda^{t+1}- \lambda}
        \nonumber\\
        &\quad - \theta_t[\inprod{\Delta F_t}{x^{t+1}-x^t} 
        + \inprod{\delta^{t}- \delta^{t-1}}{x^{t+1}-x^t}]
        - \inprod{\delta^{t+1} + \delta_\Gfrak^t}{x^{t+1}-x}.
        \label{eq:int_rel44}
    \end{align}
    It is easy to see that
    \begin{align}
    \inprod{\lambda^{t+1}}{\ell_g(x^{t+1})} -\inprod{s^t}{\lambda^{t+1} -\lambda} &= \inprod{\lambda^{t+1} -\lambda}{\ell_g(x^{t+1}) -s^t}+ \inprod{\lambda}{\ell_g(x^{t+1})} \nonumber\\
    &= \inprod{\lambda^{t+1} -\lambda}{q_{t+1}} -\theta_t\inprod{\lambda^{t} -\lambda}{q_t} -\theta_t\inprod{\lambda^{t+1} -\lambda^t}{q_t} + \inprod{\lambda}{\ell_g(x^{t+1})}\label{eq:int_rel-1-45}
    \end{align}
    where $q_t:= \ell_g(x^t) - \ell_g (x^{t-1})$.
    Using Cauchy-Schwarz inequality and \eqref{eq:g-Lipschitz-prop}, we have
     \begin{align}\label{eq:lambda_relation-2}
		\inprod{\lambda}{g(x^{t+1}) -\ell_g(x^{t+1})} &\le \gnorm{\lambda}{}{} \bracket[\big]{\tfrac{L_g}{2} \gnorm{x^{t+1}-x^t}{}{2} + H_g\gnorm{x^{t+1} -x^t}{}{}} \nonumber\\
        &\le \tfrac{BL_g}{2}\gnorm{x^{t+1}-x^t}{}{2} + [\gnorm{\lambda}{}{}-B]_+\tfrac{L_g}{2}\gnorm{x^{t+1}-x^t}{}{2} + H_g\gnorm{\lambda}{}{}\gnorm{x^{t+1} -x^t}{}{}\nonumber\\
        &\le \tfrac{BL_g}{2}\gnorm{x^{t+1}-x^t}{}{2} + \bracket[\big]{\tfrac{L_gD_X}{2}[\gnorm{\lambda}{}{}-B]_+ + H_g\gnorm{\lambda}{}{}} \gnorm{x^{t+1} -x^t}{}{}
	\end{align}
    Using the equality \eqref{eq:int_rel-1-45} into \eqref{eq:int_rel44}, adding the resulting relation with \eqref{eq:lambda_relation-2}, noticing the definition of $\mathcal{H}(\lambda)$, replacing index $t$ by $i$ and rearranging terms, we have
    \begin{align}
        &\inprod{F(x)}{x^{i+1}-x} + \inprod{\lambda}{g(x^{i+1})} -\inprod{\lambda^{i+1}} {g(x)} 
        \nonumber\\
        & \le \tfrac{\eta_i}{2}\bracket[\big]{\gnorm{x-x^i}{}{2} -\gnorm{x-x^{i+1}}{}{2}} 
        +\tfrac{\tau_i}{2}\bracket[\big]{\gnorm{\lambda-\lambda^{i}}{}{2} - \gnorm{\lambda-\lambda^{i+1}}{}{2}} - \inprod{\delta^{i+1} + \delta_\Gfrak^i}{x^{i+1}-x}
        \nonumber\\
        &\quad + \inprod{\Delta \Ffrak_{i+1}}{x^{i+1}-x} 
        - \theta_i\inprod{ \Delta \Ffrak_i}{x^i-x} 
        +\inprod{q_{i+1}}{\lambda - \lambda^{i+1}} 
        -\theta_i \inprod{q_i}{\lambda - \lambda^{i}} 
        + \inprod{\delta_\sfrak^i}{\lambda^i- \lambda} 
        \nonumber\\
        &\quad -\theta_i \inprod{q_i}{\lambda^i - \lambda^{i+1}} 
        + \inprod{\delta_\sfrak^i}{\lambda^{i+1}- \lambda^i} 
        + \theta_i[\inprod{\Delta F_i}{x^{i+1}-x^i} 
        + \inprod{\delta^{i}- \delta^{i-1}}{x^{i+1}-x^i}]
        \nonumber\\
        &\quad - \tfrac{\tau_i}{2}\gnorm{\lambda^i-\lambda^{i+1}}{}{2} 
        - \tfrac{\eta_i - L_gB}{2}\gnorm{x^i-x^{i+1}}{}{2} 
        + \mathcal{H}(\lambda)\gnorm{x^{i+1}-x^i}{}{}. \label{eq:before_sum}
    \end{align}
    Multiplying the above relation by $\gamma_t$, then summing from $i= 0$ to $t$, using step-size relation \eqref{eq:step_condn_1}, and noting that $\Delta \Ffrak_0 = \Delta F_0 = q_0 = \delta^0 - \delta^{-1} = \zero$, we have for all $x\in X$ and $\lambda \ge \zero$
    \begin{align}
        &\Gamma_{t+1}[\inprod{F(x)}{\wb{x}^{t+1}-x} + \inprod{\lambda}{g(\wb{x}^{t+1})} -\inprod{\wb{\lambda}^{t+1}}{g(x)}]
        \nonumber\\
        &\le\big[ \tfrac{\gamma_0\eta_0}{2} \gnorm{x-x^0}{}{2} 
        - \tfrac{\gamma_{t}\eta_{t}}{2}\gnorm{x-x^{t+1}}{}{2}\big]
        + \big[\tfrac{\gamma_0\tau_0}{2}\gnorm{\lambda-\lambda^{0}}{}{2}
        - \tfrac{\gamma_{t}\tau_{t}}{2}\gnorm{\lambda-\lambda^{t+1}}{}{2}\big]
        \nonumber\\
        &\quad + \gamma_{t}\inprod{\Delta \Ffrak_{t+1}}{x^{t+1}-x} 
        +\gamma_{t}\inprod{q_{t+1}}{\lambda - \lambda^{t+1}} 
        + \tsum_{i=0}^{t}\Big[\gamma_i\inprod{\delta_\sfrak^i}{\lambda^i- \lambda} 
        -\gamma_i\inprod{\delta_\Gfrak^i}{x^{i}-x}
        -\gamma_i\inprod{\delta^{i+1}}{x^{i+1}-x}
        \Big]
        \nonumber\\
        &\quad + \tsum_{i = 1}^{t}\Big[ 
        \gamma_i\theta_i \inprod{q_i}{\lambda^{i+1}-\lambda^i} 
        + \gamma_i\theta_i[\inprod{\Delta F_i}{x^{i+1}-x^i} 
        + \inprod{\delta^{i}- \delta^{i-1}}{x^{i+1}-x^i}] \Big]
        \nonumber\\
        &\quad +\tsum_{i=0}^t\Big[ \gamma_i \inprod{\delta_\sfrak^i}{\lambda^{i+1} - \lambda^i}
        - \gamma_i\inprod{\delta_\Gfrak^i}{x^{i+1}-x^i}
        +\gamma_i\mathcal{H}(\lambda)\gnorm{x^{i+1}-x^i}{}{}\Big]\nonumber\\
        &\quad
        - \tsum_{i=0}^t\big[\tfrac{\gamma_i\tau_i}{2}\gnorm{\lambda^i-\lambda^{i+1}}{}{2} 
        + \tfrac{\gamma_i(\eta_i - L_gB)}{2}\gnorm{x^i-x^{i+1}}{}{2} 
        \big]. \label{eq:int_rel46_1}
    \end{align}
    Now, we provide some relation that will be useful in further upper-bounding the RHS of the above relation.
    Noting the definition of $\Delta\Ffrak_t$ and $\delta^t$, and using Young's inequality, we have
     \begin{align}
		\gamma_{t} \inprod{\Delta \Ffrak_{t+1}}{x^{t+1}-x} &- \tfrac{\gamma_{t}\eta_{t}}{2} \gnorm{x-x^{t+1}}{}{2} 
        \nonumber\\
        &= 
        \gamma_{t} \inprod{\Delta F_{t+1}}{x^{t+1}-x} + \gamma_{t} \inprod{\delta^{t+1} - \delta^{t}}{x^{t+1}-x}- \tfrac{\gamma_{t}\eta_{t}}{2} \gnorm{x-x^{t+1}}{}{2} \nonumber\\
        &\le \tfrac{\gamma_{t}}{\eta_{t}}\gnorm{\Delta F_{t+1}}{}{2} + \tfrac{\gamma_{t}}{\eta_{t}}\gnorm{\delta^{t+1} - \delta^{t}}{}{2}\nonumber\\
        &\le \tfrac{2\gamma_{t}L^2}{\eta_{t}}\gnorm{x^{t+1} - x^{t}}{}{2}+ \tfrac{2\gamma_tH^2}{\eta_t} + \tfrac{\gamma_{t}}{\eta_{t}}\gnorm{\delta^{t+1} - \delta^{t}}{}{2}\label{eq:int_rel36}
	\end{align}
    where the last relation follows by \eqref{eq:F-Lipschitz-property} and the fact that $(a+b)^2 \le 2a^2 + 2b^2$ for all $a, b\in \Rbb$.
    Similarly, using Young's inequality and \eqref{eq:F-Lipschitz-property}, we have for any $t$ and $i = 1,\dots, t$
	\begin{subequations}\label{eq:LipRel}
		\begin{align}
        &\gamma_i\theta_i \inprod{\Delta F_i}{x^{i+1}-x^i} \le \gamma_i\theta_i L\gnorm{x^i-x^{i-1}}{}{}\gnorm{x^{i+1}-x^i}{}{} + \gamma_i \theta_i H \gnorm{x^{i+1}-x^i}{}{}, \label{eq:lip_rel0},\\
        &\gamma_i \theta_i L \gnorm{x^i-x^{i-1}}{}{} \gnorm{x^{i+1}-x^i}{}{} \le \tfrac{\gamma_i(\eta_i - L_gB)}{16} \gnorm{x^{i+1} - x^i}{}{2} +\tfrac{\gamma_{i-1}(\eta_{i-1} - L_gB )}{16}\gnorm{x^i-x^{i-1}}{}{2},\label{eq:lip_rel1}\\
		&\gamma_i \theta_i H \gnorm{x^{i+1}-x^i}{}{} - \tfrac{\gamma_i(\eta_i - L_gB)}{16} \gnorm{x^{i+1} - x^i}{}{2} \le  \tfrac{4\gamma_i\theta_i^2 H^2}{\eta_i - L_gB} ,\label{eq:lip_rel2}\\ 
		&\gamma_i\theta_i\inprod{\delta^{i}- \delta^{i-1}}{x^{i+1}-x^i} - \tfrac{\gamma_i(\eta_i - L_gB)}{16} \gnorm{x^{i+1} - x^i}{}{2} \le \tfrac{4\gamma_i\theta_i^2}{\eta_i -L_gB}\gnorm{\delta^{i}- \delta^{i-1}}{}{2}\label{eq:lip_rel3}
        \end{align}
	\end{subequations}
    where \eqref{eq:lip_rel1} follows from the first relation in \eqref{eq:step_condn_dstoch-1}. Moreover, from the definition of $q_t$, we have 
    \begin{align*}
        \gnorm{q_t}{}{} &= \gnorm{\ell_g(x^t) - \ell_g(x^{t-1})}{}{} = \gnorm{g(x^{t-1}) +  g'(x^{t-1})^T(x^t-x^{t-1})  - g(x^{t-2}) -  g'(x^{t-2})^T(x^{t-1} - x^{t-2})}{}{}\\
        &\le 2M_g\gnorm{x^{t-1} - x^{t-2}}{}{} + M_g\gnorm{x^t-x^{t-1}}{}{}
    \end{align*}
    Using the above relation and the fact that $(a+b)^2 \le 3a^2 + \tfrac{3}{2}b^2$ for all $a,b\in \Rbb$, we have
    \begin{align}
        &\gamma_t\inprod{q_{t+1}}{\lambda - \lambda^{t+1}} 
        - \tfrac{\gamma_{t}\tau_{t}}{4}\gnorm{\lambda-\lambda^{t+1}}{}{2} 
        \le \tfrac{\gamma_t}{\tau_t} \gnorm{q_{t+1}}{}{2} 
        \le \tfrac{3M_g^2\gamma_t}{\tau_t}[\gnorm{x^{t+1}-x^t}{}{2} 
        + 2\gnorm{x^t-x^{t-1}}{}{2}],
        \nonumber\\
        &\gamma_i\inprod{\delta_\Gfrak^i}{x^{i+1}-x^i} 
        - \tfrac{\gamma_i(\eta_i- L_gB)}{16}\gnorm{x^{i+1}-x^i}{}{2} 
        \le \tfrac{4\gamma_i}{\eta_i-L_gB}\gnorm{\delta_\Gfrak^i}{}{2}, \quad \forall  i =0,\dots, t,
        \nonumber\\
        &\gamma_i\inprod{\delta_\sfrak^i}{\lambda^{i+1}-\lambda^i} 
        - \tfrac{\gamma_i\tau_i}{4}\gnorm{\lambda^{i+1}-\lambda^i}{}{2} 
        \le \tfrac{\gamma_i}{\tau_i}\gnorm{\delta_\sfrak^i}{}{2}, \quad \forall  i =0,\dots, t,
        \nonumber\\
        &\gamma_i\mathcal{H}(\lambda)  \gnorm{x^{i+1}-x^i}{}{} - \tfrac{\gamma_i(\eta_i - L_gB)}{16} \gnorm{x^{i+1} - x^i}{}{2} 
        \le \tfrac{4\gamma_i \mathcal{H}(\lambda)^2}{\eta_i - L_gB},\quad \forall i =0 , \dots, t
        \nonumber\\
        &\gamma_i\theta_i \inprod{q_i}{\lambda^{i+1} - \lambda^i} \le \gamma_i\theta_i \gnorm{q_i}{}{} \gnorm{\lambda^{i+1}- \lambda^i}{}{}\nonumber\\
        &\hspace{8em}\le\gamma_i\theta_i \big[2M_g \gnorm{x^{i-1}-x^{i-2}}{}{} + M_g\gnorm{x^i-x^{i-1}}{}{}\big]\gnorm{\lambda^{i+1}- \lambda^i}{}{}, \quad \forall  i =1,\dots, t, \nonumber\\
        &M_g\gamma_i\theta_i \gnorm{x^i - x^{i-1}}{}{}
        \gnorm{\lambda^{i+1} - \lambda^i}{}{} 
        - \tfrac{\gamma_i\tau_i}{12}\gnorm{\lambda^{i+1}- \lambda^i}{}{2} 
        - \tfrac{\gamma_{i-1}(\eta_{i-1} - L_gB)}{24}\gnorm{x^i-x^{i-1}}{}{2} 
        \le 0, \quad \forall  i =1,\dots, t,
        \label{eq:int_rel46}
        \\
		&2M_g\gamma_i\theta_i \gnorm{x^{i-1} - x^{i-2}}{}{}
        \gnorm{\lambda^{i+1} - \lambda^i}{}{} 
        - \tfrac{\gamma_i\tau_i}{6}\gnorm{\lambda^{i+1}- \lambda^i}{}{2} 
        - \tfrac{\gamma_{i-2}(\eta_{i-2} - L_gB)}{12}
        \gnorm{x^{i-1} - x^{i-2}}{}{2} 
        \le 0, \quad \forall  i =1,\dots, t, 
        \label{eq:int_rel47}
    \end{align}
    where note that \eqref{eq:int_rel46} and \eqref{eq:int_rel47} follow by the second relation in \eqref{eq:step_condn_dstoch-1} and \eqref{eq:step_condn_dstoch-2}, respectively.
    Summing up all of the above relations along with \eqref{eq:int_rel36} and all relations in \eqref{eq:LipRel} for $i = 1, \dots, t$, we have for all $x\in X$ and $\lambda \ge \zero$,
    \begin{align*}
        &\gamma_{t} \inprod{\Delta \Ffrak_{t+1}}{x^{t+1}-x} + \gamma_t\inprod{q_{t+1}}{\lambda - \lambda^{t+1}} - \tfrac{\gamma_{t}\eta_{t}}{2} \gnorm{x-x^{t+1}}{}{2} - \tfrac{\gamma_t\tau_t}{4}\gnorm{\lambda - \lambda^{t+1}}{}{2} \\
        &+ \tsum_{i=1}^t\big[\gamma_i\theta_i \inprod{\Delta F_i }{x^{i+1}-x^i} + \gamma_i\theta_i\inprod{\delta^i - \delta^{i-1}}{x^{i+1}-x^i} + \gamma_i\theta_i\inprod{q_i}{\lambda^{i+1}-\lambda^i} \big]\\
        &+\tsum_{i=0}^t \big[ \gamma_i\inprod{\delta_\Gfrak^i}{x^{i+1}-x^i} + \gamma_i\inprod{\delta_\sfrak^i}{\lambda^{i+1} -\lambda^i} + \gamma_i\mathcal{H}(\lambda)\gnorm{x^{i+1}-x^i}{}{} \big]\\
        & - \tsum_{i=0}^t\big[ \tfrac{2\gamma_i(\eta_i-L_gB)}{16}\gnorm{x^{i+1}-x^i}{}{2} + \tfrac{\gamma_i\tau_i}{4}\gnorm{\lambda^{i+1}-\lambda^i}{}{2}\big] -\tsum_{i=1}^t \tfrac{\gamma_i\tau_i}{4}\gnorm{\lambda^{i+1}-\lambda^i}{}{2}\\
        &-\tsum_{i=1}^t\big[\tfrac{3\gamma_i(\eta_i - L_gB)}{16}\gnorm{x^{i+1}-x^i}{}{2}+ \tfrac{5\gamma_{i-1}(\eta_{i-1}-L_gB)}{48}\gnorm{x^i-x^{i-1}}{}{2} + \tfrac{\gamma_{i-2}(\eta_{i-2}-L_gB)}{12}\gnorm{x^{i-1}-x^{i-2}}{}{2}\big]\\
        &\le \big[\tfrac{3M_g^2\gamma_t}{\tau_t} + \tfrac{2\gamma_tL^2}{\eta_t}\big]\gnorm{x^{t+1}-x^t}{}{2} + \tfrac{6M_g^2\gamma_t}{\tau_t}\gnorm{x^t-x^{t-1}}{}{2} + \tfrac{2\gamma_tH^2}{\eta_t} + \tfrac{\gamma_{t}}{\eta_{t}}\gnorm{\delta^{t+1} - \delta^{t}}{}{2} \\
        &\quad + \tsum_{i=1}^{t}\big[ \tfrac{4\gamma_i\theta_i^2 H^2}{\eta_i - L_gB} + \tfrac{4\gamma_i\theta_i^2}{\eta_i -L_gB}\gnorm{\delta^{i}- \delta^{i-1}}{}{2} \big] + \tsum_{i=0}^t\big[ \tfrac{4\gamma_i}{\eta_i-L_gB}\gnorm{\delta_\Gfrak^i}{}{2} + \tfrac{\gamma_i}{\tau_i}\gnorm{\delta_\sfrak^i}{}{2} + \tfrac{4\gamma_i \mathcal{H}(\lambda)^2}{\eta_i - L_gB} \big]
    \end{align*}
    Summing the above relation with \eqref{eq:int_rel46_1},
    we have for all $x \in X$ and $\lambda \ge \zero$,
    \begin{align*}
        &\Gamma_{t+1}[\inprod{F(x)}{\wb{x}^{t+1}-x} + \inprod{\lambda}{g(\wb{x}^{t+1})} -\inprod{\wb{\lambda}^{t+1}}{g(x)}]
        \nonumber\\
        &\le\tfrac{\gamma_0\eta_0}{2} \gnorm{x-x^0}{}{2} 
        + \big[\tfrac{\gamma_0\tau_0}{2}\gnorm{\lambda-\lambda^{0}}{}{2}
        - \tfrac{\gamma_{t}\tau_{t}}{4}\gnorm{\lambda-\lambda^{t+1}}{}{2}\big]
        \nonumber\\
        &\quad + [\tfrac{2\gamma_{t}L^2}{\eta_{t}} 
        + \tfrac{3M_g^2\gamma_{t}}{\tau_{t}} - \tfrac{3\gamma_t(\eta_t-L_gB)}{16}]
        \gnorm{x^{t+1}-x^{t}}{}{2} + \big[\tfrac{6M_g^2\gamma_{t}}{\tau_{t}} 
        -\tfrac{\gamma_{t-1}(\eta_{t-1}-L_gB)}{12} \big]
        \gnorm{x^{t}-x^{t-1}}{}{2} \nonumber\\
        &\quad
        + \tsum_{i=0}^{t}\Big[
        \gamma_i\inprod{\delta_\sfrak^i}{\lambda^i- \lambda} 
        -\gamma_i\inprod{\delta_\Gfrak^i}{x^{i}-x}
        -\gamma_i\inprod{\delta^{i+1}}{x^{i+1}-x}
        \Big] + \tfrac{2\gamma_tH^2}{\eta_t} + \tfrac{\gamma_t}{\eta_t}\gnorm{\delta^{t+1} - \delta^{t}}{}{2}
        \nonumber\\
         &\quad + \tsum_{i= 0}^t \bracket[\big]{
        \tfrac{\gamma_i}{\tau_i} \gnorm{\delta_\sfrak^i}{}{2} 
        + \tfrac{4\gamma_i}{\eta_i-L_gB}\gnorm{\delta_\Gfrak^i}{}{2}
        + \tfrac{4\gamma_i \mathcal{H}(\lambda)^2}{\eta_i - L_gB}} + \tsum_{i=1}^t \big[\tfrac{4\gamma_i\theta_i^2 H^2}{\eta_i - L_gB}
        + \tfrac{4\gamma_i\theta_i^2}{\eta_i - L_gB}\gnorm{\delta^i-\delta^{i-1}}{}{2}\big] \nonumber\\
        &\quad -\tfrac{9\gamma_0(\eta_0-L_gB)}{48}\gnorm{x^1-x^0}{}{2} - \tfrac{\gamma_0\tau_0}{4}\gnorm{\lambda_1-\lambda^0}{}{2}.
    \end{align*} 
    Substituting the step-size condition~\eqref{eq:step_condion_dstoc-3} in the above relation, noting that $\eta_0 > L_gB$ and dropping the last line of the above bound, we obtain \eqref{eq:main_rel_dstoch}. Hence, we conclude the proof.
\end{Proof}
In the following proposition, we show bounds on some important quantities.
\begin{lem}\label{lem:dstoch_useful_bounds}
    For stochastic oracles satisfying \eqref{eq:SO_F_operator-mo} and \eqref{eq:SO_g_operator-mo}, we have 
    \begin{subequations}\label{eq:bd_var_dstoch}
        \begin{align}
            \Ebb[\gnorm{\delta^{i} - \delta^{i-1}}{}{2}] 
            &\le 2\sigma^2, \label{eq:bd_var_F}\\
            \Ebb[\gnorm{\delta_\Gfrak^i}{}{2}] 
            &\le \gnorm{\sigma_\Gfrak}{}{2}\Ebb[\gnorm{\lambda^{i+1}}{}{2}],
            \label{eq:bd_var_Gfrak}\\
             \Ebb[\gnorm{\delta_\sfrak^i}{}{2}] 
            &\le 2(1+2\theta_i)^2[\sigma_\gfrak^2 +D_X^2\gnorm{\sigma_\Gfrak}{}{2}], \label{eq:bd_var_sfrak},
        \end{align}
    \end{subequations}
    where recall that $\sigma_\Gfrak = [\sigma_1, \dots, \sigma_m]^T$.
    When $x,\lambda$ is fixed (i.e., non-random), we have
    \begin{align}\label{eq:useful-expectation2}
        \Ebb[\inprod{\delta^{i+1}}{x^{i+1}-x}]=0, \qquad
        \Ebb[\inprod{\delta_\Gfrak^i}{x^{i}-x}] = \zero, \qquad 
        \Ebb[\inprod{\delta_\sfrak^i}{\lambda^i- \lambda}] = 0
    \end{align}
    Moreover, assume \eqref{eq:step_condn_1} is satisfied and consider the sequences $\{x^i_a\}$, $\{\wt{x}^i_a\}$ and $\{\lambda^i_a\}$ as follows
    \begin{subequations}
        \begin{align}
            x^1_a = x^0, \quad &x^{i+1}_a := \argmin_{x \in X} -\inprod{\tfrac{\delta^i}{\eta_{i-1}}}{x} + \tfrac{1}{2}\gnorm{x-x^i_a}{}{2},\quad \forall  i\ge1,  \label{eq:x-aug-seq_def}\\
            \wt{x}^0_a = x^0, 
            \quad &\wt{x}^i_a := \argmin_{x\in X} -\inprod{\tfrac{\delta_\Gfrak^{i-1}}{\eta_{i-1}}}{x} 
            + \tfrac{1}{2}\gnorm{x- \wt{x}^{i-1}_a}{}{2},
            \quad \forall  i\ge1, \label{eq:x-aug-seq_def2}\\
            \lambda^0_a = \lambda^0,
            \quad &\lambda_a^i := \argmin_{\lambda \ge 0} \inprod{\tfrac{\delta_\sfrak^{i-1}}{{\eta_{i-1}}}}{\lambda} + \tfrac{1}{2}\gnorm{\lambda-\lambda_a^{i-1}}{}{2}, \qquad\forall i\ge 1.
            \label{eq:lambda-aug-seq_def2}
        \end{align}
    \end{subequations}\vspace{-1em}
    Then, we have 
    \begin{subequations}
        \begin{align}
        \forall x\in X,
        \quad \tsum_{i=0}^t \gamma_i\inprod{\delta^{i+1}}{x-{x}^{i+1}_a} 
        &\le \tfrac{\gamma_0\eta_0}{2}\gnorm{x- x^0}{}{2} 
        + \tsum_{i=0}^t \tfrac{\gamma_i}{2\eta_i}\gnorm{\delta^{i+1}}{}{2}.
        \label{eq:stoch_innerprod_delta}\\
        \forall x\in X, 
        \quad \hspace{1.4em}\tsum_{i=0}^t \gamma_i\inprod{\delta_\Gfrak^i}{x-\wt{x}^i_a} 
        &\le \tfrac{\gamma_0\eta_0}{2}\gnorm{x- x^0}{}{2} 
        + \tsum_{i=0}^t \tfrac{\gamma_i}{2\eta_i}\gnorm{\delta_\Gfrak^i}{}{2}, \label{eq:stoch_innerprod_delta_G}\\
        \forall\lambda \ge \zero,  
        \quad\hspace{1.6em} \tsum_{i=0}^t \gamma_i\inprod{\delta_\sfrak^i}{\lambda^i_a-\lambda}
        &\le \tfrac{\gamma_0\tau_0}{2}\gnorm{\lambda- \lambda^{0}}{}{2} 
        + \tsum_{i=0}^t \tfrac{\gamma_i}{2\tau_i}\gnorm{\delta_\sfrak^i}{}{2}.
        \label{eq:stoch_innerprod_delta_s}        
    \end{align}
    \end{subequations}
\end{lem}
\begin{Proof}
    Denoting $\xi^{[t]} := (\xi^0, \dots, \xi^t)$ and borrowing a similar notation for $\xibar^{[t]}$, it is easy to see that $x^{t+1}$ depends on $(\xi^{[t]}, \xibar^{[t]})$, and $s^t$ as well as $\lambda^{t+1}$ depends on $(\xi^{[t-1]}, \xibar^{[t]})$. Then, we can show
    \eqref{eq:bd_var_F}, \eqref{eq:bd_var_Gfrak}, first two relations of \eqref{eq:useful-expectation2}, \eqref{eq:stoch_innerprod_delta} and \eqref{eq:stoch_innerprod_delta_G} by reproducing the proof of Lemma \ref{lem:dstoch_useful_bounds-mo} for \eqref{eq:bd_var_F-mo}, \eqref{eq:bd_var_Gfrak-mo}, \eqref{eq:useful_expectations-mo}, \eqref{eq:stoch_innerprod_delta-mo} and \eqref{eq:stoch_innerprod_delta_G-mo}, respectively. We focus on proving \eqref{eq:bd_var_sfrak}, third relation in \eqref{eq:useful-expectation2} and \eqref{eq:stoch_innerprod_delta_s}.
    \begin{align*}
        \Ebb[\ell_\gfrak^t(x^t)]
        &= \Ebb[ \gfrak(x^{t-1}, \xibar^t) + \Gfrak(x^{t-1}, \xibar^t)^T(x^t-x^{t-1}) ],
        \\
        &=\Ebb\big[ \Ebb[\gfrak(x^{t-1}, \xibar^t)|\xi^{[t-1]},\xibar^{[t-1]}] + \Ebb[\Gfrak(x^{t-1}, \xibar^t)|\xi^{[t-1]}, \xibar^{[t-1]}]^T(x^t-x^{t-1})\big] = \ell_g(x^t),\\
        \Ebb[\ell_\gfrak^t(x^{t-1})]
        &= \Ebb[ \gfrak(x^{t-2}, \xibar^t) + \Gfrak(x^{t-2}, \xibar^t)^T(x^{t-1}-x^{t-2}) ],
        \\
        &=\Ebb\big[ \Ebb[\gfrak(x^{t-2}, \xibar^t)|\xi^{[t-1]},\xibar^{[t-1]}] + \Ebb[\Gfrak(x^{t-2}, \xibar^t)|\xi^{[t-1]}, \xibar^{[t-1]}]^T(x^{t-1}-x^{t-2})\big] = \ell_g(x^{t-1}),\\
        \Ebb[\delta_\sfrak^t] &= \Ebb[\sfrak^t -s^t]
        =\Ebb\big[(1+\theta_t)\Ebb[\ell_\gfrak(x^t)] 
        - \theta_t\Ebb[\gfrak(x^{t-1})] -s^t\big]= \zero,
    \end{align*}
    where the last relation is a consequence of the preceding two relations, final equality in the first two relations follows from the properties of stochastic oracles in \eqref{eq:SO_g_operator-mo}, and the rest of the equalities follow from the linearity and tower law of the expectation. Using the above relations and the fact that $\lambda^i$ depends on $(\xi^{[i-2]}, \xibar^{i-1})$, we have
    \begin{align*}
        \Ebb[\inprod{\delta_\sfrak^i}{\lambda^i- \lambda}] &= \Ebb[\inprod{\Ebb[\delta_\sfrak^i|\xi^{[i-1], \xibar^{[i-1]}}]}{\lambda^i- \lambda}] = \zero,
    \end{align*}
    for any non-random $x$ and $\lambda$. Hence, we prove \eqref{eq:useful-expectation2}.

    
   Noting the definition of the sequence $\{\lambda^i_a\}_{i\ge 0}$ in \eqref{eq:lambda-aug-seq_def2} 
    and applying Proposition~\ref{prop:tech_res1-mo}, we have for all $\lambda \ge \zero$
    \[\inprod{\tfrac{\delta_\sfrak^i}{\tau_i}}{\lambda^i_a-\lambda} 
    \le \tfrac{1}{2}\gnorm{\lambda- \lambda^{i}_a}{}{2} 
    - \tfrac{1}{2}\gnorm{\lambda- \lambda^{i+1}_a}{}{2}
    + \tfrac{1}{2\tau_i^2}\gnorm{\delta_\sfrak^i}{}{2}.\]
    Multiplying the above relation by $\gamma_i\tau_i$, summing it from $i = 0$ to $t$ and noting the second relation in \eqref{eq:step_condn_1}, we get \eqref{eq:stoch_innerprod_delta_s}.
    

    Finally, we prove the bound on $\Ebb[\gnorm{\delta_\sfrak^i}{}{2}]$. 
    In view of SO relation in \eqref{eq:SO_g_operator-mo}, we have
    \begin{align}
        \Ebb[\gnorm{\ell_\gfrak^i(x^i)-\ell_g(x^i)}{}{2}]
        &\le 2\Ebb[\gnorm{\gfrak(x^{i-1}, \xibar^i) - g(x^{i-1})}{}{2}]
        + 2 \Ebb[\gnorm{[\Gfrak(x^{i-1}, \xibar^i) -  g'(x^{i-1})]^T(x^i-x^{i-1})}{}{2}] 
        \nonumber\\
        &\le 2\sigma_\gfrak^2 
        + 2\Ebb\big[
        \tsum_{j=1}^m \big( \Gfrak_j(x^{i-1}, \xibar^i) -  g'_j(x^{i-1})^T(x^i-x^{i-1})\big)^2
        \big]
        \nonumber\\
        &\le  2\sigma_\gfrak^2 
        + 2\Ebb\big[
        \tsum_{j=1}^m \gnorm{\Gfrak_j(x^{i-1}, \xibar^i) -  g'_j(x^{i-1})}{}{2} \gnorm{x^i-x^{i-1}}{}{2}
        \big]
        \nonumber\\
        &\le 2\sigma_\gfrak^2 + 2D_X^2\gnorm{\sigma_\Gfrak}{}{2} \label{eq:bd_var_ell_g1}
    \end{align}
    where the first inequality follows due to the relation $\gnorm{a+b}{}{2} \le 2\gnorm{a}{}{2} +2\gnorm{b}{}{2}$, and the last inequality follows by the definition of $\sigma_\Gfrak$.
    Note that the same arguments can be used to get 
    \begin{equation}
        \Ebb[\gnorm{\ell_\gfrak^i(x^{i-1})-\ell_g(x^{i-1})}{}{2}] \le 2\sigma_\gfrak^2 + 2D_X^2\gnorm{\sigma_\Gfrak}{}{2}\label{eq:bd_var_ell_g2}
    \end{equation}
    Now, we have 
    \begin{align*}
        \Ebb[\gnorm{\delta_\sfrak^i}{}{2}] 
        &= \Ebb\big[ \gnorm{(1+\theta_i)[\ell_\gfrak^i(x^i)-\ell_g(x^i)]
        -\theta_i[\ell_\gfrak^i(x^{i-1}) - \ell_g(x^{i-1})]}{}{2}\big]
        \nonumber\\
        &\le (1+\tfrac{\theta_i}{1+\theta_i})(1+\theta_i)^2 \Ebb[\gnorm{\ell_\gfrak^i(x^i)-\ell_g(x^i)}{}{2}] + (1+ \tfrac{1+\theta_i}{\theta_i})\theta_i^2\Ebb[\gnorm{\ell_\gfrak^i(x^{i-1}) - \ell_g(x^{i-1})}{}{2}]
        \nonumber\\
        &\le 2(1+2\theta_i)^2[\sigma_\gfrak^2 +D_X^2\gnorm{\sigma_\Gfrak}{}{2}]
    \end{align*}
    where in the first inequality, we used $\gnorm{a+b}{}{2} \le (1+k)\gnorm{a}{}{2} + (1+ \tfrac{1}{k})\gnorm{b}{}{2}$ and the second inequality follows from \eqref{eq:bd_var_ell_g1} and \eqref{eq:bd_var_ell_g2}.
    Hence, we conclude the proof.
\end{Proof}
Now, we show bounds related to the convergence of the \fstopconex~method.
\begin{lem}\label{lem:dstoch_final_bound}
    Suppose $\{\gamma_t, \theta_t, \eta_t, \tau_t\}$ satisfy \eqref{eq:step_condn_1}, \eqref{eq:step_condn_dstoch-1}, \eqref{eq:step_condn_dstoch-2} and \eqref{eq:step_condion_dstoc-3} are satisfied and let $B \ge 1$ be a constant. Then, we have
    \begin{align}
        \Gamma_{T}\Ebb[\sup_{x \in \wt{X}}\inprod{F(x)}{\wb{x}^{T}-x}] 
        &\le \tfrac{3\gamma_0\eta_0}{2} D_X^2
        + \tsum_{i= 0}^{T-1} \bracket[\Big]{ \tfrac{\gamma_i}{\tau_i} \Ebb[\gnorm{\delta_\sfrak^i}{}{2}] + 
        (\tfrac{4\gamma_i}{\eta_i-L_gB}+\tfrac{\gamma_i}{2\eta_i}) \Ebb[\gnorm{\delta_\Gfrak^i}{}{2}]
        + \tfrac{\gamma_i}{2\eta_i} \Ebb[\gnorm{\delta^{i+1}}{}{2}]
        }
        \nonumber\\
        +&\tsum_{i=1}^{T-1}\bracket[\Big]{ \tfrac{4\gamma_i\theta_i^2 H^2}{\eta_i - L_gB} + \tfrac{4\gamma_i\theta_i^2}{\eta_i - L_gB} \Ebb\gnorm{\delta^i-\delta^{i-1}}{}{2}
        } + \tfrac{2\gamma_{T-1}^2H^2}{\eta_{T-1}} + \tfrac{\gamma_{T-1}}{\eta_{T-1}}\Ebb\gnorm{\delta^T - \delta^{T-1}}{}{2},
        \label{eq:dstoch-bd-opt}\\
        \Gamma_T \Ebb[\gnorm{[g(\wb{x}^T]_+}{}{}] 
        &\le\tfrac{\gamma_0\eta_0}{2} \gnorm{x^*-x^0}{}{2} 
        + {\gamma_0\tau_0}(\gnorm{\lambda^*}{}{}+1)^2
        \nonumber\\
        +& \tsum_{i= 0}^{T-1}\bracket[\big]{
        \tfrac{3\gamma_i}{2\tau_i} \Ebb[\gnorm{\delta_\sfrak^i}{}{2}] 
        + \tfrac{4\gamma_i}{\eta_i-L_gB}\Ebb[\gnorm{\delta_\Gfrak^i}{}{2}]
        + \tfrac{4\gamma_i \mathcal{H}_*^2}{\eta_i - L_gB}}
        \nonumber\\
        +&\tsum_{i = 1}^{T-1}\bracket[\Big]{\tfrac{4\gamma_i\theta_i^2 H^2}{\eta_i - L_gB}
        + \tfrac{4\gamma_i\theta_i^2}{\eta_i - L_gB}\Ebb\gnorm{\delta^i-\delta^{i-1}}{}{2} 
        }+ \tfrac{2\gamma_{T-1}^2H^2}{\eta_{T-1}} + \tfrac{\gamma_{T-1}}{\eta_{T-1}}\Ebb\gnorm{\delta^T - \delta^{T-1}}{}{2},
        \label{eq:dstoch-bd-feas}
	\end{align}
 where recall that $\mathcal{H}_* :=  H_g(\gnorm{\lambda^*}{}{}+1) + \tfrac{L_gD_X}{2}[\gnorm{\lambda^*}{}{} + 1 -B]_+$.
\end{lem}
\begin{Proof}
    Note that all conditions in Lemma \ref{prop:dstoch_1st_result} are satisfied. Setting $\lambda = \lambda^0 = \zero$ in \eqref{eq:main_rel_dstoch} with index $t$ replaced by last iteration $T-1$, noting that $\inprod{\wb{\lambda}^T}{g(x)} \le 0$ for all $x \in \wt{X}$ and using \eqref{eq:stoch_innerprod_delta_G} and \eqref{eq:stoch_innerprod_delta}, we have for all $x \in \wt{X}$
       \begin{align*}
        \Gamma_{T}\inprod{F(x)}{\wb{x}^{T}-x}
        &\le\tfrac{3\gamma_0\eta_0}{2} \gnorm{x-x^0}{}{2} 
        + \tsum_{i=0}^{T-1}\Big[
        \gamma_i\inprod{\delta_\sfrak^i}{\lambda^i- \zero} 
        -\gamma_i\inprod{\delta_\Gfrak^i}{x^{i}-\wt{x}^i_a}
        -\gamma_i\inprod{\delta^{i+1}}{x^{i+1}-x^{i+1}_a}
        \Big]
        \nonumber\\
         &\quad+ \tsum_{i= 0}^{T-1} \bracket[\big]{
        \tfrac{\gamma_i}{\tau_i} \gnorm{\delta_\sfrak^i}{}{2} 
        +(\tfrac{4\gamma_i}{\eta_i-L_gB} + \tfrac{\gamma_i}{2\eta_i}) \gnorm{\delta_\Gfrak^i}{}{2}
        + \tfrac{\gamma_i}{2\eta_i} \gnorm{\delta^{i+1}}{}{2}
        }
        \nonumber\\
        &\quad+\tsum_{i=1}^{T-1}\bracket[\big]{
        \tfrac{4\gamma_i\theta_i^2 H^2}{\eta_i - L_gB}
        + \tfrac{4\gamma_i\theta_i^2}{\eta_i - L_gB}\gnorm{\delta^i-\delta^{i-1}}{}{2} 
        } + \tfrac{2\gamma_{T-1}^2H^2}{\eta_{T-1}} + \tfrac{\gamma_{T-1}}{\eta_{T-1}}\gnorm{\delta^T - \delta^{T-1}}{}{2},
    \end{align*}
    where we used $\mathcal{H}(\zero) = 0$. 
    Taking supremum over $x \in \wt{X}$ on both sides of the above relation, then taking expectation, noting the relations in \eqref{eq:useful-expectation2} with $\lambda = \zero$ and finally observing that $\Ebb[\inprod{\delta_\Gfrak^i}{\wt{x}^i_a - x^i}] = \Ebb[\inprod{\delta^{i+1}}{x^{i+1}_a - x^{i+1}}] = 0$, we obtain \eqref{eq:dstoch-bd-opt}. Hence, we get the bound on the expected optimality of $\wb{x}^T$.

    Now, we prove bound on the expected infeasibility of $\wb{x}^T$. Note that
	\begin{align}
		\inprod{F(x^*)}{\wb{x}^T-x^*} &+ (\gnorm{\lambda^*}{}{}+1)\gnorm{[g(\wb{x}^T)]_+}{}{} - \inprod{\wb{\lambda}^T}{g(x^*)} \nonumber\\
        &\ge \inprod{F(x^*)}{\wb{x}^T-x^*} + \inprod{\lambda^*}{g(\wb{x}^T)}   - \inprod{\wb{\lambda}^T}{g(x^*)}+ \gnorm{[g(\wb{x}^T)]_+}{}{}
		\ge \gnorm{[g(\wb{x}^T)]_+}{}{} \label{eq:int_rel8}
	\end{align}
    where the first inequality follows due to Cauchy-Schwarz inequality and  $[g(\wb{x}^T)]_+ \ge g(\wb{x}^T)$, and the second inequality follows from \eqref{eq:sp_condition_Lag_func} with $(x, \lambda) = (\wb{x}^T, \wb{\lambda}^T) \in X \times \Rbb^m_+$ 
    implying that 
    $\Lcal(\wb{x}^T, \lambda^*) \ge \Lcal(x^*, \wb{\lambda}^T)$. 
    Setting $x =x^*$ and $\lambda = \wt{\lambda}^T := (\gnorm{\lambda^*}{}{}+1)\tfrac{[g(\wb{x}^T)]_+}{\gnorm{[g(\wb{x}^T)]_+}{}{}}$ in \eqref{eq:main_rel_dstoch} and using \eqref{eq:int_rel8} to lower bound the resulting LHS and applying \eqref{eq:stoch_innerprod_delta_s}, we have 
    \begin{align*}\allowdisplaybreaks
        \Gamma_T \gnorm{[g(\wb{x}^T)]_+}{}{}
        &\le\tfrac{\gamma_0\eta_0}{2} \gnorm{x^*-x^0}{}{2} 
        + {\gamma_0\tau_0}(\gnorm{\lambda^*}{}{}+1)^2
        \nonumber\\
        &\quad 
        + \tsum_{i=0}^{T-1}\Big[
        \gamma_i\inprod{\delta_\sfrak^i}{\lambda^i- \lambda^i_a} 
        -\gamma_i\inprod{\delta_\Gfrak^i}{x^{i}-x^*}
        -\gamma_i\inprod{\delta^{i+1}}{x^{i+1}-x^*}
        \Big]
        \nonumber\\
         &\quad + \tsum_{i= 0}^{T-1} \bracket[\big]{
        \tfrac{3\gamma_i}{2\tau_i} \gnorm{\delta_\sfrak^i}{}{2} 
        + \tfrac{4\gamma_i}{\eta_i-L_gB}\gnorm{\delta_\Gfrak^i}{}{2}
        +\tfrac{4\gamma_i \mathcal{H}_*^2}{\eta_i - L_gB}
        }
        \nonumber\\
        &\quad+\tsum_{i = 1}^{T-1}\bracket[\Big]{ +\tfrac{4\gamma_i\theta_i^2 H^2}{\eta_i - L_gB}
        + \tfrac{4\gamma_i\theta_i^2}{\eta_i - L_gB}\gnorm{\delta^i-\delta^{i-1}}{}{2} 
        }+ \tfrac{2\gamma_{T-1}^2H^2}{\eta_{T-1}} + \tfrac{\gamma_{T-1}}{\eta_{T-1}}\gnorm{\delta^T - \delta^{T-1}}{}{2},
    \end{align*}
    where we note that $\mathcal{H}(\wt{\lambda}^T) = \mathcal{H}_*$. Now taking, expectation on both sides, noting \eqref{eq:useful-expectation2} with $x = x^*$ and observing that $\Ebb[\inprod{\delta_\sfrak^i}{\lambda^i- \lambda^i_a}] = 0$, we obtain \eqref{eq:dstoch-bd-feas}. Hence, we conclude the proof.
\end{Proof}
To use Lemma \ref{lem:dstoch_final_bound}, we need to use bounds on the stochastic noise terms, which are proved in Lemma \ref{lem:dstoch_useful_bounds}. 
However, as observed in \eqref{eq:bd_var_Gfrak}, the bound on $\Ebb[\gnorm{\delta_\Gfrak^i}{}{2}]$ depends on $\lambda^{i+1}$ which itself is an unbounded random variable. 
Hence, it is clear that to bound $\Ebb[\gnorm{\delta_\Gfrak^i}{}{2}]$, we must show that $\lambda^{i+1}$ generated by \fstopconex~method satisfies a bound on $\Ebb[\gnorm{\lambda^{i+1}}{}{2}]$. We prove such a result in the next lemma. 
\begin{lem}\label{lem:delta_Gfrak_bound}
    Suppose parameters $\{\gamma_t, \theta_t, \eta_t, \tau_t\}$ satisfy the requirements in Lemma \ref{prop:dstoch_1st_result}. Moreover, assume that $\{\eta_t, \tau_t\}$ satisfy 
    \begin{equation}\label{eq:step_condn_dstoch4}
        \tfrac{32\gnorm{\sigma_\Gfrak}{}{2}}{\tau_t(\eta_t-L_gB)} < 1,
    \end{equation}
    for all $t \le T-1$ and there exist two constants $U_1$ and $U_2$ satisfying
    \begin{align}
        U_1 &\ge 
        \big(1 - \tfrac{32\gnorm{\sigma_\Gfrak}{}{2}}{\tau_t(\eta_t-L_gB)}\big)^{-1}
        \Big[ 
        \tfrac{8\gnorm{\sigma_\Gfrak}{}{2}}{\gamma_t\tau_t}
        \Big\{
        \tfrac{\gamma_0\eta_0}{2} \gnorm{x^*-x^0}{}{2} 
        + \tsum_{i= 0}^t 
        \big[
        \tfrac{2\gamma_i(1+2\theta_i)^2}{\tau_i}(\sigma_\gfrak^2+ D_X^2\gnorm{\sigma_\Gfrak}{}{2})  
        + \tfrac{4\gamma_i \mathcal{H}_*^2}{\eta_i - L_gB} 
        \big]
        \nonumber\\
        &\hspace{13em}+ \tfrac{3\gamma_0\tau_0}{4}\gnorm{\lambda^*}{}{2}+\tsum_{i=1}^t \big[ \tfrac{4\gamma_i\theta_i^2 H^2}{\eta_i - L_gB}
        + \tfrac{8\gamma_i\theta_i^2\sigma^2 }{\eta_i - L_gB}
        \big]
        +\tfrac{2\gamma_tH^2}{\eta_t} + \tfrac{2\gamma_t\sigma^2}{\eta_t}
        \Big\}
        \Big]
        \label{eq:U1_relation}\\
        U_2 &\ge 
         \big(1 - \tfrac{32\gnorm{\sigma_\Gfrak}{}{2}}{\tau_t(\eta_t-L_gB)}\big)^{-1}\tfrac{32\gamma_i\gnorm{\sigma_\Gfrak}{}{2}}{\gamma_t\tau_t(\eta_i-L_gB)},
         \label{eq:U2_relation}
    \end{align}
    for all $t \le T-1$ and $i \le t-1$. 
    Then, we have
    \begin{equation}\label{eq:final_U1U2_bound_delta_G}
        \Ebb[\gnorm{\delta_\Gfrak^t}{}{2}] \le U_1(1+U_2)^t.
    \end{equation}
\end{lem}
\begin{Proof}
    As per Proposition \ref{prop:FCVI_SPP_reformulation}, we have $(x^*, \lambda^*)$ is saddle point of $\min_{x \in X}\max_{\lambda \ge \zero} \Lcal(x, \lambda)$. This implies
    \[
        \inprod{F(x^*)}{\wb{x}^{t+1}-x^*} + \inprod{\lambda^*}{g(\wb{x}^{t+1})} -\inprod{\wb{\lambda}^{t+1}}{g(x^*)} \ge 0.
    \]
    Using $x = x^*$ and $\lambda = \lambda^*$ in Lemma \ref{prop:dstoch_1st_result} (more specifically in relation \eqref{eq:main_rel_dstoch}),  replacing $\lambda^0 = \zero$, noting the above relation and the fact that $\mathcal{H}(\lambda^*) = H_g\gnorm{\lambda^*}{}{} + \tfrac{L_gD_X}{2}[\gnorm{\lambda^*}{}{} - B]_+\le \mathcal{H}_*$ (recall $\mathcal{H}_*$ from Lemma \ref{lem:dstoch_final_bound}), we get
    \begin{align}
        \tfrac{\gamma_{t}\tau_{t}}{4}\gnorm{\lambda^*-\lambda^{t+1}}{}{2}
        &\le\tfrac{\gamma_0\eta_0}{2} \gnorm{x^*-x^0}{}{2} 
        + \tfrac{\gamma_0\tau_0}{2}\gnorm{\lambda^*}{}{2}
        \nonumber\\
        &\quad 
        + \tsum_{i=0}^{t}\Big[
        \gamma_i\inprod{\delta_\sfrak^i}{\lambda^i- \lambda^*} 
        -\gamma_i\inprod{\delta_\Gfrak^i}{x^{i}-x^*}
        -\gamma_i\inprod{\delta^{i+1}}{x^{i+1}-x^*}
        \Big]
        \nonumber\\
         &\quad + \tsum_{i= 0}^t \bracket[\big]{
        \tfrac{\gamma_i}{\tau_i} \gnorm{\delta_\sfrak^i}{}{2} 
        + \tfrac{4\gamma_i}{\eta_i-L_gB}\gnorm{\delta_\Gfrak^i}{}{2}
        + \tfrac{4\gamma_i \mathcal{H}_*^2}{\eta_i - L_gB}}
        \nonumber\\
        &\quad +\tsum_{i=1}^t \bracket[\big]{\tfrac{4\gamma_i\theta_i^2 H^2}{\eta_i - L_gB}
        + \tfrac{4\gamma_i\theta_i^2}{\eta_i - L_gB}\gnorm{\delta^i-\delta^{i-1}}{}{2} 
        }  + \tfrac{2\gamma_tH^2}{\eta_t} + \tfrac{\gamma_t}{\eta_t}\gnorm{\delta^{t+1} - \delta^{t}}{}{2}.\label{eq:using_main_rel_dstoch}
    \end{align}
    In view of \eqref{eq:bd_var_Gfrak}, note that
    \begin{equation}\label{eq:delta_Gfrak_recursion_init}
        \Ebb[\gnorm{\delta_\Gfrak^t}{}{2}] 
        \le \gnorm{\sigma_\Gfrak}{}{2}\Ebb[\gnorm{\lambda^{t+1}}{}{2}] 
        \le 2\gnorm{\sigma_\Gfrak}{}{2}\bracket[\big]{\gnorm{\lambda^*}{}{2} + \Ebb[\gnorm{\lambda^{t+1} -\lambda^*}{}{2}]},
    \end{equation}
    where the second inequality follows by $\gnorm{a+b}{}{2} \le 2[\gnorm{a}{}{2} + \gnorm{b}{}{2}]$. Adding $\tfrac{\gamma_{t}\tau_{t}}{4}\gnorm{\lambda^*}{}{2}$ to \eqref{eq:using_main_rel_dstoch}, noting that second relation in \eqref{eq:step_condn_1} implies $\gamma_t\tau_t \le \gamma_0\tau_0$, multiplying the resulting relation by $\tfrac{8}{\gamma_t\tau_t}\gnorm{\sigma_\Gfrak}{}{2}$, then taking expectation on both sides,  employing \eqref{eq:bd_var_sfrak}, \eqref{eq:bd_var_F}, \eqref{eq:useful-expectation2} to substitute bounds for expectation terms on the RHS, and adding the resulting relation with \eqref{eq:delta_Gfrak_recursion_init}, we obtain
    \begin{align*}
        \Ebb[\gnorm{\delta_\Gfrak^t}{}{2}] 
        &\le \tfrac{8\gnorm{\sigma_\Gfrak}{}{2}}{\gamma_t\tau_t}\Big\{
        \tfrac{\gamma_0\eta_0}{2} \gnorm{x^*-x^0}{}{2} 
        + \tfrac{3\gamma_0\tau_0}{4}\gnorm{\lambda^*}{}{2} 
        + \tsum_{i= 0}^t \big[
        \tfrac{2\gamma_i(1+2\theta_i)^2}{\tau_i}(\sigma_\gfrak^2+ D_X^2\gnorm{\sigma_\Gfrak}{}{2}) 
        + \tfrac{4\gamma_i \mathcal{H}_*^2}{\eta_i - L_gB} 
        \big]
         \nonumber\\
         &\hspace{5em} 
         +\tsum_{i=1}^t\bracket[\big]{\tfrac{4\gamma_i\theta_i^2 H^2}{\eta_i - L_gB}
        + \tfrac{8\gamma_i\theta_i^2 \sigma^2}{\eta_i - L_gB} }
         +\tfrac{2\gamma_tH^2}{\eta_t} + \tfrac{2\gamma_t\sigma^2}{\eta_t}
        \Big\}
        + \tsum_{i=0}^t\tfrac{32\gamma_i\gnorm{\sigma_\Gfrak}{}{2}}{\gamma_t\tau_t(\eta_i-L_gB)}\Ebb[\gnorm{\delta_\Gfrak^i}{}{2}].
    \end{align*}
    Given \eqref{eq:step_condn_dstoch4}, the coefficient of $\delta_\Gfrak^t$ term in the right-hand side above is strictly less than 1. After shifting this term to the left-hand side, rearranging terms and noting the conditions \eqref{eq:U1_relation}, \eqref{eq:U2_relation} on $U_1, U_2$, respectively, we see the following recursive upper bound on $\Ebb[\gnorm{\delta_\Gfrak^t}{}{2}]$
    \[\Ebb[\gnorm{\delta_\Gfrak^t}{}{2}] \le U_1 + U_2 \tsum_{i=0}^{t-1} \Ebb[\gnorm{\delta_\Gfrak^i}{}{2}].\]
    Applying Proposition \ref{prop:recursive_bound-mo} on the above relation, we obtain \eqref{eq:final_U1U2_bound_delta_G}. Hence, we conclude the proof.
\end{Proof}
Note that even though we obtain a bound on $\Ebb[\gnorm{\delta_\Gfrak^t}{}{2}]$ in \eqref{eq:final_U1U2_bound_delta_G}, this bound is still a function of algorithmic parameters $\gamma_t, \theta_t, \eta_t, \tau_t$. Hence, to get the final result, we need to provide a step-size policy such that relations \eqref{eq:step_condn_1}-
\eqref{eq:step_condion_dstoc-3} and \eqref{eq:step_condn_dstoch4} are satisfied. Moreover, the step-size policy allows us to get useful bounds on $U_1$ and $U_2$ as required in \eqref{eq:U1_relation} and \eqref{eq:U2_relation}, respectively. Below, we show that stepsize policy of Theorem \ref{thm:feasib-fully-stoch} and prove the final convergence guarantee.\\
{\bf Proof of Theorem \ref{thm:full_stoch_opconex_conv}.}\\
    It is easy to see that \eqref{eq:step_policy_dstoch} satisfies \eqref{eq:step_condn_1}. Since $\eta \ge 8L$, we have the first relation in \eqref{eq:step_condn_dstoch-1}. Also, $\eta \ge \tfrac{8M_gB}{D_X}$ and $\tau \ge \tfrac{9M_gD_X}{B}$ implies the second relation in \eqref{eq:step_condn_dstoch-1}, \eqref{eq:step_condn_dstoch-2} and the first relation in \eqref{eq:step_condion_dstoc-3}. It is easy to see that $\eta \ge 8L + \tfrac{8M_gB}{D_X}$ and $\tau \ge \tfrac{9M_gD_X}{B}$ implies the second relation in \eqref{eq:step_condion_dstoc-3}. Finally, $\eta \ge \tfrac{8B\gnorm{\sigma_\Gfrak}{}{}}{D_X}$ and $\tau \ge \tfrac{9D_X\gnorm{\sigma_\Gfrak}{}{}}{B}$ implies $\tfrac{32\gnorm{\sigma_\Gfrak}{}{2}}{\tau_t(\eta_t-L_gB)} < \tfrac{4}{9} < 1$ which verifies \eqref{eq:step_condn_dstoch4}. 

    Now, we show bounds on $U_2$. As per \eqref{eq:U2_relation}, we see that $ \big(1 - \tfrac{32\gnorm{\sigma_\Gfrak}{}{2}}{\tau_t(\eta_t-L_gB)}\big)^{-1}\tfrac{32\gamma_i\gnorm{\sigma_\Gfrak}{}{2}}{\gamma_t\tau_t(\eta_i-L_gB)} < \tfrac{9}{5}.\tfrac{32\gnorm{\sigma_\Gfrak}{}{2}}{\eta\tau} < \tfrac{0.9\gnorm{\sigma_\Gfrak}{}{}D_X}{\sigma_{X,g}T} < \tfrac{0.9}{T}$ for all $t \le T-1$ and $i \le t-1$. Here, first inequality follows due to $\tfrac{32\gnorm{\sigma_\Gfrak}{}{2}}{\tau_t(\eta_t-L_gB)} < \tfrac{4}{9}$, second inequality follows due to $\eta \tau \ge \tfrac{64\gnorm{\sigma_\Gfrak}{}{}\sigma_{X,g}T}{D_X}$, and the last inequality follows by definition of $\sigma_{X,g}$. In view of the earlier relation and noting the RHS of \eqref{eq:U2_relation}, we can set $U_2 = \tfrac{0.9}{T}$.
    
    Noting the four relations 
    (i) $\tau^2 \ge \tfrac{64\sigma_{X,g}^2T}{B^2}$, 
    (ii) $\eta\tau \ge \tfrac{16\sigma_{X,g}(H + H_gB + \sqrt{2}\sigma + 4B\gnorm{\sigma_\Gfrak}{}{})T}{BD_X}$, 
    (iii)$\tfrac{\eta+L_gB}{\tau}D_X^2 
    \le \tfrac{(8L+L_gB)BD_X}{9M_g} + \tfrac{8}{9}B^2 + \tfrac{(H+H_gB+ \sqrt{2}\sigma + 4B\gnorm{\sigma_\Gfrak}{}{})BD_X}{4\sigma_{X,g}}$, and 
    (iv) $\tfrac{\gnorm{\sigma_\Gfrak}{}{}D_X}{\sigma_{X, g}} < 1$, we can see that RHS of \eqref{eq:U1_relation} for all $t \le T-1$ is at most 
    \begin{align*}
        &\tfrac{9}{5}\bracket[\big]{8\gnorm{\sigma_\Gfrak}{}{2}
        \braces[\big]{ \tfrac{L_gB + \eta}{\tau}D_X^2 
        + \tfrac{3}{4}\gnorm{\lambda^*}{}{2}
        +\tfrac{18\sigma_{X,g}^2T}{\tau^2}
        + \tfrac{4(H^2 + \mathcal{H}_*^2 + 2\sigma^2)T}{\eta\tau}
        }
        }\\
        &\le \tfrac{9}{5}\bracket[\big]{8\gnorm{\sigma_\Gfrak}{}{2}
        \braces[\big]{ \tfrac{L_gB + \eta}{\tau}D_X^2 
        + \tfrac{3}{4}\gnorm{\lambda^*}{}{2}
        +\tfrac{9}{32}B^2
        + \tfrac{4(H^2 + \mathcal{H}_*^2 + 2\sigma^2)T}{\eta\tau}
        }
        }\\
        &\le \tfrac{9}{5}\bracket[\big]{8\gnorm{\sigma_\Gfrak}{}{2}
        \braces[\big]{ \tfrac{L_gB + \eta}{\tau}D_X^2 
        + \tfrac{3}{4}\gnorm{\lambda^*}{}{2}
        +\tfrac{9}{32}B^2
        }
        + 2(H + \sqrt{2}\sigma + \tfrac{\mathcal{H}_*^2}{H+ H_gB + \sqrt{2}\sigma})B\gnorm{\sigma_\Gfrak}{}{}
        }\\
        &\le \tfrac{9}{5}\bracket[\big]{
        6\gnorm{\lambda^*}{}{2}\gnorm{\sigma_\Gfrak}{}{2}
        +18B^2\gnorm{\sigma_\Gfrak}{}{2}
        + \tfrac{(8L+L_gB)BD_X}{M_g}\gnorm{\sigma_\Gfrak}{}{2} 
        + 2(2H + 2\sqrt{2}\sigma + \tfrac{\mathcal{H}_*^2}{H+ H_gB + \sqrt{2}\sigma} + H_gB 
        )B\gnorm{\sigma_\Gfrak}{}{}
        }
    \end{align*}
    where first inequality follows by (i), second follows by (ii) and (iv), third follows by (iii) and (iv). Hence, we can set $U_1 := \tfrac{9}{5}\bracket[\big]{ 
    (6\gnorm{\lambda^*}{}{2} + 18B^2 + \tfrac{(8L+L_gB)BD_X}{M_g})\gnorm{\sigma_\Gfrak}{}{2} + 4(H + \sqrt{2}\sigma + \tfrac{\mathcal{H}_*^2}{2(H+ H_gB + \sqrt{2}\sigma)}  + \tfrac{H_gB}{2})B\gnorm{\sigma_\Gfrak}{}{}
    }$. Now applying Lemma \ref{lem:delta_Gfrak_bound} (particularly \eqref{eq:final_U1U2_bound_delta_G}), we have 
    \begin{equation}\label{eq:omega_bound_delta_G}
         \Ebb[\gnorm{\delta_\Gfrak^t}{}{2}] \le U_1(1+U_2)^{T-1} \le U_1(1+ \tfrac{0.9}{T})^{T-1} < e^{0.9}U_1 \le \omega^2,
    \end{equation}
    where the last relation follows from the definition of $\omega$. 

    Applying Lemma \ref{lem:dstoch_final_bound}, noting  \eqref{eq:bd_var_F}, \eqref{eq:bd_var_sfrak} and \eqref{eq:omega_bound_delta_G}, and substituting the step-size scheme in \eqref{eq:step_policy_dstoch}, we obtain \eqref{eq:conv_dstoch_opt} and \eqref{eq:conv_dstoch_feas}. Hence, we conclude the proof.

\section{Saddle point problems with coupling constraints}\label{sec:saddle}
In this section, we consider a general class of saddle point problems as follows
\begin{equation}\label{eq:FCSP}
	\begin{split}
		\min_{u \in U} \max_{v \in V}\ &f(u,v) \\
		\text{s.t.} \ &g(u,v) \le \zero,
	\end{split}
\end{equation} 
where $f(\cdot, v)$ is a convex function of $u$ for all $v \in V$, $f(u, \cdot)$ is a concave function of $v$ for all $u \in U$, $U  \subset \Rbb^{n_u}$ and $V  \subset \Rbb^{n_v}$ are feasible set constraints, and $g(u,v)$ is a coupling constraint that is jointly convex in $u$ and $v$. Notably, this convex-concave saddle point problem with a convex feasible domain has been a significant topic in the literature. From the complexity point of view, all algorithms in the literature consider only the case where the feasible domain of $(u,v)$ is a product set:  $U \times V$.  Recently, \cite{yang2024data} considered an advanced saddle point problem where they imposed convex function constraints $g_1(u) \le \zero$ and $g_2(u) \le \zero$. Notwithstanding the generality of this model compared to traditional set constraints, the feasible set $(U \cap \{g_1(u) \le \zero\}) \times (V \cap \{g_2(v) \le \zero\})$ remains a product set, which does not precisely capture the coupling constraint $g(u,v) \le\zero$. Conversely, the coupling constraint $g(u,v)$ can easily capture the function constraint setting by letting $g(u,v) = \begin{bmatrix}
	g_1^\top(u),
	g_2^\top(v)
\end{bmatrix}^\top$. Consequently, \eqref{eq:FCSP} seems to be the most general model of convex-concave saddle point problems in the literature. To the best of our knowledge, there are no algorithmic convergence complexity results for cases involving coupling constraints in the domain. In this section, we illustrate that the \opconex{}~method can be easily adapted to address this problem. As a result, we can solve \eqref{eq:FCSP} in various settings, including smooth/nonsmooth/stochastic $f$ and smooth/nonsmooth/stochastic $g$. Furthermore, by applying the same arguments used to extend the adaptive method from Section~\ref{sec:adaptive FCVI}, we can develop adaptive methods for \eqref{eq:FCSP} when $f$ and $g$ are Lipschitz smooth and deterministic functions. 

Throughout this section, we use $w:= (u,v)$ to denote the combined variable, $W:= U \times V$ for the product-form set constraint, and $g(w)$ to denote the coupling constraint $g(u,v)$. This notation carries in a consistent manner, e.g., $\wh{w} = (\wh{u}, \wh{v})$ and $w^t = (u^t, v^t)$. We define $\wt{W} := W \cap \{(u,v):g(u,v) \le \zero\}$. 

\paragraph{Convergence criterion.}
We start with a brief discussion on the convergence criterion for \eqref{eq:FCSP}. Note that a standard approach of constructing the Lagrangian can be deceiving here. Indeed, the usual ``Lagrangian'' $f(u,v) + \lambda g(u,v)$ is neither convex nor concave in $v$. However, this approach is flawed because it transforms a convex-concave SP problem into a convex-nonconcave problem, causing the loss of tractability. To rectify this notion, we first look at the convergence criterion for set-constrained problems (where $g(w)\le \zero$ is dropped): $\wh{w} = (\wh{u}, \wh{v}) \in U \times V$ is an $\epsilon$-approximate SP solution if 
\[ \max_{w \in W}\ G(\wh{w} ; w) [:= f(\wh{u}, v) -f(u, \wh{v})] \le \epsilon.\]
It is well-known that for any feasible $\wh{w}$, the primal-dual gap $\max_{w \in W}\ G(\wh{w}; w)$ is always nonnegative and is 0 iff $\wh{w}$ is a saddle point. Since $G(\wh{w}; w)$ is convex with respect to $\wh{w}$ for all $w \in W$, hence the primal-dual gap is convex with respect to the variable $\wh{w}$. Essentially, even if the problem is concave in variable $v$, the convergence criterion is convex with respect to candidate solution $\wh{v}$. We will exploit this fact in the convergence analysis. A natural extension of the above criterion for \eqref{eq:FCSP} is as follows:\vspace{-0.5em}
\begin{Defi}\label{def:appx_soln_FCSP}
	A point $\wh{w} \in W$ is an $\epsilon$-approximate saddle point of \eqref{eq:FCSP} if it satisfies
	\begin{equation*}
		\max_{w \in \wt{W}}\ G(\wh{w}; w) \le \epsilon,\qquad
		g(\wh{w}) \le \epsilon.
	\end{equation*}
\end{Defi}
\vspace{-5pt}Without coupling constraint $g$, the above criterion reduces to standard criterion, and when $\epsilon \to 0$, it exactly characterizes saddle points. This is equivalent to Nash equilibrium. Indeed, $\epsilon = 0$ implies (i) feasibility $g(\wh{w}) \le \zero, \wh{w} \in W$, and (ii) primal-dual gap is equal to $0$. Putting these two facts together implies $\wh{w}$ is the saddle point solution of \eqref{eq:FCSP}. 
With the coupling constraint $g$, the above criterion is more strict than the generalized Nash equilibrium, which incorporates shared constraints between the players in the Nash equilibrium (NE). Recall the definition that for the two-person game,  $w^\star=(u^\star,v^\star)$ is a generalized Nash equilibrium (GNE) of \eqref{eq:FCSP} if 
\[ f(u^\star, v)\le f(u^\star, v^\star) \le f(u, v^\star), \quad \forall u\in U\cap \{ g(u,v^\star) \le \zero\},\ \forall v\in V\cap \{ g(u^\star, v)\le \zero\}.
\]
Suppose that $(\wh{u},\wh{v})$ is a saddle point, (i.e., $\epsilon=0$). By definition we have  $(\wh{u},\wh{v})\in \wt{W}$ and 
\begin{equation}\label{eq:sd-defi}
\max_{(u,v)\in W} \quad   f(\wh{u}, v)- f(u,\wh{v}) \le 0, \qquad 
\st   \quad g(u,v) \le \zero.
\end{equation}
Then, for any $u\in U\cap \{g(u,\wh{v})\le \zero\}$, due to the feasibility, plugging $(u,v)\leftarrow (u, \wh{v})$ in \eqref{eq:sd-defi}, we have $f(\wh{u},\wh{v})-f(u, \wh{v})\le 0$. Similarly, plugging $ (\wh{u}, {v})$ in \eqref{eq:sd-defi} for any $v\in V\cap \{g(\wh{u},{v})\le \zero\}$, we have $f(\wh{u},v) \le f(\wh{u},\wh{v})$. This result shows that $(\wh{u},\wh{v})$ is indeed the optimal strategy when the opponent is fixed, and hence is a GNE solution.

Now, we are ready to present the algorithmic details. For brevity, we only present the application of \opconex{} for the fully stochastic version \eqref{eq:FCSP}. It is worthwhile to mention that \adopconex{} and \sopconex{} methods apply to problem~\eqref{eq:FCSP} with the same convergence rate guarantee in terms of criterion in Definition \ref{def:appx_soln_FCSP} in the deterministic smooth, fully-stochastic nonsmooth settings, respectively. Moreover, our \opconex{} method improves upon \sopconex{} method in the criterion in Definition \ref{def:appx_soln_FCSP}. Hence, our methods provide the most competitive convergence guarantees for convex-concave saddle point problems in the most comprehensive setup. 

\subsection{\opconex~method for the saddle point problem}
To utilize the \opconex~method for solving \eqref{eq:FCSP}, we reformulate \eqref{eq:FCSP} as an FCVI problem \eqref{eq:FCVI} as follows. 
Let the set constraint $X = W$ and the variable $x = w$. We define $F(x) = F(w) := \begin{bmatrix}
	f'_u(u,v)\\ - f'_v(u,v)
\end{bmatrix}$ as the operator throughout this section. Here, $f'_u(u,v) \in \partial f(\cdot, v)(u)$ is a subgradient of the function $f(\cdot,v)$ at $u$. Similarly, $- f'_v(u,v)$ is the subgradient of the function $-f(u, \cdot)$ at $v$. It is clear that since $f$ is a convex-concave function, the resulting operator $F(w)$ is monotone. As per this reformulation, convex coupled constraint $\{g(u,v) \le \zero\}$ in \eqref{eq:FCSP} is equivalent to $g(x) \le \zero$ in \eqref{eq:FCVI}. 
Then, it is clear that $\wt{W}$ corresponds to $\wt{X}$. 

Now we discuss the correspondence of fully-stochastic FCVI with fully-stochastic \eqref{eq:FCSP}. We denote the stochastic operator 
$$\Ffrak^t = \Ffrak(w^t, \xi^t) = \begin{bmatrix}
    f'_u(u^t,v^t, \xi)\\ - f'_v(u^t,v^t, \xi) 
\end{bmatrix} = \begin{bmatrix}
    \Ffrak_u(w^t, \xi^t) \\ \Ffrak_v(w^t, \xi^t)
\end{bmatrix} = \begin{bmatrix}
    \Ffrak_u^t \\ \Ffrak_v^t
\end{bmatrix}.$$ 
Moreover, we denote 
\[\Gfrak^t_i = \Gfrak_i(w^t, \xi^t) = \begin{bmatrix}
    g'_u(u^t,v^t, \xi)\\ g'_v(u^t,v^t, \xi) 
\end{bmatrix} = \begin{bmatrix}
    \Gfrak_{u,i}(w^t, \xi^t) \\ \Gfrak_{v,i}(w^t, \xi^t)
\end{bmatrix} = \begin{bmatrix}
    \Gfrak^t_{u,i} \\ \Gfrak^t_{v,i}
\end{bmatrix}\]
Hence, we can apply the \opconex~method (Algorithm \ref{alg:alg4}) for problem \eqref{eq:FCSP} as follows:
In line 4 of \opconex~method, $s^t$ will be evaluated as 
$(1+ \theta_t)\ell_\gfrak^t(w^t) - \theta_t\ell_\gfrak^t(w^{t-1})$ in the similar manner.  Moreover, line 6 
of Algorithm \ref{alg:alg4} decomposes into two separate objective  for variables $u$ and $v$:
\begin{align*}
	&\inprod*{(1+\theta_t)\begin{bmatrix}
			\Ffrak^t_u\\\Ffrak^t_v
		\end{bmatrix} -\theta_t \begin{bmatrix}
			\Ffrak^{t-1}_u\\\Ffrak^{t-1}_v
	\end{bmatrix} + \tsum_{i \in [m] } \lambda^{t+1}_i \begin{bmatrix}
	\Gfrak^t_{u,i} \\ \Gfrak^t_{v,i}
\end{bmatrix}}{\begin{bmatrix}
u\\v
\end{bmatrix}} + \tfrac{\eta_t}{2} \gnorm*{\begin{bmatrix}
u\\v
\end{bmatrix} - \begin{bmatrix}
u^t\\v^t
\end{bmatrix}}{}{2}\\
&= \inprod{(1+\theta_t)\Ffrak^t_u -\theta_t \Ffrak^{t-1}_u +\tsum_{i \in [m] } \lambda^{t+1}_i \Gfrak^t_{u,i} } {u} + \tfrac{\eta_t}{2} \gnorm{u-u^t}{}{2}\\
&\quad+\inprod{(1+\theta_t)\Ffrak^t_v -\theta_t \Ffrak^{t-1}_v +\tsum_{i \in [m] } \lambda^{t+1}_i \Gfrak^t_{v,i} } {v} + \tfrac{\eta_t}{2} \gnorm{v-v^t}{}{2}.
\end{align*}
The above decomposition, along with the fact that set $W = U \times V$ is a product set, implies that the optimization problem can be decomposed itself. This also implies that $u^{t+1}$ and $v^{t+1}$ of this algorithm can be computed in parallel. Then, the update scheme of \opconex{} method for \eqref{eq:FCSP} looks as follows:
\begin{align*}
    \sfrak^t &\leftarrow (1+\theta_t)\ell_\gfrak^t(x^t) - \theta_t\ell_\gfrak^t(x^{t-1})\\
    \lambda^{t+1} &\leftarrow \argmin_{\lambda \ge \zero} \inprod{-\sfrak^t}{\lambda} + \tfrac{\tau_t}{2}\gnorm{\lambda-\lambda^t}{}{2}\\
    u^{t+1} &\leftarrow \argmin_{u \in U} \inprod{(1+\theta_t)\Ffrak^t_u -\theta_t \Ffrak^{t-1}_u +\tsum_{i \in [m] } \lambda^{t+1}_i \Gfrak^t_{u,i} } {u} + \tfrac{\eta_t}{2} \gnorm{u-u^t}{}{2}\\
    v^{t+1} &\leftarrow \argmin_{v\in V} \inprod{(1+\theta_t)\Ffrak^t_v -\theta_t \Ffrak^{t-1}_v +\tsum_{i \in [m] } \lambda^{t+1}_i \Gfrak^t_{v,i} } {v} + \tfrac{\eta_t}{2} \gnorm{v-v^t}{}{2}
\end{align*}Now, we provide a convergence analysis of this algorithm.

\paragraph{Convergence analysis} It is a well-known $w^* = (u^*,v^*)$ is a saddle point of \eqref{eq:FCSP} if and only if $w^*$ is the solution for corresponding VI reformulation \eqref{eq:FCVI}. 

\begin{lem}\label{lem:basic_for_SP}
	Under the assumptions of Lemma \ref{prop:dstoch_1st_result}, we get similar bounds on the criterion
	\begin{align*}
		\Gamma_{t+1}\big[G(\wb{w}^{t+1}, w) + \inprod{\lambda}{g(\wb{w}^{t+1})} - \inprod{\wb{\lambda}^{t+1}}{g(w)}\big].
	\end{align*}
\end{lem}
\begin{Proof}
	Our convergence analysis follows the proof of Lemma \ref{prop:dstoch_1st_result}.
	We modify \eqref{eq:int_rel4} as follows:
	\begin{equation}\label{eq:SP_proof_relation1}
		\inprod{F(w^{t+1})}{w^{t+1}-w} \ge G(w^{t+1}, w) [=f(u^{t+1},v) -f(u, v^{t+1})], \qquad \forall \ w\in W.
	\end{equation}
	This follows due to the new definition of $F$ and the convex-concave property of $f$. Replacing this bound in the analysis and the following proof of Lemma \ref{prop:dstoch_1st_result} up to  \eqref{eq:before_sum}, we obtain bound on the 
    \[G(w^{i+1},w) + \inprod{\lambda^{i+1}}{g(w^{i+1})} - \inprod{\lambda^{i+1}}{g(w)}.\]
    The rest of the proof follows by noting the fact that 
    \[\tsum_{i = 0}^t\gamma_i G(w^{i+1},w) \ge \Gamma_{t+1}G(\wb{w}^{t+1},w),\]
    where the inequality follows from convexity of $G(\cdot, w)$.
\end{Proof}
\begin{lem}\label{lem:adavance_for_FCSP}
    Under the assumptions of Lemma \ref{lem:dstoch_final_bound}, we get similar bounds on $\Gamma_T \Ebb[\sup_{w \in \wt{W}} G(\wb{w}^T,w)]$ and $\Gamma_T \Ebb[[g(\wb{w}^T)]_+]$.
\end{lem}
\begin{Proof}
    The proof follows by modifying \eqref{eq:int_rel8}. We have 
    \begin{align}
		G\paran{\wb{w}^T,w^*} &+ (\gnorm{\lambda^*}{}{}+1)\gnorm{[g(\wb{w}^T)]_+}{}{} - \inprod{\wb{\lambda}^T}{g(x^*)} \nonumber\\
        &\ge G\paran{\wb{w}^T,w^*} + \inprod{\lambda^*}{g(\wb{w}^T)}   - \inprod{\wb{\lambda}^T}{g(w^*)}+ \gnorm{[g(\wb{w}^T)]_+}{}{} \nonumber\\
        &\ge \inprod{F(w^*)}{\wb{w}^T-w^*} + \inprod{\lambda^*}{g(\wb{w}^T)} - \inprod{\wb{\lambda}^T}{g(w^*)} + \gnorm{[g(\wb{w}^T)]_+}{}{}] \nonumber\\
		&\ge \gnorm{[g(\wb{x}^T)]_+}{}{}, \label{eq:int_rel8_SP}
	\end{align}
    where the first inequality follows from the Cauchy-Schwarz inequality and the fact that $\inprod{\lambda^*}{[g(\wb{w}^T)]_+ }\ge \inprod{\lambda^*}{g(\wb{w}^T)}$, second inequality follows due to the new definition of $F$ and convex-concave nature of $f(u,v)$, and third inequality follows since $(w^*, \lambda^*)$ is a saddle point of $\Lcal(w, \lambda)$ defined on the new monotone operator $F$ associated with \eqref{eq:FCSP}. Then, using $w = w^*$, $\lambda = \wt{\lambda}^T$ and following the proof of Lemma \ref{lem:dstoch_final_bound}, we conclude this proof.
\end{Proof}
Using the above two lemmas and following the proof of Theorem \ref{thm:full_stoch_opconex_conv}, we have a similar convergence guarantee of \opconex{} method for \eqref{eq:FCSP}.
\begin{thm}
    Under the conditions of Theorem \ref{thm:full_stoch_opconex_conv}, we have the same convergence guarantees on the convergence criteria
    $\Ebb[\sup_{w \in \wt{W}} G(\wb{w}^T,w)]$ and $\Ebb[[g(\wb{w}^T)]_+]$.
\end{thm}
Given the similarity of the analysis for the FCVI problem \eqref{eq:FCVI} and saddle point problem \eqref{eq:FCSP}, we can easily transpose convergence results of Theorem \ref{thm:feasib-fully-stoch} and Corollary \ref{cor:conv_adaptive_OE} for \sopconex{} and \adopconex{} method, respectively, for \eqref{eq:FCSP} in terms of convergence criterion of Definition \ref{def:appx_soln_FCSP}. We skip the details for the sake of brevity.

\section{Conclusion}\label{sec:conclusion}
This paper focused on monotone function constrained variational inequality (FCVI) problems. We showed an equivalent reformulation of this problem to a VI problem on the KKT operator. Since the KKT operator is not uniformly Lipschitz continuous, we showed an adaptive version of the Operator Extrapolation method yields the optimal $\mathcal{O}(\tfrac{1}{T})$ convergence rate for smooth deterministic FCVIs. When the problem is nonsmooth and fully-stochastic, we showed that a modified version of the operator extrapolation method, called \sopconex{}, converges to the solution at the rate of $\mathcal{O}(\tfrac{1}{\sqrt{T}})$. This method has suboptimal dependence on some noise and Lipschtiz constant terms. To improve on this convergence behavior, we propose a fast \opconex{} method that overcomes the bottlenecks of \sopconex{} method through an innovative combination of the operator and constraint extrapolation terms and incorporating a stable momentum update. All of these convergence results are new in the literature. 
Finally, we show that all of our methods can be applied to convex-concave saddle point problems with coupled function constraints while maintaining the same complexity results for all of our methods. This shows the first such comprehensive complexity result for the most general formulation of saddle point problems. 

\bibliographystyle{abbrvnat}
\bibliography{ref.bib}

\end{document}